# WELL-POSEDNESS OF THE PLASMA-VACUUM INTERFACE PROBLEM

Paolo Secchi

Dipartimento di Matematica, Facoltà di Ingegneria, Università di Brescia
Via Valotti, 9, 25133 Brescia, Italy

Yuri Trakhinin

Sobolev Institute of Mathematics, Koptyug av. 4, 630090 Novosibirsk, Russia

Abstract. We consider the free boundary problem for the plasma-vacuum interface in ideal compressible magnetohydrodynamics (MHD). In the plasma region the flow is governed by the usual compressible MHD equations, while in the vacuum region we consider the *pre-Maxwell dynamics* for the magnetic field. At the free-interface, driven by the plasma velocity, the total pressure is continuous and the magnetic field on both sides is tangent to the boundary. The plasma-vacuum system is not isolated from the outside world, because of a given surface current on the fixed boundary that forces oscillations.

Under a suitable stability condition satisfied at each point of the initial interface, stating that the magnetic fields on either side of the interface are not collinear, we show the existence and uniqueness of the solution to the nonlinear plasma-vacuum interface problem in suitable anisotropic Sobolev spaces. The proof is based on the results proved in the companion paper [35], about the well-posedness of the homogeneous linearized problem and the proof of a basic a priori energy estimate. The proof of the resolution of the nonlinear problem given in the present paper follows from the analysis of the elliptic system for the vacuum magnetic field, a suitable tame estimate in Sobolev spaces for the full linearized equations, and a Nash-Moser iteration.

## 1. Introduction

Consider the equations of ideal compressible MHD:

$$\begin{cases} \partial_t \rho + \operatorname{div}(\rho v) = 0, \\ \partial_t(\rho v) + \operatorname{div}(\rho v \otimes v - H \otimes H) + \nabla q = 0, \\ \partial_t H - \nabla \times (v \times H) = 0, \\ \partial_t(\rho e + \frac{1}{2}|H|^2) + \operatorname{div}\big((\rho e + p)v + H \times (v \times H)\big) = 0, \end{cases} \quad (1)$$

where $\rho$ denotes density, $v \in \mathbb{R}^3$ plasma velocity, $H \in \mathbb{R}^3$ magnetic field, $p = p(\rho, S)$ pressure, $q = p + \frac{1}{2}|H|^2$ total pressure, $S$ entropy, $e = E + \frac{1}{2}|v|^2$ total energy, and $E = E(\rho, S)$ internal energy. With a state equation of gas, $\rho = \rho(p, S)$, and the first principle of thermodynamics, (1) is a closed system.

System (1) is supplemented by the divergence constraint

$$\operatorname{div} H = 0 \quad (2)$$

Date: January 22, 2013.

The first author PS is supported by the national research project PRIN 2009 "Equations of Fluid Dynamics of Hyperbolic Type and Conservation Laws". Part of this work was done during the fellowship of the second author YT at the Landau Network-Centro Volta-Cariplo Foundation spent at the Department of Mathematics of the University of Brescia in Italy. YT would like to warmly thank the Department of Mathematics of the University of Brescia for its kind hospitality during the visiting period.

on the initial data. As is known, taking into account (2), we can easily symmetrize system (1) by rewriting it in the nonconservative form

$$
\begin{cases}
\dfrac{\rho_p}{\rho}\dfrac{\mathrm{d}p}{\mathrm{d}t} + \operatorname{div} v = 0, \qquad \rho\,\dfrac{\mathrm{d}v}{\mathrm{d}t} - (H\cdot\nabla)H + \nabla q = 0, \\[2mm]
\dfrac{\mathrm{d}H}{\mathrm{d}t} - (H\cdot\nabla)v + H\operatorname{div} v = 0, \qquad \dfrac{\mathrm{d}S}{\mathrm{d}t} = 0,
\end{cases}
\tag{3}
$$

where $\rho_p \equiv \partial\rho/\partial p$ and $\mathrm{d}/\mathrm{d}t = \partial_t + (v\cdot\nabla)$. A different symmetrization is obtained if we consider $q$ instead of $p$. In terms of $q$ the equation for the pressure in (3) takes the form

$$
\frac{\rho_p}{\rho}\left\{\frac{\mathrm{d}q}{\mathrm{d}t} - H\cdot\frac{\mathrm{d}H}{\mathrm{d}t}\right\} + \operatorname{div} v = 0,
\tag{4}
$$

where it is understood that now $\rho = \rho(q - |H|^2/2, S)$ and similarly for $\rho_p$. Then we derive $\operatorname{div} v$ from (4) and rewrite the equation for the magnetic field in (3) as

$$
\frac{\mathrm{d}H}{\mathrm{d}t} - (H\cdot\nabla)v - \frac{\rho_p}{\rho}H\left\{\frac{\mathrm{d}q}{\mathrm{d}t} - H\cdot\frac{\mathrm{d}H}{\mathrm{d}t}\right\} = 0.
\tag{5}
$$

Substituting (4), (5) in (3) then gives the following symmetric system

$$
\begin{pmatrix}
\rho_p/\rho & \underline{0} & -(\rho_p/\rho)H & 0 \\
\underline{0}^T & \rho I_3 & 0_3 & \underline{0}^T \\
-(\rho_p/\rho)H^T & 0_3 & I_3 + (\rho_p/\rho)H\otimes H & \underline{0}^T \\
0 & \underline{0} & \underline{0} & 1
\end{pmatrix}
\partial_t
\begin{pmatrix} q \\ v \\ H \\ S \end{pmatrix}
+
$$

$$
+
\begin{pmatrix}
(\rho_p/\rho)v\cdot\nabla & \nabla\cdot & -(\rho_p/\rho)Hv\cdot\nabla & 0 \\
\nabla & \rho v\cdot\nabla I_3 & -H\cdot\nabla I_3 & \underline{0}^T \\
-(\rho_p/\rho)H^Tv\cdot\nabla & -H\cdot\nabla I_3 & (I_3 + (\rho_p/\rho)H\otimes H)v\cdot\nabla & \underline{0}^T \\
0 & \underline{0} & \underline{0} & v\cdot\nabla
\end{pmatrix}
\begin{pmatrix} q \\ v \\ H \\ S \end{pmatrix}
= 0,
\tag{6}
$$

where $\underline{0} = (0,0,0)$. Given this symmetrization, as the unknown we can choose the vector $U = U(t,x) = (q,v,H,S)$. For the sake of brevity we write system (6) in the form

$$
A_0(U)\partial_t U + \sum_{j=1}^{3} A_j(U)\partial_j U = 0,
\tag{7}
$$

which is symmetric hyperbolic provided the hyperbolicity condition $A_0 > 0$ holds:

$$
\rho > 0, \quad \rho_p > 0.
\tag{8}
$$

Plasma-vacuum interface problems for system (1) appear in the mathematical modeling of plasma confinement by magnetic fields (see, e.g., [14]). In this model the plasma is confined inside a perfectly conducting rigid wall and separated from it by a vacuum region, due to the effect of strong magnetic fields. However, the plasma-vacuum system is not isolated from the outside world because energy flows into the system. This can be modeled by a given surface current which forces oscillations onto the system.

This subject is very popular since the 1950–70's, but most of theoretical studies are devoted to finding stability criteria of equilibrium states. The typical work in this direction is the classical paper of Bernstein et al. [5]. In astrophysics, the plasma-vacuum interface problem can be used for modeling the motion of a star or the solar corona when magnetic fields are taken into account.

According to our knowledge there are still no well-posedness results for full (*non-stationary*) plasma-vacuum models. More precisely, a basic energy a priori estimate in Sobolev spaces for the linearization of a plasma-vacuum interface problem (see its description just below) was derived in [38], and the existence of solutions to this problem was recently proved in [35]. The proof of the existence and uniqueness of smooth solutions of the original nonlinear free boundary problem is the main goal of the present paper. Note that the a priori estimate for the linearized problem obtained in [35] somewhat improves the similar estimate firstly deduced in [38]. We also use the same notations and functional spaces as in [35].

Let $\Omega^+(t)$ and $\Omega^-(t)$ be space-time domains occupied by the plasma and the vacuum respectively. That is, in the domain $\Omega^+(t)$ we consider system (1) (or (7)) governing the motion of an ideal plasma and in the domain $\Omega^-(t)$, as in [5, 14], we consider the so-called *pre-Maxwell dynamics*

$$\nabla \times \mathcal{H} = 0, \qquad \operatorname{div} \mathcal{H} = 0, \tag{9}$$

$$\nabla \times E = -\partial_t \mathcal{H}, \qquad \operatorname{div} E = 0, \tag{10}$$

describing the vacuum magnetic field $\mathcal{H} \in \mathbb{R}^3$ and electric field $E \in \mathbb{R}^3$. That is, as usual in nonrelativistic MHD, in the Maxwell equations we neglect the displacement current $(1/c)\,\partial_t E$, where $c$ is the speed of light.

From (10) the electric field $E$ is a secondary variable that may be computed from the magnetic field $\mathcal{H}$. Hence, in the vacuum only one basic variable is needed, viz. $\mathcal{H}$, satisfying the elliptic (div-curl) system (9).

Let us assume that the interface between plasma and vacuum is given by a hypersurface $\Gamma(t) = \{F(t,x) = 0\}$. It is to be determined and moves with the velocity of plasma particles at the boundary:

$$\frac{\mathrm{d}F}{\mathrm{d}t} = 0 \quad \text{on } \Gamma(t) \tag{11}$$

(for all $t \in [0,T]$). As $F$ is an unknown of the problem, this is a free-boundary problem. The plasma variable $U$ is connected with the vacuum magnetic field $\mathcal{H}$ through the relations [5, 14]

$$[q] = 0, \quad H \cdot N = 0, \quad \mathcal{H} \cdot N = 0, \quad \text{on } \Gamma(t), \tag{12}$$

where $N = \nabla F$ and $[q] = q|_\Gamma - \frac{1}{2}|\mathcal{H}|^2_{|\Gamma}$ denotes the jump of the total pressure across the interface. These relations together with (11) are the boundary conditions at the interface $\Gamma(t)$.

As in [19, 37], we will assume that for problem (1), (9), (11), (12) the hyperbolicity conditions (8) are assumed to be satisfied in $\Omega^+(t)$ up to the boundary $\Gamma(t)$, i.e., the plasma density does not go to zero continuously, but has a jump (clearly in the vacuum region $\Omega^-(t)$ the density is identically zero). This assumption is compatible with the continuity of the total pressure in (12).

Since the interface moves with the velocity of plasma particles at the boundary, by introducing the Lagrangian coordinates one can reduce the original problem to that in a fixed domain. This approach has been recently employed with success in a series of papers on the Euler equations in vacuum, see [9, 10, 11, 12, 19]. However, as, for example, for contact discontinuities in various models of fluid dynamics (e.g., for current-vortex sheets [7, 36]), this approach seems hardly applicable for problem (1), (9), (11), (12). Therefore, we will work in the Eulerian coordinates and for technical simplicity we will assume that the space-time domains $\Omega^\pm(t)$ have the following form.

1.1. **The reference domain $\Omega$.** To avoid using local coordinate charts necessary for arbitrary geometries, and for the sake of simplicity, we will assume that the space domain $\Omega$ occupied by plasma and vacuum is given by

$$\Omega := \{(x_1, x_2, x_3) \in \mathbb{R}^3 \mid , x_1 \in (-1,1)\, x' = (x_2, x_3) \in \mathbb{T}^2\},$$

where $\mathbb{T}^2$ denotes the 2-torus, which can be thought of as the unit square with periodic boundary conditions. This permits the use of *one* global Cartesian coordinates system. We also set

$$\Omega^\pm := \Omega \cap \{x_1 \gtrless 0\}, \qquad \Gamma := \Omega \cap \{x_1 = 0\}.$$

Let us assume that the moving interface $\Gamma(t)$ takes the form

$$\Gamma(t) := \{(x_1, x') \in \mathbb{R} \times \mathbb{T}^2,\, x_1 = \varphi(t,x')\} \qquad t \in [0,T],$$

where it is assumed that $-1 < \varphi(t,\cdot) < 1$. Then we have $\Omega^\pm(t) = \{x_1 \gtrless \varphi(t,x')\} \cap \Omega$. With our parametrization of $\Gamma(t)$, an equivalent formulation of the boundary conditions (11), (12) at the free interface is

$$\partial_t \varphi = v_N, \quad [q] = 0, \quad H_N = 0, \quad \mathcal{H}_N = 0 \quad \text{on } \Gamma(t), \tag{13}$$

where $v_N = v \cdot N,\, H_N = H \cdot N,\, \mathcal{H}_N = \mathcal{H} \cdot N,\, N = (1, -\partial_2 \varphi, -\partial_3 \varphi)$.

On the fixed *top* and *bottom* boundaries

$$\Gamma_\pm := \{(\pm 1, x'),\, x' \in \mathbb{T}^2\}$$

of the domain $\Omega$, we prescribe the boundary conditions

$$v_1 = H_1 = 0 \quad \text{on } [0, T] \times \Gamma_+\,, \qquad \nu \times \mathcal{H} = \mathfrak{J} \quad \text{on } [0, T] \times \Gamma_-\,. \tag{14}$$

In the last equation $\nu = (-1, 0, 0)$ is the outward normal vector at $\Gamma_-$ and $\mathfrak{J}$ represents a given surface current which forces oscillations onto the plasma-vacuum system. The effect of such an outer boundary is that the system is not isolated from the outside world because energy flows into the system. In laboratory plasmas this external excitation may be caused by a system of coils. The model can also be exploited for the analysis of waves in astrophysical plasmas, e.g. by mimicking the effects of excitation of MHD waves by an external plasma by means of a localized set of "coils", when the response of the internal plasma is the main issue (e.g. in the problem of sunspot oscillations excited by sound waves in the photosphere). For a more complete discussion we refer the reader to [14].

When the system is isolated from the outside world, the natural boundary condition at $\Gamma_-$ for the vacuum magnetic field is $\mathcal{H} \cdot \nu = \mathcal{H}_1 = 0$, for perfectly conducting wall, i.e. the same we are prescribing at $\Gamma_+$ for $H$. For a simply connected domain as in the above choice, $\mathcal{H}$ is then necessarily zero ($\mathcal{H}$ is the unique solution of a homogeneous problem) and one solves the plasma equations with a vanishing total pressure $q$ on $\Gamma(t)$. The problem becomes meaningful for non simply connected vacuum regions, as in most of interesting applications, see [14]. In that case the null space associated to the homogeneous equations (9), under the boundary conditions $\mathcal{H}_N = 0$ on $\Gamma(t)$, $\mathcal{H}_1 = 0$ on $\Gamma_-$, is finite-dimensional, see [3]. One looks for a vacuum magnetic field from this finite-dimensional subspace, interacting with the plasma solution through relations (13). We postpone to a future work the analysis of interaction of a plasma with vacuum magnetic field in a non simply connected domain.

System (7), (9), (13), (14) is supplemented with initial conditions

$$U(0, x) = U_0(x), \quad x \in \Omega^+(0), \qquad \varphi(0, x) = \varphi_0(x), \quad x \in \Gamma(0),$$
$$\mathcal{H}(0, x) = \mathcal{H}_0(x), \quad x \in \Omega^-(0), \tag{15}$$

From the mathematical point of view, a natural wish is to find conditions on the initial data providing the existence and uniqueness on some time interval $[0, T]$ of a solution $(U, \mathcal{H}, \varphi)$ to problem (7), (9), (13)–(15) in Sobolev spaces. Since (1) is a system of hyperbolic conservation laws that can produce shock waves and other types of strong discontinuities (e.g., current-vortex sheets [36]), it is natural to expect to obtain only local-in-time existence theorems.

We must regard the boundary conditions on $H$ in (13), (14) as the restriction on the initial data (15). More precisely, we can prove that a solution of (7), (13), (14) (if it exists for all $t \in [0, T]$) satisfies

$$\text{div } H = 0 \quad \text{in } \Omega^+(t) \quad \text{and} \quad H_N = 0 \quad \text{on } \Gamma(t),$$

for all $t \in [0, T]$, if the latter were satisfied at $t = 0$, i.e., for the initial data (15). In particular, the fulfillment of div $H = 0$ implies that systems (1) and (7) are equivalent on solutions of problem (7), (13)–(15).

## 1.2. **An equivalent formulation in the fixed domain.**

We want to reduce the free boundary problem (7), (9), (13)–(15) to the fixed domain $\Omega$. For this purpose we introduce a suitable change of variables that is inspired by [18]. In all what follows, $H^s(\omega)$ denotes the Sobolev space of order $s$ on a domain $\omega$. We recall that on the torus $\mathbb{T}^2$, $H^s(\mathbb{T}^2)$ can be defined by means of the Fourier coefficients and coincides with the set of distributions $u$ such that

$$\sum_{k \in \mathbb{Z}^2} \left(1 + |k|^2\right)^s |c_k(u)|^2 < +\infty\,,$$

$c_k(u)$ denoting the $k$-th Fourier coefficient of $u$. In the following $\mathbb{T}^2$ is always identified with $\Gamma$. The following Lemma shows how to lift functions from $\Gamma$ to $\Omega$.

**Lemma 1.** *Let $m \geq 1$ be an integer. Then there exists a continuous linear map $\varphi \in H^{m-0.5}(\Gamma) \mapsto \Psi \in H^m(\Omega)$ such that $\Psi(0, x') = \varphi(x')$ on $\Gamma$, $\Psi(\pm 1, x') = 0$ on $\Gamma_\pm$, and moreover $\partial_1 \Psi(0, x') = 0$ if $m \geq 2$.*

The proof of Lemma 1 is given in [7]. The following Lemma gives the time-dependent version of Lemma 1.

**Lemma 2.** *Let $m \geq 1$ be an integer and let $T > 0$. Then there exists a continuous linear map*

$$\varphi \in \cap_{j=0}^{m-1} \mathcal{C}^j([0,T]; H^{m-j-0.5}(\Gamma)) \mapsto \Psi \in \cap_{j=0}^{m-1} \mathcal{C}^j([0,T]; H^{m-j}(\Omega))$$

*such that $\Psi(t, 0, x') = \varphi(t, x')$, $\Psi(t, \pm 1, x') = 0$, and moreover $\partial_1 \Psi(t, 0, x') = 0$ if $m \geq 2$. Furthermore, there exists a constant $C > 0$ that is independent of $T$ and only depends on $m$, such that*

$$\forall \varphi \in \cap_{j=0}^{m-1} \mathcal{C}^j([0,T]; H^{m-j-0.5}(\Gamma)), \quad \forall j = 0, \ldots, m-1, \quad \forall t \in [0,T],$$
$$\|\partial_t^j \Psi(t, \cdot)\|_{H^{m-j}(\Omega)} \leq C \|\partial_t^j \varphi(t, \cdot)\|_{H^{m-j-0.5}(\Gamma)}.$$

The proof of Lemma 2 is also given in [7]. The diffeomorphism that reduces the free boundary problem (7), (9), (13)–(14) to the fixed domain $\Omega$ is given in the following Lemma.

**Lemma 3.** *Let $m \geq 3$ be an integer. Then there exists a numerical constant $\epsilon_0 > 0$ such that for all $T > 0$, for all $\varphi \in \cap_{j=0}^{m-1} \mathcal{C}^j([0,T]; H^{m-j-0.5}(\Gamma))$ satisfying $\|\varphi\|_{\mathcal{C}([0,T]; H^{2.5}(\Gamma))} \leq \epsilon_0$, the function*

$$\Phi(t, x) := (x_1 + \Psi(t, x), x'), \qquad (t, x) \in [0, T] \times \Omega, \tag{16}$$

*with $\Psi$ as in Lemma 2, defines an $H^m$-diffeomorphism of $\Omega$ for all $t \in [0, T]$. Moreover, there holds $\partial_t^j \Phi \in \mathcal{C}([0,T]; H^{m-j}(\Omega))$ for $j = 0, \ldots, m-1$, $\Phi(t, 0, x') = (\varphi(t, x'), x')$, $\Phi(t, \pm 1, x') = (\pm 1, x')$, $\partial_1 \Phi(t, 0, x') = (1, 0, 0)$, and*

$$\forall t \in [0, T], \quad \|\Psi(t, \cdot)\|_{W^{1,\infty}(\Omega)} \leq \frac{1}{2}.$$

*Proof.* The proof follows directly from Lemma 2 and the Sobolev imbedding Theorem, because

$$\partial_1 \Phi_1(t, x) = 1 + \partial_1 \Psi(t, x) \geq 1 - \|\Psi(t, \cdot)\|_{\mathcal{C}([0,T]; W^{1,\infty}(\Omega))}$$
$$\geq 1 - C \|\varphi\|_{\mathcal{C}([0,T]; H^{2.5}(\Gamma))} \geq 1/2,$$

provided that $\varphi$ is taken sufficiently small in $\mathcal{C}([0,T]; H^{2.5}(\Gamma))$. In the latter inequality, $C$ denotes a numerical constant. The other properties of $\Psi$ follow directly from Lemma 2. $\qquad\square$

We introduce the change of independent variables defined by (16) by setting

$$\widetilde{U}(t, x) := U(t, \Phi(t, x)), \quad \widetilde{\mathcal{H}}(t, x) := \mathcal{H}(t, \Phi(t, x)).$$

Dropping for convenience tildes in $\widetilde{U}$, $\widetilde{\mathcal{H}}$, problem (7), (9) (13)–(15) can be reformulated on the fixed reference domain $\Omega$ as

$$\mathbb{P}(U, \Psi) = 0 \quad \text{in } [0,T] \times \Omega^+, \quad \mathbb{V}(\mathcal{H}, \Psi) = 0 \quad \text{in } [0,T] \times \Omega^-, \tag{17}$$

$$\mathbb{B}(U, \mathcal{H}, \varphi) = \bar{\mathfrak{J}} \quad \text{on } [0,T] \times (\Gamma^3 \times \Gamma_+ \times \Gamma_-), \tag{18}$$

$$(U, \mathcal{H})|_{t=0} = (U_0, \mathcal{H}_0) \quad \text{in } \Omega^+ \times \Omega^-, \qquad \varphi|_{t=0} = \varphi_0 \quad \text{on } \Gamma, \tag{19}$$

where $\mathbb{P}(U, \Psi) = P(U, \Psi)U$,

$$P(U, \Psi) = A_0(U)\partial_t + \widetilde{A}_1(U, \Psi)\partial_1 + A_2(U)\partial_2 + A_3(U)\partial_3,$$

$$\widetilde{A}_1(U, \Psi) = \frac{1}{\partial_1 \Phi_1}\Big(A_1(U) - A_0(U)\partial_t \Psi - \sum_{k=2}^{3} A_k(U)\partial_k \Psi\Big),$$

$$\mathbb{V}(\mathcal{H}, \Psi) = \begin{pmatrix} \nabla \times \mathfrak{H} \\ \text{div } \mathfrak{h} \end{pmatrix},$$

$$\mathfrak{H} = (\mathcal{H}_1 \partial_1 \Phi_1, \mathcal{H}_{\tau_2}, \mathcal{H}_{\tau_3}), \quad \mathfrak{h} = (\mathcal{H}_N, \mathcal{H}_2 \partial_1 \Phi_1, \mathcal{H}_3 \partial_1 \Phi_1),$$

$$\mathcal{H}_N = \mathcal{H}_1 - \mathcal{H}_2 \partial_2 \Psi - \mathcal{H}_3 \partial_3 \Psi, \quad \mathcal{H}_{\tau_i} = \mathcal{H}_1 \partial_i \Psi + \mathcal{H}_i, \quad i = 2, 3,$$

$$\mathbb{B}(U, \mathcal{H}, \varphi) = \begin{pmatrix} \partial_t \varphi - v_N \\ [q] \\ \mathcal{H}_N \\ v_1 \\ \nu \times \mathcal{H} \end{pmatrix}, \quad [q] = q_{|x_1=0+} - \frac{1}{2}|\mathcal{H}|^2_{|x_1=0-},$$

$$v_N = v_1 - v_2 \partial_2 \Psi - v_3 \partial_3 \Psi, \qquad \bar{\mathfrak{J}} = (0, 0, 0, 0, \mathfrak{J})^T.$$

In (18) the notation $[0, T] \times (\Gamma^3 \times \Gamma_+ \times \Gamma_-)$ means that the first three components of this vector equation are taken on $[0, T] \times \Gamma$, the fourth one on $[0, T] \times \Gamma_+$, and the fifth one on $[0, T] \times \Gamma_-$. To avoid an overload of notation we have denoted by the same symbols $v_N, \mathcal{H}_N$ here above and $v_N, \mathcal{H}_N$ as in (13). Notice that $v_{N|x_1=0} = v_1 - v_2\partial_2\varphi - v_3\partial_3\varphi$, $\mathcal{H}_{N|x_1=0} = \mathcal{H}_1 - \mathcal{H}_2\partial_2\varphi - \mathcal{H}_3\partial_3\varphi$, as in the previous definition in (13).

We did not include in problem (17)–(19) the equation

$$\operatorname{div} h = 0 \quad \text{in } [0, T] \times \Omega^+, \tag{20}$$

and the boundary conditions

$$H_N = 0 \quad \text{on } [0, T] \times \Gamma, \qquad H_1 = 0 \quad \text{on } [0, T] \times \Gamma_+, \tag{21}$$

where $h = (H_N, H_2\partial_1\Phi_1, H_3\partial_1\Phi_1)$, $H_N = H_1 - H_2\partial_2\Psi - H_3\partial_3\Psi$, because they are just restrictions on the initial data (19). More precisely, referring to [36] for the proof, we have the following proposition.

**Proposition 4.** *Let the initial data* (19) *satisfy* (20) *and* (21) *for $t = 0$. If $(U, \mathcal{H}, \varphi)$ is a solution of problem* (17)–(19), *then this solution satisfies* (20) *and* (21) *for all $t \in [0, T]$.*

Note that Proposition 4 stays valid if we replace (17) by system (1) in the straightened variables. This means that these systems are equivalent on solutions of our plasma-vacuum interface problem and we may justifiably replace the conservation laws (1) by their nonconservative form (7).

In [35] we proved the well-posedness of the linearized problem associated to the nonlinear problem (17)–(19) in Sobolev spaces[1] provided that the "unperturbed flow" (a basic state) satisfies the hyperbolicity condition (8) and the stability condition

$$|H \times \mathcal{H}|_{x_1=0} \geq \delta_0 > 0, \tag{22}$$

where $\delta_0$ is a fixed constant. Since the basic state in [35] was also assumed to satisfy (21) and the third boundary condition in (18), one can show that the stability condition (22) is equivalently rewritten as

$$|H_2\mathcal{H}_3 - H_3\mathcal{H}_2|_{x_1=0} \geq \delta_0 > 0. \tag{23}$$

Now our main goal is to prove the well-posedness of the nonlinear problem (17)–(19) in suitable anisotropic Sobolev spaces (see Sect. 2) provided that the initial data (19) satisfy the hyperbolicity condition (8) and the stability condition (22) (together with appropriate compatibility conditions).

## 2. FUNCTION SPACES

Now we introduce the main function spaces to be used in the following. Let us denote

$$\begin{aligned} Q_T &:= (-\infty, T] \times \Omega, \quad Q_T^\pm := (-\infty, T] \times \Omega^\pm, \\ \omega_T &:= (-\infty, T] \times \Gamma, \quad \omega_T^\pm := (-\infty, T] \times \Gamma_\pm. \end{aligned} \tag{24}$$

2.1. **Weighted Sobolev spaces.** For $\gamma \geq 1$ and $s \in \mathbb{N}$, $H_\gamma^s(\Omega)$ will denote the Sobolev space of order $s$, equipped with the $\gamma$–depending norm defined by

$$||u||_{H_\gamma^s(\Omega)}^2 := \sum_{|\alpha| \leq s} \gamma^{2(s-|\alpha|)} ||\partial^\alpha u||_{L^2(\Omega)}^2 \,.$$

For functions defined over $Q_T$ we will consider the weighted Sobolev spaces $H_\gamma^m(Q_T)$ equipped with the $\gamma$–depending norm

$$||u||_{H_\gamma^m(Q_T)}^2 := \sum_{|\alpha| \leq m} \gamma^{2(m-|\alpha|)} ||\partial^\alpha u||_{L^2(Q_T)}^2 \,.$$

Similar weighted Sobolev spaces will be considered for functions defined on $\Omega^\pm, Q_T^\pm$.

---

[1]More precisely, the well-posedness in so-called conormal Sobolev spaces was proved (see Sect. 2 for their definition).

## 2.2. Conormal Sobolev spaces.

Let us consider functions defined over $Q_T^+$. For $j = 0, \dots, 3$, we set

$$Z_0 = \partial_t, \quad Z_1 := \sigma(x_1)\partial_1, \quad Z_j := \partial_j, \text{ for } j = 2, 3,$$

where $\sigma(x_1) \in C^\infty(0, 1)$ is a positive function such that $\sigma(x_1) = x_1$ in a neighborhood of the origin and $\sigma(x_1) = 1 - x_1$ in a neighborhood of $x_1 = 1$. Then, for every multi-index $\alpha = (\alpha_0, \dots, \alpha_3) \in \mathbb{N}^4$, the *conormal* derivative $Z^\alpha$ is defined by

$$Z^\alpha := Z_0^{\alpha_0} \dots Z_3^{\alpha_3};$$

we also write $\partial^\alpha = \partial_0^{\alpha_0} \dots \partial_3^{\alpha_3}$ for the usual partial derivative corresponding to $\alpha$.

Given an integer $m \geq 1$, the *conormal Sobolev space* $H_{tan}^m(Q_T^+)$ is defined as the set of functions $u \in L^2(Q_T^+)$ such that $Z^\alpha u \in L^2(Q_T^+)$, for all multi-indices $\alpha$ with $|\alpha| \leq m$ (see [22, 23]). Agreeing with the notations set for the usual Sobolev spaces, for $\gamma \geq 1$, $H_{tan,\gamma}^m(Q_T^+)$ will denote the conormal space of order $m$ equipped with the $\gamma$−depending norm

$$||u||_{H_{tan,\gamma}^m(Q_T^+)}^2 := \sum_{|\alpha| \leq m} \gamma^{2(m-|\alpha|)} ||Z^\alpha u||_{L^2(Q_T^+)}^2 \qquad (25)$$

and we have $H_{tan}^m(Q_T^+) := H_{tan,1}^m(Q_T^+)$. Similar conormal Sobolev spaces with $\gamma$-depending norms will be considered for functions defined on $\Omega^\pm$ (disregarding $Z_0$ derivatives), $Q_T^{-}$ [2].

## 2.3. Anisotropic Sobolev spaces.

Keeping the same notations used above, for every positive integer $m$ the *anisotropic Sobolev space* $H_*^m(\Omega^+)$ is defined as

$$H_*^m(\Omega^+) := \{w \in L^2(\Omega^+) : Z^\alpha \partial_1^k w \in L^2(\Omega^+), |\alpha| + 2k \leq m\}.$$

For the sake of convenience we also set $H_*^0(\Omega^+) = H_{tan}^0(\Omega^+) = L^2(\Omega^+)$. For an extensive study of the anisotropic spaces $H_*^m(\Omega^+)$ we refer the reader to [22, 24, 32] and references therein. We observe that

$$H^m(\Omega^+) \hookrightarrow H_*^m(\Omega^+) \hookrightarrow H_{tan}^m(\Omega^+) \subset H_{loc}^m(\Omega^+), \qquad (26)$$
$$H_*^m(\Omega^+) \hookrightarrow H^{[m/2]}(\Omega^+), \quad H_*^1(\Omega^+) = H_{tan}^1(\Omega^+)$$

(except for $H_{loc}^m(\Omega^+)$ all imbeddings are continuous). The anisotropic space $H_{*,\gamma}^m(\Omega^+)$ is the same space equipped with the $\gamma$-depending norm

$$||w||_{H_{*,\gamma}^m(\Omega^+)}^2 := \sum_{|\alpha| + 2k \leq m} \gamma^{2(m-|\alpha|-2k)} ||Z^\alpha \partial_1^k w||_{L^2(\Omega^+)}^2. \qquad (27)$$

We have $H_*^m(\Omega^+) = H_{*,1}^m(\Omega^+)$. Some useful properties of the $\gamma$-dependent space $H_{*,\gamma}^m(\Omega^+)$, used in this paper, are listed in Appendix A.

In a similar way we define the anisotropic space $H_{*,\gamma}^m(Q_T^+)$, equipped with its natural norm.

We will use the same notation for spaces of scalar and vector-valued functions.

## 3. Main result and discussion

Let us now state our main result:

**Theorem 5.** *Let $m \in \mathbb{N}$, $m \geq 13$, and $\mathfrak{J} \in H^{m+9}([0, T_0] \times \Gamma_-)$ for some $T_0 > 0$. Consider initial data $U_0 \in H^{m+9.5}(\Omega^+)$, $\mathcal{H}^0 \in H^{m+9.5}(\Omega^-)$, and $\varphi_0 \in H^{m+10}(\Gamma)$. Moreover, the initial data satisfy (8), (22), (129) and (130) and are compatible up to order $m + 9$ in the sense of Definition 20. Then there exists $0 < T \leq T_0$, $\epsilon_1 > 0$ such that, if $\|\varphi_0\|_{H^{2.5}(\Gamma)} \leq \epsilon_1$, problem (17)–(19) has a unique solution $(U, \mathcal{H}, \varphi)$ in $[0, T]$, with*

$$U \in H_{*,\gamma}^m(]0, T[\times\Omega^+), \quad \mathcal{H} \in H_\gamma^m(]0, T[\times\Omega^-), \quad \varphi \in H_\gamma^{m+1/2}(]0, T[\times\Gamma).$$

**Remark 6.** *The initial vacuum magnetic field $\mathcal{H}^0 \in H^{m+9.5}(\Omega^-)$ is not given independently of the other initial data. In fact, as it is assumed to satisfy (130), which is an uniquely solvable elliptic system, $\mathcal{H}^0$ is uniquely determined from $\varphi_0$ (i.e. the initial space domain) and $\mathfrak{J}(0)$ (the external density current at initial time) by Theorem 13. In this sense, for a given $\mathfrak{J}$, the actual initial data of the problem may be considered only $U_0$, $\varphi_0$.*

---

[2] On $\Omega^-$ or $Q_T^-$ $Z_1$ is defined by $Z_1 := \sigma(-x_1)\partial_1$.

Theorem 5 shows that the stability condition (22) is sufficient for nonlinear well-posedness. As far as we know, it is not known what happens if (22) is violated, whether there is some form of strong/weak stability[3] or a transition to instability implying the ill-posedness of the original nonlinear problem. The study of the well-posedness of the plasma-vacuum interface problem for the case when condition (22) is violated, i.e., at some points of the initial interface the plasma and vacuum magnetic fields may be parallel to each other (or one of them is zero), is postponed to the future.

The remainder of the paper is organized as follows. In the next Section 4 we formulate the linearized problem associated to (17)–(19) and introduce suitable decompositions of the magnetic fields to reduce it to that with homogeneous boundary conditions and homogeneous linearized "vacuum" equations. In fact, for proving the basic a priori energy estimate, in [35] it is convenient to have the vacuum magnetic field satisfying homogeneous equations and boundary conditions as in (41), and plasma magnetic field satisfying homogeneous constraints (53), (54). Thus we introduce the decomposition $\dot{\mathcal{H}} = \mathcal{H}' + \mathcal{H}''$ in the vacuum side, with $\mathcal{H}'$ solution of (41), and $\mathcal{H}''$ taking all the nonhomogeneous part (40), and the decomposition (49) in the plasma side.

In Section 5, for convenience of the reader we recall the well-posedness result of [35] for the reduced homogeneous linearized problem. In Section 6, for each fixed time $t$, we study the nonhomogeneous elliptic system (40) for the component $\mathcal{H}''$ of the vacuum magnetic field. Here we work in suitable function spaces taking account of the particular geometry and the different conditions on the upper and lower boundary of $\Omega^-$. An important point is the direct $L^2$ estimate of the solution by negative $H^{-1}$ norms of the data, inspired from [3], which will be crucial in Section 8 when dealing with commutators.

In Section 7 we obtain the final form of the $H^1$ estimate for the full linearized problem (39), see Theorem 15. In Section 8 we deduce a so-called tame estimate and show the well-posedness of the linearized problem in anisotropic Sobolev spaces of an arbitrary fixed order of smoothness, see Theorem 16. In the plasma side the main difficulty comes from the characteristic boundary, forcing to work in anisotropic Sobolev spaces $H_*^m$ with different regularity in the normal and tangential directions. Moreover, due to the loss of one derivative with respect to the source terms in the basic $H^1$ estimate (86), special care is needed in the estimate of commutators containing the solution itself, as they can't be treated in the usual way as zero order terms. Here we use the calculus tools developed in [22] (see Appendix A), in particular for the cases when only conormal regularity is involved. In the vacuum side the main difficulty comes from the commutators with conormal derivatives. Here we use the a priori estimates of Section 6 with negative norms $H^{-1}$ (in space) to compensate the too higher $H^2$ regularity (in time) appearing in the right-hand side of the basic $H^1$ estimate (86).

In Sections 9–13 we give the proof of our main Theorem 5. In particular, Section 11 is devoted to the description of the Nash-Moser iteration scheme, similar to that of [2, 8, 36], and Section 12 to the proof of its convergence by induction. The main difficulty is that at each iteration the inversion of the operator $(\mathbb{L}', \mathbb{V}', \mathbb{B}')$ requires the linearization around a state satisfying the constraints (29)–(34), (61), that is the constraints of the basic state given in Section 4. We thus need to introduce a smooth modified state, denoted $V_{n+1/2}, \mathcal{K}^{n+1/2}, \psi_{n+1/2}$, that satisfies the above mentioned constraints; the exact definition of this intermediate state is detailed in subsection 12.4. A similar difficulty was found in [8, 36]. Section 13 is devoted to the proof of the uniqueness of a smooth solution.

At last, in Appendix A we recall some technical results about the anisotropic Sobolev spaces and some useful calculus inequalities for them whereas in Appendix B we recall for the reader's convenience some well-known commutator estimates and Moser-type inequalities for standard Sobolev spaces. Moreover, in Appendix C we briefly explain minor modifications necessary to adapt the energy a priori estimate obtained in [35] for the linearized problem for the case with the whole space domain to our present case with the added fixed top and bottom boundaries $\Gamma_\pm$.

## 4. THE LINEARIZED PROBLEM

---

[3] Strictly speaking, in this paper by stability we mean the *well-posedness* of the problem resulting from the linearization about a given (generally speaking, non-stationary) basic state.

### 4.1. **Basic state.** Let

$$(\widehat{U}(t,x), \widehat{\mathcal{H}}(t,x), \hat{\varphi}(t,x')) \tag{28}$$

be a given sufficiently smooth vector-function with $\widehat{U} = (\hat{q}, \hat{v}, \widehat{H}, \widehat{S})$, respectively defined on $Q_T^+, Q_T^-, \omega_T$, with

$$\|\widehat{U}\|_{H^9_{*,\gamma}(Q_T^+)} + \|\widehat{\mathcal{H}}\|^2_{H^9_\gamma(Q_T^-)} + \|\hat{\varphi}\|_{H^{9.5}_\gamma(\omega_T)} \le K, \qquad \|\hat{\varphi}\|_{\mathcal{C}([0,T];H^{2.5}(\Gamma))} \le \epsilon_0, \tag{29}$$

where $K > 0$ is a constant and $\epsilon_0$ is the arbitrary constant introduced in Lemma 3. Corresponding to the given $\hat{\varphi}$ we construct $\widehat{\Psi}$ and the diffeomorphism $\widehat{\Phi}$ as in Lemmata 2 and 3, such that

$$\partial_1 \widehat{\Phi}_1 \ge 1/2.$$

Notice that (29) yields[4]

$$\|\widehat{U}\|_{W^{2,\infty}(Q_T^+)} + \|\partial_1 \widehat{U}\|_{W^{2,\infty}(Q_T^+)} + \|\widehat{\mathcal{H}}\|_{W^{2,\infty}(Q_T^-)} + \|\nabla_{t,x}\widehat{\Psi}\|_{W^{2,\infty}(Q_T)} \le C(K),$$

where $\nabla_{t,x} = (\partial_t, \nabla)$ and $C = C(K) > 0$ is a constant depending on $K$.

We assume that the basic state (28) satisfies (for some positive $\rho_0, \rho_1 \in \mathbb{R}$)

$$\rho(\hat{p}, \widehat{S}) \ge \rho_0 > 0, \quad \rho_p(\hat{p}, \widehat{S}) \ge \rho_1 > 0 \quad \text{in } \overline{Q}_T^+, \tag{30}$$

$$\partial_t \widehat{H} + \frac{1}{\partial_1 \widehat{\Phi}_1} \left\{ (\hat{w} \cdot \nabla)\widehat{H} - (\hat{h} \cdot \nabla)\hat{v} + \widehat{H}\operatorname{div}\hat{u} \right\} = 0 \quad \text{in } Q_T^+, \tag{31}$$

$$\operatorname{div} \hat{\mathfrak{h}} = 0 \quad \text{in } Q_T^-, \tag{32}$$

$$\partial_t \hat{\varphi} - \hat{v}_N = 0, \quad \widehat{\mathcal{H}}_N = 0 \quad \text{on } \omega_T, \quad \hat{v}_1 = 0 \quad \text{on } \omega_T^+, \quad \nu \times \hat{\mathcal{H}} = \mathfrak{J} \quad \text{on } \omega_T^-, \tag{33}$$

where all the "hat" values are determined like corresponding values for $(U, \mathcal{H}, \varphi)$, i.e.

$$\widehat{\mathfrak{H}} = (\widehat{\mathcal{H}}_1 \partial_1 \widehat{\Phi}_1, \widehat{\mathcal{H}}_{\tau_2}, \widehat{\mathcal{H}}_{\tau_3}), \quad \hat{\mathfrak{h}} = (\widehat{\mathcal{H}}_N, \widehat{\mathcal{H}}_2 \partial_1 \widehat{\Phi}_1, \widehat{\mathcal{H}}_3 \partial_1 \widehat{\Phi}_1), \quad \hat{h} = (\widehat{H}_N, \widehat{H}_2 \partial_1 \widehat{\Phi}_1, \widehat{H}_3 \partial_1 \widehat{\Phi}_1),$$

$$\widehat{H}_N = \widehat{H}_1 - \widehat{H}_2 \partial_2 \widehat{\Psi} - \widehat{H}_3 \partial_3 \widehat{\Psi}, \quad \widehat{\mathcal{H}}_N = \widehat{\mathcal{H}}_1 - \widehat{\mathcal{H}}_2 \partial_2 \widehat{\Psi} - \widehat{\mathcal{H}}_3 \partial_3 \widehat{\Psi},$$

$$\hat{p} = \hat{q} - |\widehat{H}|^2/2, \quad \hat{v}_N = \hat{v}_1 - \hat{v}_2 \partial_2 \widehat{\Psi} - \hat{v}_3 \partial_3 \widehat{\Psi},$$

$$\hat{u} = (\hat{v}_N, \hat{v}_2 \partial_1 \widehat{\Phi}_1, \hat{v}_3 \partial_1 \widehat{\Phi}_1), \quad \hat{w} = \hat{u} - (\partial_t \widehat{\Psi}, 0, 0).$$

It follows from (31) that the constraints

$$\operatorname{div} \hat{h} = 0 \quad \text{in } Q_T^+, \quad \widehat{H}_N = 0 \quad \text{on } \omega_T, \quad \widehat{H}_1 = 0 \quad \text{on } \omega_T^+, \tag{34}$$

are satisfied for the basic state (28) if they hold at $t = 0$ (see [36] for the proof). Thus, for the basic state we also require the fulfillment of conditions (34) at $t = 0$.

### 4.2. **Linearized problem.** The linearized equations for (17), (18) read:

$$\mathbb{P}'(\widehat{U}, \widehat{\Psi})(\delta U, \delta \Psi) := \frac{\mathrm{d}}{\mathrm{d}\varepsilon} \mathbb{P}(U_\varepsilon, \Psi_\varepsilon)|_{\varepsilon=0} = f \quad \text{in } Q_T^+,$$

$$\mathbb{V}'(\widehat{\mathcal{H}}, \widehat{\Psi})(\delta \mathcal{H}, \delta \Psi) := \frac{\mathrm{d}}{\mathrm{d}\varepsilon} \mathbb{V}(\mathcal{H}_\varepsilon, \Psi_\varepsilon)|_{\varepsilon=0} = \mathcal{G}' \quad \text{in } Q_T^-,$$

$$\mathbb{B}'(\widehat{U}, \widehat{\mathcal{H}}, \hat{\varphi})(\delta U, \delta \mathcal{H}, \delta \varphi) := \frac{\mathrm{d}}{\mathrm{d}\varepsilon} \mathbb{B}(U_\varepsilon, \mathcal{H}_\varepsilon, \varphi_\varepsilon)|_{\varepsilon=0} = g \quad \text{on } \omega_T^3 \times \omega_T^\pm,$$

where $U_\varepsilon = \widehat{U} + \varepsilon \delta U$, $\mathcal{H}_\varepsilon = \widehat{\mathcal{H}} + \varepsilon \delta \mathcal{H}$, $\varphi_\varepsilon = \hat{\varphi} + \varepsilon \delta \varphi$; $\delta \Psi$ is constructed from $\delta \varphi$ as in Lemma 2 and $\Psi_\varepsilon = \widehat{\Psi} + \varepsilon \delta \Psi$. In the above boundary equation the first three components are taken on $\omega_T$, the fourth one on $\omega_T^+$, and the fifth one on $\omega_T^-$. Here we introduce the source terms $f = (f_1, \ldots, f_8)$, $\mathcal{G}' = (\chi, \Xi)$, $\chi = (\chi_1, \chi_2, \chi_3)$, and $g = (g_1, g_2, g_3)$ to make the interior equations and the boundary conditions inhomogeneous.

We compute the exact form of the linearized equations (below we drop $\delta$):

$$\mathbb{P}'(\widehat{U}, \widehat{\Psi})(U, \Psi) = P(\widehat{U}, \widehat{\Psi})U + \mathcal{C}(\widehat{U}, \widehat{\Psi})U - \left\{ P(\widehat{U}, \widehat{\Psi})\Psi \right\} \frac{\partial_1 \widehat{U}}{\partial_1 \widehat{\Phi}_1} = f,$$

---

[4]This inequality is taken as an assumption in [35].

$$\mathbb{V}'(\widehat{\mathcal{H}}, \widehat{\Psi})(\mathcal{H}, \Psi) = \mathbb{V}(\mathcal{H}, \widehat{\Psi}) + \begin{pmatrix} \nabla\widehat{\mathcal{H}}_1 \times \nabla\Psi \\ \nabla \times \begin{pmatrix} 0 \\ -\widehat{\mathcal{H}}_3 \\ \widehat{\mathcal{H}}_2 \end{pmatrix} \cdot \nabla\Psi \end{pmatrix} = \mathcal{G}',$$

$$\mathbb{B}'(\widehat{U}, \widehat{\mathcal{H}}, \hat{\varphi})(U, \mathcal{H}, \varphi) = \begin{pmatrix} \partial_t\varphi + \hat{v}_2\partial_2\varphi + \hat{v}_3\partial_3\varphi - v_N \\ q - \widehat{\mathcal{H}} \cdot \mathcal{H} \\ \mathcal{H}_N - \widehat{\mathcal{H}}_2\partial_2\varphi - \widehat{\mathcal{H}}_3\partial_3\varphi \\ v_1 \\ \nu \times \mathcal{H} \end{pmatrix} = g,$$

where $v_N := v_1 - v_2\partial_2\widehat{\Psi} - v_3\partial_3\widehat{\Psi}$ and the matrix $\mathcal{C}(\widehat{U}, \widehat{\Psi})$ is determined as follows:

$$\mathcal{C}(\widehat{U}, \widehat{\Psi})Y = (Y, \nabla_y A_0(\widehat{U}))\partial_t\widehat{U} + (Y, \nabla_y\widetilde{A}_1(\widehat{U}, \widehat{\Psi}))\partial_1\widehat{U}$$
$$+ (Y, \nabla_y A_2(\widehat{U}))\partial_2\widehat{U} + (Y, \nabla_y A_3(\widehat{U}))\partial_3\widehat{U},$$

$$(Y, \nabla_y A(V)) := \sum_{i=1}^{8} y_i \left(\frac{\partial A(Y)}{\partial y_i}\Big|_{Y=V}\right), \quad Y = (y_1, \dots, y_8).$$

Since the differential operators $\mathbb{P}'(\widehat{U}, \widehat{\Psi})$ and $\mathbb{V}'(\widehat{\mathcal{H}}, \widehat{\Psi})$ are first-order operators in $\Psi$, as in [1] the linearized problem is rewritten in terms of the "good unknown"

$$\dot{U} := U - \frac{\Psi}{\partial_1\widehat{\Phi}_1}\partial_1\widehat{U}, \quad \dot{\mathcal{H}} := \mathcal{H} - \frac{\Psi}{\partial_1\widehat{\Phi}_1}\partial_1\widehat{\mathcal{H}}. \tag{35}$$

Taking into account assumptions (33) and omitting detailed calculations, we rewrite our linearized equations in terms of the new unknowns (35):

$$\mathbb{P}'(\widehat{U}, \widehat{\Psi})(U, \Psi) = P(\widehat{U}, \widehat{\Psi})\dot{U} + \mathcal{C}(\widehat{U}, \widehat{\Psi})\dot{U} + \frac{\Psi}{\partial_1\widehat{\Phi}_1}\partial_1\{\mathbb{P}(\widehat{U}, \widehat{\Psi})\} = f,$$
$$\mathbb{V}'(\widehat{\mathcal{H}}, \widehat{\Psi})(\mathcal{H}, \Psi) = \mathbb{V}(\dot{\mathcal{H}}, \widehat{\Psi}) + \frac{\Psi}{\partial_1\widehat{\Phi}_1}\partial_1\{\mathbb{V}(\widehat{\mathcal{H}}, \widehat{\Psi})\} = \mathcal{G}', \tag{36}$$

$$\mathbb{B}'_e(\widehat{U}, \widehat{\mathcal{H}}, \hat{\varphi})(\dot{U}, \dot{\mathcal{H}}, \varphi) := \mathbb{B}'(\widehat{U}, \widehat{\mathcal{H}}, \hat{\varphi})(U, \mathcal{H}, \varphi)$$
$$= \begin{pmatrix} \partial_t\varphi + \hat{v}_2\partial_2\varphi + \hat{v}_3\partial_3\varphi - \dot{v}_N - \varphi\,\partial_1\hat{v}_N \\ \dot{q} - \widehat{\mathcal{H}} \cdot \dot{\mathcal{H}} + [\partial_1\hat{q}]\varphi \\ \dot{\mathcal{H}}_N - \partial_2(\widehat{\mathcal{H}}_2\varphi) - \partial_3(\widehat{\mathcal{H}}_3\varphi) \\ \dot{v}_1 \\ \nu \times \dot{\mathcal{H}} \end{pmatrix} = g, \tag{37}$$

where $\dot{v}_N = \dot{v}_1 - \dot{v}_2\partial_2\widehat{\Psi} - \dot{v}_3\partial_3\widehat{\Psi}$, $\dot{\mathcal{H}}_N = \dot{\mathcal{H}}_1 - \dot{\mathcal{H}}_2\partial_2\widehat{\Psi} - \dot{\mathcal{H}}_3\partial_3\widehat{\Psi}$, and

$$[\partial_1\hat{q}] = (\partial_1\hat{q})|_{x_1=0} - (\widehat{\mathcal{H}} \cdot \partial_1\widehat{\mathcal{H}})|_{x_1=0}.$$

We used assumption (32), taken at $x_1 = 0$, while writing down the third boundary condition in (37).

As in [1, 8, 36], we drop the zeroth-order terms in $\Psi$ in (36) and consider the effective linear operators

$$\mathbb{P}'_e(\widehat{U}, \widehat{\Psi})\dot{U} := P(\widehat{U}, \widehat{\Psi})\dot{U} + \mathcal{C}(\widehat{U}, \widehat{\Psi})\dot{U} = f,$$
$$\mathbb{V}'_e(\widehat{\mathcal{H}}, \widehat{\Psi})\dot{\mathcal{H}} := \mathbb{V}(\dot{\mathcal{H}}, \widehat{\Psi}) = \mathcal{G}'. \tag{38}$$

In the future nonlinear analysis of Section 12 the dropped terms in (36) should be considered as error terms. With the new form (38), (37) of the linearized equations, our linearized problem for $(\dot{U}, \dot{\mathcal{H}}, \varphi)$

reads in explicit form:

$$\widehat{A}_0 \partial_t \dot{U} + \sum_{j=1}^{3} \widehat{A}_j \partial_j \dot{U} + \widehat{\mathcal{C}} \dot{U} = f \qquad \text{in } Q_T^+, \tag{39a}$$

$$\nabla \times \dot{\mathfrak{H}} = \chi, \quad \text{div } \dot{\mathfrak{h}} = \Xi \qquad \text{in } Q_T^-, \tag{39b}$$

$$\partial_t \varphi = \dot{v}_N - \hat{v}_2 \partial_2 \varphi - \hat{v}_3 \partial_3 \varphi + \varphi \, \partial_1 \hat{v}_N + g_1, \tag{39c}$$

$$\dot{q} = \widehat{\mathcal{H}} \cdot \dot{\mathcal{H}} - [\partial_1 \hat{q}] \varphi + g_2, \tag{39d}$$

$$\dot{\mathcal{H}}_N = \partial_2 (\widehat{\mathcal{H}}_2 \varphi) + \partial_3 (\widehat{\mathcal{H}}_3 \varphi) + g_3 \qquad \text{on } \omega_T, \tag{39e}$$

$$\dot{v}_1 = g_4 \quad \text{on } \omega_T^+, \qquad \nu \times \dot{\mathcal{H}} = g_5 \quad \text{on } \omega_T^-, \tag{39f}$$

$$(\dot{U}, \dot{\mathcal{H}}, \varphi) = 0 \qquad \text{for } t < 0, \tag{39g}$$

where

$$\widehat{A}_\alpha =: A_\alpha(\widehat{U}), \quad \alpha = 0, 2, 3, \quad \widehat{A}_1 =: \widetilde{A}_1(\widehat{U}, \widehat{\Psi}), \quad \widehat{\mathcal{C}} := \mathcal{C}(\widehat{U}, \widehat{\Psi}),$$

$$\dot{\mathfrak{H}} = (\dot{\mathcal{H}}_1 \partial_1 \widehat{\Phi}_1, \dot{\mathcal{H}}_{\tau_2}, \dot{\mathcal{H}}_{\tau_3}), \quad \dot{\mathfrak{h}} = (\dot{\mathcal{H}}_N, \dot{\mathcal{H}}_2 \partial_1 \widehat{\Phi}_1, \dot{\mathcal{H}}_3 \partial_1 \widehat{\Phi}_1),$$

$$\dot{\mathcal{H}}_N = \dot{\mathcal{H}}_1 - \dot{\mathcal{H}}_2 \partial_2 \widehat{\Psi} - \dot{\mathcal{H}}_3 \partial_3 \widehat{\Psi}, \quad \dot{\mathcal{H}}_{\tau_i} = \dot{\mathcal{H}}_1 \partial_i \widehat{\Psi} + \dot{\mathcal{H}}_i, \quad i = 2, 3.$$

For the resolution of the elliptic problem (39b), (39e), (39f) the data $\chi, g_5$ must satisfy necessary compatibility conditions described in (66).

We assume that the source terms $f, \chi, \Xi$ and the boundary data $g$ vanish in the past and consider the case of zero initial data. We postpone the case of nonzero initial data to the nonlinear analysis (see e.g. [8, 36]).

### 4.3. Reduction to homogeneous constraints in the "vacuum part".

We decompose $\dot{\mathcal{H}}$ in (39) as $\dot{\mathcal{H}} = \mathcal{H}' + \mathcal{H}''$ (and accordingly $\dot{\mathfrak{H}} = \mathfrak{H}' + \mathfrak{H}''$, $\dot{\mathfrak{h}} = \mathfrak{h}' + \mathfrak{h}''$), where $\mathcal{H}''$ is required to solve for each $t$ the elliptic problem

$$\begin{aligned} \nabla \times \mathfrak{H}'' &= \chi, \quad \text{div } \mathfrak{h}'' = \Xi \quad \text{in } Q_T^-, \\ \mathfrak{h}_1'' &= \mathcal{H}_N'' = g_3 \qquad \text{on } \omega_T, \\ \nu \times \mathcal{H}'' &= g_5 \qquad \text{on } \omega_T^-. \end{aligned} \tag{40}$$

For the resolution of (40) the data $\chi, g_5$ must satisfy the necessary compatibility conditions (66). The resolution of (40) is given in Section 6.

Given $\mathcal{H}''$, we look for $\mathcal{H}'$ such that

$$\begin{aligned} \nabla \times \mathfrak{H}' &= 0, \quad \text{div } \mathfrak{h}' = 0 \qquad \text{in } Q_T^-, \\ \dot{q} &= \widehat{\mathcal{H}} \cdot \mathcal{H}' - [\partial_1 \hat{q}] \varphi + g_2', \\ \mathcal{H}_N' &= \partial_2 (\widehat{\mathcal{H}}_2 \varphi) + \partial_3 (\widehat{\mathcal{H}}_3 \varphi) \quad \text{on } \omega_T, \\ \nu \times \mathcal{H}' &= 0 \qquad \text{on } \omega_T^-, \end{aligned} \tag{41}$$

where we have denoted $g_2' = g_2 + \widehat{\mathcal{H}} \cdot \mathcal{H}''$. If $\mathcal{H}''$ solves (40) and $\mathcal{H}'$ is a solution of (41) then $\dot{\mathcal{H}} = \mathcal{H}' + \mathcal{H}''$ clearly solves (39b), (39d), (39e), (39f).

From (39), (41), the new form of the reduced linearized problem with unknowns $(\dot{U}, \mathcal{H}')$ reads

$$\widehat{A}_0 \partial_t \dot{U} + \sum_{j=1}^{3} \widehat{A}_j \partial_j \dot{U} + \widehat{\mathcal{C}} \dot{U} = f \qquad \text{in } Q_T^+, \tag{42a}$$

$$\nabla \times \mathfrak{H}' = 0, \quad \text{div } \mathfrak{h}' = 0 \qquad \text{in } Q_T^-, \tag{42b}$$

$$\partial_t \varphi = \dot{v}_N - \hat{v}_2 \partial_2 \varphi - \hat{v}_3 \partial_3 \varphi + \varphi \, \partial_1 \hat{v}_N + g_1, \tag{42c}$$

$$\dot{q} = \widehat{\mathcal{H}} \cdot \mathcal{H}' - [\partial_1 \hat{q}] \varphi + g_2', \tag{42d}$$

$$\mathcal{H}'_N = \partial_2(\widehat{H}_2 \varphi) + \partial_3(\widehat{H}_3 \varphi) \qquad \text{on } \omega_T, \tag{42e}$$

$$\dot{v}_1 = g_4 \quad \text{on } \omega_T^+, \qquad \nu \times \mathcal{H}' = 0 \quad \text{on } \omega_T^-, \tag{42f}$$

$$(\dot{U}, \mathcal{H}', \varphi) = 0 \qquad \text{for } t < 0. \tag{42g}$$

4.4. **Reduction to homogeneous constraints in the "plasma part".** From problem (42) we can deduce nonhomogeneous equations associated with the divergence constraint div $\dot{h} = 0$ and the "redundant" boundary conditions $\dot{H}_N|_{x_1=0} = 0, \dot{H}_1|_{x_1=1} = 0$ for the nonlinear problem. More precisely, with reference to [36, Proposition 2] for the proof, we have the following.

**Proposition 7** ([36]). *Let the basic state* (28) *satisfies assumptions* (29)–(34). *Then solutions of problem* (42) *satisfy*

$$\text{div } \dot{h} = r \quad \text{in } Q_T^+, \tag{43}$$

$$\widehat{H}_2 \partial_2 \varphi + \widehat{H}_3 \partial_3 \varphi - \dot{H}_N - \varphi \, \partial_1 \widehat{H}_N = R \quad \text{on } \omega_T, \qquad \dot{H}_1 = R^+ \quad \text{on } \omega_T^+. \tag{44}$$

*Here*

$$\dot{h} = (\dot{H}_N, \dot{H}_2 \partial_1 \widehat{\Phi}_1, \dot{H}_3 \partial_1 \widehat{\Phi}_1), \quad \dot{H}_N = \dot{H}_1 - \dot{H}_2 \partial_2 \widehat{\Psi} - \dot{H}_3 \partial_3 \widehat{\Psi}.$$

*The functions* $r = r(t, x)$, $R = R(t, x')$ *and* $R^+ = R^+(t, x')$, *which vanish in the past, are determined by the source terms and the basic state as solutions to the linear inhomogeneous equations*

$$\begin{cases} \partial_t a + \frac{1}{\partial_1 \widehat{\Phi}_1} \{ \hat{w} \cdot \nabla a + a \, \text{div } \hat{u} \} = \mathcal{F}_H & \text{in } Q_T^+, \\ \partial_t R + \partial_2(\hat{v}_2 R) + \partial_3(\hat{v}_3 R) = \mathcal{Q} & \text{on } \omega_T, \\ \partial_t R^+ + \partial_2(\hat{v}_2 R^+) + \partial_3(\hat{v}_3 R^+) = \mathcal{Q}^+ & \text{on } \omega_T^+, \end{cases} \tag{45}$$

*where* $a = r/\partial_1 \widehat{\Phi}_1$, $\mathcal{F}_H = (\text{div } f_H)/\partial_1 \widehat{\Phi}_1$,

$$f_H = (f_N, f_6, f_7), \quad f_N = f_5 - f_6 \partial_2 \widehat{\Psi} - f_7 \partial_3 \widehat{\Psi},$$

$$\mathcal{Q} = \{ \partial_2(\widehat{H}_2 g_1) + \partial_3(\widehat{H}_3 g_1) - f_N \}|_{x_1=0}, \qquad \mathcal{Q}^+ = \{ \partial_2(\widehat{H}_2 g_4) + \partial_3(\widehat{H}_3 g_4) + f_5 \}|_{x_1=1}.$$

Let us reduce (42) to a problem with homogeneous boundary conditions (42c), (42d), (42f) (i.e. $g_1 = g_2' = g_4 = 0$) and homogeneous constraints (43) and (44) (i.e. $r = R = R^+ = 0$). More precisely, we describe a "lifting" function as follows:

$$\widetilde{U} = (\tilde{q}, \tilde{v}_1, 0, 0, \widetilde{H}, 0),$$

where $\tilde{q} = g_2', \tilde{v}_1 = -g_1$ on $\omega_T, \tilde{v}_1 = g_4$ on $\omega_T^+$, and where $\widetilde{H}$ solves the equation for $\dot{H}$ contained in (42a) with $\dot{v} = 0$:

$$\partial_t \widetilde{H} + \frac{1}{\partial_1 \widehat{\Phi}_1} \left\{ (\hat{w} \cdot \nabla) \widetilde{H} - (\tilde{h} \cdot \nabla) \hat{v} + \widetilde{H} \text{div } \hat{w} \right\} = f_H \qquad \text{in } Q_T^+, \tag{46}$$

where $\tilde{h} = (\widetilde{H}_1 - \widetilde{H}_2 \partial_2 \widehat{\Psi} - \widetilde{H}_3 \partial_3 \widehat{\Psi}, \widetilde{H}_2, \widetilde{H}_3)$, $f_H = (f_5, f_6, f_7)$. It is very important that, in view of (33), we have $\hat{w}_1|_{x_1=0} = \hat{w}_1|_{x_1=1} = 0$; therefore the linear equation (46) does not need any boundary condition and we easily get the estimates

$$\|\widetilde{H}\|_{H^k_{tan,\gamma}(Q_T^+)} \leq C \|f\|_{H^k_{tan,\gamma}(Q_T^+)}, \qquad k = 1, 2. \tag{47}$$

Here and after $C$ is a constant that can change from line to line, and sometimes we show the dependence of $C$ from other constants. In particular, in (47) the constant $C$ depends on $K$ and $T$. From (47) we obtain

$$\|\widetilde{U}\|_{H^1_{tan,\gamma}(Q_T^+)} \leq C(\|f\|_{H^1_{tan,\gamma}(Q_T^+)} + \|g_1, g_2'\|_{H^{1/2}_\gamma(\omega_T)} + \|g_4\|_{H^{1/2}_\gamma(\omega_T^+)})$$
$$\leq \frac{C}{\gamma}(\|f\|_{H^2_{tan,\gamma}(Q_T^+)} + \|g_1, g_2'\|_{H^{3/2}_\gamma(\omega_T)} + \|g_4\|_{H^{3/2}_\gamma(\omega_T^+)}). \quad (48)$$

Then the new unknown

$$U^\natural = \dot{U} - \widetilde{U}, \quad \mathcal{H}^\natural = \mathcal{H}' \quad (49)$$

satisfies problem (42) with $f = F$, where

$$F = (F_1, \dots, F_8) = f - \mathbb{P}'_e(\widehat{U}, \widehat{\Psi})\widetilde{U}. \quad (50)$$

In view of (46), $(F_5, F_6, F_7) = 0$, and it follows from Proposition 7 that $U^\natural$ satisfies (43) and (44) with $r = R = R^+ = 0$. Moreover, again taking into account (47), for the new source term $F$ we get the estimate

$$\|F\|_{H^1_{tan,\gamma}(Q_T^+)} \leq C\{\|f\|_{H^1_{tan,\gamma}(Q_T^+)} + \|\widetilde{H}\|_{H^2_{tan,\gamma}(Q_T^+)} + \|(\tilde{q}, \tilde{v}_1)\|_{H^2_\gamma(Q_T^+)}\}$$
$$\leq C\{\|f\|_{H^2_{tan,\gamma}(Q_T^+)} + \|g_1, g_2'\|_{H^{3/2}_\gamma(\omega_T)} + \|g_4\|_{H^{3/2}_\gamma(\omega_T^+)}\}. \quad (51)$$

Concerning the above inequality it is worth noticing that in $\mathbb{P}'_e(\widehat{U}, \widehat{\Psi})\widetilde{U}$ the only normal derivatives involved are those of the noncharacteristic part of $\widetilde{U}$, namely $\tilde{q}, \tilde{v}_1$; in this regard see the equivalent system (58) that will be introduced below.

Dropping for convenience the indices $\natural$ in (49), the new form of our reduced linearized problem now reads

$$\widehat{A}_0 \partial_t U + \sum_{j=1}^3 \widehat{A}_j \partial_j U + \widehat{C} U = F \quad \text{in } Q_T^+, \quad (52a)$$

$$\nabla \times \mathfrak{H} = 0, \quad \text{div } \mathfrak{h} = 0 \quad \text{in } Q_T^-, \quad (52b)$$

$$\partial_t \varphi = v_N - \hat{v}_2 \partial_2 \varphi - \hat{v}_3 \partial_3 \varphi + \varphi \, \partial_1 \hat{v}_N, \quad (52c)$$

$$q = \widehat{\mathcal{H}} \cdot \mathcal{H} - [\partial_1 \hat{q}] \varphi, \quad (52d)$$

$$\mathcal{H}_N = \partial_2(\widehat{\mathcal{H}}_2 \varphi) + \partial_3(\widehat{\mathcal{H}}_3 \varphi) \quad \text{on } \omega_T, \quad (52e)$$

$$v_1 = 0 \quad \text{on } \omega_T^+, \quad \nu \times \mathcal{H} = 0 \quad \text{on } \omega_T^-, \quad (52f)$$

$$(U, \mathcal{H}, \varphi) = 0 \quad \text{for } t < 0. \quad (52g)$$

and solutions should satisfy

$$\text{div } h = 0 \quad \text{in } Q_T^+, \quad (53)$$

$$H_N = \widehat{H}_2 \partial_2 \varphi + \widehat{H}_3 \partial_3 \varphi - \varphi \, \partial_1 \widehat{H}_N \quad \text{on } \omega_T, \quad H_1 = 0 \quad \text{on } \omega_T^+. \quad (54)$$

All the notations here for $U$ and $\mathcal{H}$ (e.g., $h$, $\mathfrak{H}$, $\mathfrak{h}$, etc.) are analogous to the corresponding ones for $\dot{U}$ and $\dot{\mathcal{H}}$ introduced above.

4.5. **An equivalent formulation of** (52). In the following analysis it is convenient to make use of different "plasma" variables and an equivalent form of equations (52a). We define the matrix

$$\hat{\eta} = \begin{pmatrix} 1 & -\partial_2 \widehat{\Psi} & -\partial_3 \widehat{\Psi} \\ 0 & \partial_1 \widehat{\Phi}_1 & 0 \\ 0 & 0 & \partial_1 \widehat{\Phi}_1 \end{pmatrix}. \quad (55)$$

It follows that

$$u = (v_N, v_2 \partial_1 \widehat{\Phi}_1, v_3 \partial_1 \widehat{\Phi}_1) = \hat{\eta} \, v, \quad h = (H_N, H_2 \partial_1 \widehat{\Phi}_1, H_3 \partial_1 \widehat{\Phi}_1) = \hat{\eta} \, H. \quad (56)$$

Multiplying (52a) on the left side by the matrix

$$\widehat{R} = \begin{pmatrix} 1 & \underline{0} & \underline{0} & 0 \\ \underline{0}^T & \hat{\eta} & 0_3 & \underline{0}^T \\ \underline{0}^T & 0_3 & \hat{\eta} & \underline{0}^T \\ 0 & \underline{0}^T & \underline{0}^T & 1 \end{pmatrix},$$

after some calculations we get the symmetric hyperbolic system for the new vector of unknowns $\mathcal{U} = (q, u, h, S)$ (compare with (6), (52a)):

$$\partial_1 \widehat{\Phi}_1 \begin{pmatrix} \hat{\rho}_p/\hat{\rho} & \underline{0} & -(\hat{\rho}_p/\hat{\rho})\hat{h} & 0 \\ \underline{0}^T & \hat{\rho}\hat{a}_0 & 0_3 & \underline{0}^T \\ -(\hat{\rho}_p/\hat{\rho})\hat{h}^T & 0_3 & \hat{a}_0 + (\hat{\rho}_p/\hat{\rho})\hat{h} \otimes \hat{h} & \underline{0}^T \\ 0 & \underline{0} & \underline{0} & 1 \end{pmatrix} \partial_t \begin{pmatrix} q \\ u \\ h \\ S \end{pmatrix} + \begin{pmatrix} 0 & \nabla\cdot & \underline{0} & 0 \\ \nabla & 0_3 & 0_3 & \underline{0}^T \\ \underline{0}^T & 0_3 & 0_3 & \underline{0}^T \\ 0 & \underline{0} & \underline{0} & 0 \end{pmatrix} \begin{pmatrix} q \\ u \\ h \\ S \end{pmatrix}$$

$$+ \partial_1 \widehat{\Phi}_1 \begin{pmatrix} (\hat{\rho}_p/\hat{\rho})\hat{w} \cdot \nabla & \nabla\cdot & -(\hat{\rho}_p/\hat{\rho})\hat{h}\hat{w} \cdot \nabla & 0 \\ \nabla & \hat{\rho}\hat{a}_0\hat{w} \cdot \nabla & -\hat{a}_0\hat{h} \cdot \nabla & \underline{0}^T \\ -(\hat{\rho}_p/\hat{\rho})\hat{h}^T\hat{w} \cdot \nabla & -\hat{a}_0\hat{h} \cdot \nabla & (\hat{a}_0 + (\hat{\rho}_p/\hat{\rho})\hat{h} \otimes \hat{h})\hat{w} \cdot \nabla & \underline{0}^T \\ 0 & \underline{0} & \underline{0} & \hat{w} \cdot \nabla \end{pmatrix} \begin{pmatrix} q \\ u \\ h \\ S \end{pmatrix} + \widehat{\mathcal{C}}'\mathcal{U} = \mathcal{F}, \tag{57}$$

where $\hat{\rho} := \rho(\hat{p}, \widehat{S})$, $\hat{\rho}_p := \rho_p(\hat{p}, \widehat{S})$, and $\hat{a}_0$ is the symmetric and positive definite matrix

$$\hat{a}_0 = (\hat{\eta}^{-1})^T \hat{\eta}^{-1},$$

with a new matrix $\widehat{\mathcal{C}}'$ in the zero-order term (whose precise form has no importance) and where we have set $\mathcal{F} = \partial_1 \widehat{\Phi}_1 \, \widehat{R} F$. We write system (57) in compact form as

$$\widehat{\mathcal{A}}_0 \partial_t \mathcal{U} + \sum_{j=1}^{3} (\widehat{\mathcal{A}}_j + \mathcal{E}_{1j+1}) \partial_j \mathcal{U} + \widehat{\mathcal{C}}' \mathcal{U} = \mathcal{F}, \tag{58}$$

where

$$\mathcal{E}_{12} = \begin{pmatrix} 0 & 1 & 0 & 0 & \cdots & 0 \\ 1 & 0 & 0 & 0 & \cdots & 0 \\ 0 & 0 & 0 & 0 & \cdots & 0 \\ 0 & 0 & 0 & 0 & \cdots & 0 \\ \vdots & \vdots & \vdots & \vdots & & \vdots \\ 0 & 0 & 0 & 0 & \cdots & 0 \end{pmatrix}, \qquad \mathcal{E}_{13} = \begin{pmatrix} 0 & 0 & 1 & 0 & \cdots & 0 \\ 0 & 0 & 0 & 0 & \cdots & 0 \\ 1 & 0 & 0 & 0 & \cdots & 0 \\ 0 & 0 & 0 & 0 & \cdots & 0 \\ \vdots & \vdots & \vdots & \vdots & & \vdots \\ 0 & 0 & 0 & 0 & \cdots & 0 \end{pmatrix},$$

$$\mathcal{E}_{14} = \begin{pmatrix} 0 & 0 & 0 & 1 & \cdots & 0 \\ 0 & 0 & 0 & 0 & \cdots & 0 \\ 0 & 0 & 0 & 0 & \cdots & 0 \\ 1 & 0 & 0 & 0 & \cdots & 0 \\ \vdots & \vdots & \vdots & \vdots & & \vdots \\ 0 & 0 & 0 & 0 & \cdots & 0 \end{pmatrix}.$$

The formulation (58) has the advantage of the form of the boundary matrix of the system $\widehat{\mathcal{A}}_1 + \mathcal{E}_{12}$, with

$$\widehat{\mathcal{A}}_1 = 0 \qquad \text{on } \omega_T \cup \omega_T^+, \tag{59}$$

because $\hat{w}_1 = \hat{h}_1 = 0$, and $\mathcal{E}_{12}$ is a constant matrix. Thus system (58) is symmetric hyperbolic with characteristic boundary of constant multiplicity (see [26, 28, 30] for maximally dissipative boundary conditions). Thus, the final form of our reduced linearized problem is

$$\widehat{\mathcal{A}}_0 \partial_t \mathcal{U} + \sum_{j=1}^{3} (\widehat{\mathcal{A}}_j + \mathcal{E}_{1j+1}) \partial_j \mathcal{U} + \widehat{\mathcal{C}} \mathcal{U} = \mathcal{F}, \qquad \text{in } Q_T^+, \tag{60a}$$

$$\nabla \times \mathfrak{H} = 0, \quad \text{div } \mathfrak{h} = 0 \qquad \text{in } Q_T^-, \tag{60b}$$

$$\partial_t \varphi = u_1 - \hat{v}_2 \partial_2 \varphi - \hat{v}_3 \partial_3 \varphi + \varphi \, \partial_1 \hat{v}_N, \tag{60c}$$

$$q = \widehat{\mathcal{H}} \cdot \mathcal{H} - [\partial_1 \hat{q}] \varphi, \tag{60d}$$

$$\mathcal{H}_N = \partial_2 (\widehat{\mathcal{H}}_2 \varphi) + \partial_3 (\widehat{\mathcal{H}}_3 \varphi) \qquad \text{on } \omega_T, \tag{60e}$$

$$v_1 = 0 \quad \text{on } \omega_T^+, \qquad \nu \times \mathcal{H} = 0 \quad \text{on } \omega_T^-, \tag{60f}$$

$$(\mathcal{U}, \mathcal{H}, \varphi) = 0 \qquad \text{for } t < 0, \tag{60g}$$

under the constraints (53), (54).

## 5. The main result of [35]

We recall the main result of the paper [35]. Actually, in [35] we considered the case with the whole space domain, i.e. the case $\Omega^\pm = \mathbb{R}^3 \cap \{x_1 \gtrless 0\}$ and $\Gamma = \mathbb{R}^3 \cap \{x_1 = 0\}$. The result of [35] can be readily extended to the present space domain with minor modifications (see Appendix C), in particular, when treating the added fixed top boundary $\Gamma_+$, under the standard rigid wall boundary conditions on $\Gamma_+$ in (14), see [29], [33].

Recall that $\mathcal{U} = (q, u, h, S)$, where $u$ and $h$ were defined in (56).

**Theorem 8.** *Let $T > 0$. Let the basic state (28) satisfies assumptions (29)–(34) and*

$$|\widehat{H} \times \widehat{\mathcal{H}}| \geq \delta_0/2 > 0 \qquad on \; \omega_T, \tag{61}$$

*where $\delta_0$ is a fixed constant. There exists $\gamma_0 \geq 1$ such that for all $\gamma \geq \gamma_0$ and for all $\mathcal{F}_\gamma \in H^1_{tan,\gamma}(Q_T^+)$, vanishing in the past, namely for $t < 0$, problem (60) has a unique solution $(\mathcal{U}_\gamma, \mathcal{H}_\gamma, \varphi_\gamma) \in H^1_{tan,\gamma}(Q_T^+) \times H^1_\gamma(Q_T^-) \times H^{3/2}_\gamma(\omega_T)$ with traces $(q_\gamma, u_{1\gamma}, h_{1\gamma})|_{\omega_T} \in H^{1/2}_\gamma(\omega_T)$, $\mathcal{H}_\gamma|_{\omega_T} \in H^{1/2}_\gamma(\omega_T)$, and $(q_\gamma, u_{1\gamma}, h_{1\gamma})|_{\omega_T^\pm} \in H^{1/2}_\gamma(\omega_T^\pm)$. Moreover, the solution obeys the a priori estimate*

$$\gamma \left( \|\mathcal{U}_\gamma\|^2_{H^1_{tan,\gamma}(Q_T^+)} + \|\mathcal{H}_\gamma\|^2_{H^1_\gamma(Q_T^-)} + \|(q_\gamma, u_{1\gamma}, h_{1\gamma})\|^2_{\omega_T \cup \omega_T^+ \, H^{1/2}_\gamma(\omega_T \cup \omega_T^+)} + \|\varphi_\gamma\|^2_{H^{3/2}_\gamma(\omega_T)} \right)$$
$$\leq \frac{C}{\gamma} \|\mathcal{F}_\gamma\|^2_{H^1_{tan,\gamma}(Q_T^+)}, \tag{62}$$

*where we have set $\mathcal{U}_\gamma = e^{-\gamma t} \mathcal{U}, \mathcal{H}_\gamma = e^{-\gamma t} \mathcal{H}, \varphi_\gamma = e^{-\gamma t} \varphi$ and so on. Here $C = C(K, T, \delta_0) > 0$ is a constant independent of the data $\mathcal{F}$ and $\gamma$.*

We observe that in the above a priori estimate there is no loss of regularity from the data $\mathcal{F}$ to the solution $(\mathcal{U}, \mathcal{H})$.

**Remark 9.** *Strictly speaking, the uniqueness of the solution to problem (60) follows from (62), provided that our solution belongs to $H^2$. The existence of solutions with a higher degree of regularity (in particular, $H^2$) is given in Theorem 16.*

**Remark 10.** *Differently from the above statement, in [35] it is proved that $\varphi_\gamma \in H^1_\gamma(\omega_T)$ with corresponding a priori estimate. Actually, the regularity of $\varphi_\gamma$ can be easily improved to $H^{3/2}_\gamma(\omega_T)$ as above, because, under the stability condition (61), $\nabla_{t,x'} \varphi_\gamma$ is estimated by the traces $(u_{1\gamma}, h_{1\gamma}, \mathfrak{h}_{1\gamma})|_{\omega_T} \in H^{1/2}_\gamma(\omega_T)$.*

## 6. The elliptic problem (40)

Let us first study the elliptic system (40) for freezed time $t$. We rewrite it as (we drop for convenience the $''$ of $\mathcal{H}''$, $\mathfrak{H}''$)

$$\begin{cases} \nabla \times \mathfrak{H} = \chi, \quad \operatorname{div}(A\mathfrak{H}) = \Xi & \text{in } \Omega^-, \\ (A\mathfrak{H})_1 = g_3 \quad \text{on } \Gamma, \quad \nu \times \mathfrak{H} = g_5 \quad \text{on } \Gamma_-, \\ (x_2, x_3) \to \mathfrak{H}(t, x_1, x_2, x_3) & \text{is 1-periodic.} \end{cases} \tag{63}$$

because $\mathfrak{h} = A\mathfrak{H}$, and where $A = A(\nabla\widehat{\Psi}) = (\partial_1\widehat{\Phi}_1)^{-1}\hat{\eta}\,\hat{\eta}^T$, with $\hat{\eta}$ defined in (55). The matrix $A$ is symmetric and positive definite. On $\Gamma_-$, since $\widehat{\Psi} = 0$ one has $\nu \times \mathcal{H} = \nu \times \mathfrak{H}$. We denote by $\nu = (\pm 1, 0, 0)$ the outward normal vector to $\Gamma_\pm$.

6.1. **Preliminaries.** Only in this subsection, for the sake of clarity the spaces of vector fields are indicated with the complete notations $L^2(\Omega^-; \mathbb{R}^3)$, $H^1(\Omega^-; \mathbb{R}^3)$, and so on.

Let us introduce the space of tangential $H^1$ vector fields on $\Gamma$,

$$H^1_{\tau\Gamma}(\Omega^-; \mathbb{R}^3) := \left\{ \eta \in H^1(\Omega^-; \mathbb{R}^3) : \ \nu \times \eta = 0 \text{ on } \Gamma \right\}.$$

We also introduce the space of $H^1$ scalar functions vanishing on $\Gamma_-$:

$$H^1_{0\Gamma_-}(\Omega^-) := \left\{ \phi \in H^1(\Omega^-) : \ \phi = 0 \text{ on } \Gamma_- \right\}, \qquad \|\phi\|_{H^1_{0\Gamma_-}} := \|\nabla\phi\|_{L^2(\Omega^-)}.$$

The definition of the norm is possible because the Poincaré inequality applies in $H^1_{0\Gamma_-}(\Omega^-)$; it also yields $\|\phi_{|\Gamma}\|_{H^{1/2}(\Gamma)} \le C\|\phi\|_{H^1_{0\Gamma_-}}$. Let us denote

$$\mathrm{Curl}_{\tau\Gamma}(\Omega^-) := \left\{ \nabla \times \eta : \ \eta \in H^1_{\tau\Gamma}(\Omega^-; \mathbb{R}^3) \right\},$$

$$G_{0\Gamma_-}(\Omega^-) := \left\{ \nabla\phi : \ \phi \in H^1_{0\Gamma_-}(\Omega^-) \right\}.$$

**Proposition 11.** *The following orthogonal decomposition holds:*

$$L^2(\Omega^-; \mathbb{R}^3) = \mathrm{Curl}_{\tau\Gamma}(\Omega^-) \oplus G_{0\Gamma_-}(\Omega^-).$$

For more general decompositions we refer the reader to [3].

*Proof.* By an integration by parts one first show that the above subspaces are orthogonal w.r.t. the $L^2$ inner product. Then, assuming there exists $u \in L^2(\Omega^-; \mathbb{R}^3)$ which is orthogonal to $\mathrm{Curl}_{\tau\Gamma}(\Omega^-)$ one obtains

$$\nabla \times u = 0 \quad \text{in } \mathcal{D}'(\Omega^-), \qquad \nu \times u = 0 \quad \text{in } H^{-1/2}(\Gamma_-).$$

If $u$ is also orthogonal to $G_{0\Gamma_-}(\Omega^-)$ one obtains

$$\operatorname{div} u = 0 \quad \text{in } \mathcal{D}'(\Omega^-), \qquad u_1 = 0 \quad \text{in } H^{-1/2}(\Gamma).$$

This shows that $u = 0$ and that the above direct sum span the whole space $L^2(\Omega^-; \mathbb{R}^3)$ as in the statement. $\qquad\square$

From the Proposition, any vector field $v \in L^2(\Omega^-; \mathbb{R}^3)$ can be uniquely decomposed as

$$v = \nabla \times \eta + \nabla\phi, \tag{64}$$

with $\eta \in H^1_{\tau\Gamma}(\Omega^-; \mathbb{R}^3)$, $\nabla\phi \in G_{0\Gamma_-}(\Omega^-)$. However, $\eta$ in (64) is not uniquely defined. To do so, we choose $\eta \in H^1_{\tau\Gamma}(\Omega^-; \mathbb{R}^3)$ such that, for an assigned $\nabla \times \eta$, it also solves

$$\operatorname{div} \eta = 0 \quad \text{in } \Omega^-, \qquad \eta_1 = 0 \quad \text{on } \Gamma_-.$$

This $\eta$ is uniquely defined; in fact, any such $\eta$ satisfies by an integration by parts and Poincaré inequality applied to each component

$$\int_{\Omega^-} |\nabla \times \eta|^2 \, dx = \int_{\Omega^-} \partial_i\eta_j\partial_i\eta_j \, dx \ge C \int_{\Omega^-} |\eta|^2 \, dx,$$

and uniqueness follows from the linearity of the problem. Summing up we have obtained the following result.

**Proposition 12.** *Given any vector field $v \in L^2(\Omega^-; \mathbb{R}^3)$ there exists a unique $\eta \in H^1_{\tau\Gamma}(\Omega^-; \mathbb{R}^3)$, $\nabla\phi \in G_{0\Gamma_-}(\Omega^-)$ such that* (64) *holds and also*

$$\|\eta\|_{H^1(\Omega^-)} + \|\nabla\phi\|_{L^2(\Omega^-)} \le C\|v\|_{L^2(\Omega^-)}.$$

We denote by $H^{-1}_{\tau\Gamma} = H^{-1}_{\tau\Gamma}(\Omega^-)$ the dual space of $H^1_{\tau\Gamma}(\Omega^-; \mathbb{R}^3)$. Denoting by $\langle\cdot, \cdot\rangle$ the pairing between $H^{-1}_{\tau\Gamma}$ and $H^1_{\tau\Gamma}$, we have

$$\|\chi\|_{H^{-1}_{\tau\Gamma}} := \sup_{\psi \in H^1_{\tau\Gamma}} \frac{|\langle\chi, \psi\rangle|}{\|\psi\|_{H^1}}, \qquad \|\chi\|_{H^{-1}_{\tau\Gamma}} \le C\|\chi\|_{L^2(\Omega^-)}, \qquad \|\partial_i\chi\|_{H^{-1}_{\tau\Gamma}} \le \|\chi\|_{L^2(\Omega^-)} \quad \text{for } i = 2, 3,$$

where the second inequality follows by an integration by parts.

We denote by $H^{-1}_{0\Gamma_-} = H^{-1}_{0\Gamma_-}(\Omega^-)$ the dual space of $H^1_{0\Gamma_-}(\Omega^-)$. Denoting by $\langle\cdot, \cdot\rangle$ the pairing between $H^{-1}_{0\Gamma_-}$ and $H^1_{0\Gamma_-}$, we have

$$\|\xi\|_{H^{-1}_{0\Gamma_-}} := \sup_{\phi \in H^1_{0\Gamma_-}} \frac{|\langle\xi, \phi\rangle|}{\|\phi\|_{H^1_{0\Gamma_-}}}, \quad \|\xi\|_{H^{-1}_{0\Gamma_-}} \le C\|\xi\|_{L^2(\Omega^-)}, \quad \|\partial_i\xi\|_{H^{-1}_{0\Gamma_-}} \le \|\xi\|_{L^2(\Omega^-)} \quad \text{for } i = 2, 3, \tag{65}$$

where the first inequality follows from the Poincaré inequality, and the second one by integrating by parts. The use of negative norms will be crucial in section 8.2 in the analysis of some commutators.

6.2. **Compatibility conditions.** For the resolution of (63) some necessary compatibility conditions are needed:

$$g_5 \cdot \nu = 0 \quad \text{on } \Gamma_-, \tag{66a}$$

$$\int_{\Omega^-} \chi \cdot \eta \, dx = \int_{\Gamma_-} g_5 \cdot \eta \, dx', \qquad \forall \eta \in H^1_{\tau\Gamma}(\Omega^-; \mathbb{R}^3) \text{ such that } \nabla\times\eta = 0 \text{ in } \Omega^-. \tag{66b}$$

The first equation on $\Gamma_-$ follows from $(\nu\times\mathfrak{H})\cdot\nu = g_5\cdot\nu = 0$. (66b) follows by multiplying the first equation of (63) by $\eta$ as above and integrating by parts. Observe that one could choose $\eta = \nabla\phi$, $\phi \in H^1_0(\Omega^-)$. In such a case from (66) one gets

$$\int_{\Omega^-} \chi \cdot \nabla\phi \, dx = 0 \qquad \forall\phi \in H^1_0(\Omega^-),$$

i.e. the weak form of the natural constraint div $\chi = 0$.

We have the following result.

**Theorem 13.** *Assume that for each fixed $t$ the data $(\chi, \Xi, g_3, g_5)$ in* (63) *satisfy $(\chi, \Xi) \in L^2(\Omega^-)$, $g_3 \in H^{1/2}(\Gamma)$, $g_5 \in H^{1/2}(\Gamma_-)$ and the compatibility conditions* (66). *Then there exists a unique solution $\mathfrak{H} \in H^1(\Omega^-)$ of* (63) *and*

$$\|\mathfrak{H}\|_{L^2(\Omega^-)} \le C(K)\big(\|\chi\|_{H^{-1}_{\tau\Gamma}} + \|\Xi\|_{H^{-1}_{0\Gamma_-}} + \|g_3\|_{H^{-1/2}(\Gamma)} + \|g_5\|_{H^{-1/2}(\Gamma_-)}\big), \tag{67a}$$

$$\|\nabla\mathfrak{H}\|_{L^2(\Omega^-)} \le C(K)\big(\|\chi, \Xi\|_{L^2(\Omega^-)} + \|g_3\|_{H^{1/2}(\Gamma)} + \|g_5\|_{H^{1/2}(\Gamma_-)}\big). \tag{67b}$$

*If $(\chi, \Xi) \in H^1(\Omega^-)$, $g_3 \in H^{3/2}(\Gamma)$, $g_5 \in H^{3/2}(\Gamma_-)$, then $\mathfrak{H} \in H^2(\Omega^-)$ and*

$$\|\mathfrak{H}\|_{H^2(\Omega^-)} \le C(K)\big(\|\chi, \Xi\|_{H^1(\Omega^-)} + \|g_3\|_{H^{3/2}(\Gamma)} + \|g_5\|_{H^{3/2}(\Gamma_-)}\big). \tag{68}$$

Clearly, this statement can be translated in a similar one for the original variable $\mathcal{H}$ in place of $\mathfrak{H}$.

*Proof.* (1) Given $\chi$, let us consider the elliptic system

$$\begin{cases} \nabla\times\zeta = \chi, \\ \text{div}\,\zeta = 0 & \text{in } \Omega^-, \\ \zeta_1 = 0 & \text{on } \Gamma, \\ \nu\times\zeta = g_5 & \text{on } \Gamma_-, \\ (x_2, x_3) \to \zeta(t, x_1, x_2, x_3) & \text{is 1-periodic.} \end{cases} \tag{69}$$

We first show by integrating by parts that a solution $\zeta$ satisfies

$$\int_{\Omega^-} \chi \cdot \eta \, dx = \int_{\Omega^-} \nabla \times \zeta \cdot \eta \, dx = \int_{\Omega^-} \zeta \cdot \nabla \times \eta \, dx + \int_{\Gamma_-} g_5 \cdot \eta \, dx' \qquad \forall \eta \in H^1_{\tau\Gamma}(\Omega^-),$$

$$0 = \int_{\Omega^-} \operatorname{div} \zeta \, \phi \, dx = -\int_{\Omega^-} \zeta \cdot \nabla \phi \, dx \qquad \forall \phi \in H^1(\Omega^-), \, \phi_{|\Gamma_-} = 0.$$

It yields

$$\int_{\Omega^-} \chi \cdot \eta \, dx = \int_{\Omega^-} \zeta \cdot (\nabla \times \eta + \nabla \phi) \, dx + \int_{\Gamma_-} g_5 \cdot \eta \, dx', \tag{70}$$

for all $\eta \in H^1_{\tau\Gamma}(\Omega^-), \nabla \phi \in G_{0\Gamma_-}(\Omega^-)$. Given any vector field $v$, let us choose $\eta \in H^1_{\tau\Gamma}(\Omega^-), \nabla \phi \in G_{0\Gamma_-}(\Omega^-)$ as in Proposition 12 and substitute in (70). This gives the weak formulation of (69)

$$\int_{\Omega^-} \chi \cdot \eta \, dx = \int_{\Omega^-} \zeta \cdot v \, dx + \int_{\Gamma_-} g_5 \cdot \eta \, dx', \qquad \forall v \in L^2(\Omega^-). \tag{71}$$

Noticing that

$$\left| \int_{\Omega^-} \chi \cdot \eta \, dx - \int_{\Gamma_-} g_5 \cdot \eta \, dx' \right| \le C \left( \|\chi\|_{H^{-1}_{\tau\Gamma}} + \|g_5\|_{H^{-1/2}(\Gamma_-)} \right) \|\eta\|_{H^1(\Omega^-)}$$

$$\le C \left( \|\chi\|_{H^{-1}_{\tau\Gamma}} + \|g_5\|_{H^{-1/2}(\Gamma_-)} \right) \|v\|_{L^2(\Omega^-)},$$

we may apply the Riesz representation theorem and find a unique solution $\zeta \in L^2(\Omega^-)$ of (71), such that

$$\|\zeta\|_{L^2(\Omega^-)} \le C \left( \|\chi\|_{H^{-1}_{\tau\Gamma}} + \|g_5\|_{H^{-1/2}(\Gamma_-)} \right). \tag{72}$$

Using again (64), (70) and (66) gives that $\zeta$ solves (69) in weak sense. Integrating by parts yields

$$\int_{\Omega^-} |\chi|^2 \, dx = \int_{\Omega^-} |\nabla \times \zeta|^2 \, dx = \int_{\Omega^-} \partial_i \zeta_j \partial_i \zeta_j \, dx + 2 \int_{\Gamma_-} \zeta \cdot \nabla \times g_5 \, dx',$$

where we have set $\nabla \times \mathfrak{J} = (\partial_2 g_{5,3} - \partial_3 g_{5,2}, 0, 0)$, and we infer $\zeta \in H^1(\Omega^-)$ with

$$\|\nabla \zeta\|_{L^2(\Omega^-)} \le C \left( \|\chi\|_{L^2(\Omega^-)} + \|g_5\|_{H^{1/2}(\Gamma_-)} \right). \tag{73}$$

Finally, by elliptic regularization, if $\chi \in H^1(\Omega^-)$ and $g_5 \in H^{3/2}(\Gamma_-)$, then $\zeta \in H^2(\Omega^-)$ and

$$\|\zeta\|_{H^2(\Omega^-)} \le C \left( \|\chi\|_{H^1(\Omega^-)} + \|g_5\|_{H^{3/2}(\Gamma_-)} \right). \tag{74}$$

(2) As $\nabla \times (\mathfrak{H} - \zeta) = 0$, we now look for $\mathfrak{H}$ in the form $\mathfrak{H} = \zeta + \nabla \xi$, where $\xi$ should satisfy the Neumann-Dirichlet problem

$$\begin{cases} \operatorname{div} (A\nabla\xi) = \Xi - \operatorname{div} (A\zeta) & \text{in } \Omega^-, \\ (A\nabla\xi)_1 = g_3 - (A\zeta)_1 & \text{on } \Gamma, \\ \xi = 0 & \text{on } \Gamma_-, \\ (x_2, x_3) \to \xi(t, x_1, x_2, x_3) & \text{is 1-periodic.} \end{cases} \tag{75}$$

We first look for a weak solution $\xi \in H^1_{0\Gamma_-}(\Omega^-)$. Multiplying (75) by $\phi \in H^1_{0\Gamma_-}(\Omega^-)$ and integrating twice by parts yields

$$\int_{\Omega^-} A\nabla\xi \cdot \nabla\phi \, dx = -\int_{\Omega^-} \Xi \phi \, dx - \int_{\Omega^-} A\zeta \cdot \nabla\phi \, dx + \int_{\Gamma} g_3 \phi \, dx', \quad \forall \phi \in H^1_{0\Gamma_-}(\Omega^-). \tag{76}$$

Recalling that the matrix $A$ is positive definite and that Poincaré inequality holds in $H^1_{0\Gamma_-}(\Omega^-)$, and noticing that the right-hand side of (76) is easily estimated by the right-hand side of (67a) times $\|\phi\|_{H^1_{0\Gamma_-}}$, by Lax-Milgram lemma we get the existence of a unique solution $\xi \in H^1_{0\Gamma_-}(\Omega^-)$ in the sense of (76), with

$$\|\nabla \xi\|_{L^2(\Omega^-)} \le C(K) \big( \|\chi\|_{H^{-1}_{\tau\Gamma}} + \|\Xi\|_{H^{-1}_{0\Gamma_-}} + \|g_3\|_{H^{-1/2}(\Gamma)} + \|g_5\|_{H^{-1/2}(\Gamma_-)} \big). \tag{77}$$

Then by standard elliptic regularity results

$$\|\nabla \xi\|_{H^1(\Omega^-)} \le C(K) \big( \|\chi, \Xi\|_{L^2(\Omega^-)} + \|g_3\|_{H^{1/2}(\Gamma)} + \|g_5\|_{H^{1/2}(\Gamma_-)} \big). \tag{78}$$

Thus we get the existence of a unique solution of (75) in $H^2(\Omega^-)$, and $\mathfrak{H} = \zeta + \nabla\xi \in H^1(\Omega^-)$ is a solution of (63). By (72), (73), (77) and (78), the solution constructed above satisfies (67). Again by elliptic regularization, if the data of (75) have one more derivative then

$$\|\nabla\xi\|_{H^2(\Omega^-)} \leq C(K)\big(\|\chi, \Xi\|_{H^1(\Omega^-)} + \|g_3\|_{H^{3/2}(\Gamma)} + \|g_5\|_{H^{3/2}(\Gamma_-)}\big). \tag{79}$$

From (74), (79) we get (68). Due to the linearity of the problem, the uniqueness of the solution to (63) follows by showing that the homogeneous problem has only the trivial solution. In fact, in the simply connected domain $\Omega^-$, $\nabla\times\mathfrak{H} = 0$ yields $\mathfrak{H} = \nabla\xi$, and the other equations of (63) give (75) with zero data, apart from an arbitrary constant as a boundary datum on $\Gamma_-$. As this problem has a unique solution $\xi$, this solution is necessarily the arbitrary constant given on $\Gamma_-$, which yields $\mathfrak{H} = \nabla\xi = 0$. □

Now we study the elliptic system (40) taking in account the time dependence. We set $\mathfrak{H}_\gamma = e^{-\gamma t}\mathfrak{H}, \chi_\gamma = e^{-\gamma t}\chi$, and so on, for all functions in (40).

**Lemma 14.** *Assume that the data $(\chi, \Xi, g_3, g_5)$ in (63) satisfy $(\chi_\gamma, \Xi_\gamma) \in L^2(Q_T^-)$, $\partial_t\chi_\gamma \in L^2(-\infty, T; H_{\tau\Gamma}^{-1})$ and $\partial_t\Xi_\gamma \in L^2(-\infty, T; H_{0\Gamma_-}^{-1})$, $g_{3\gamma} \in L^2(-\infty, T; H_\gamma^{1/2}(\Gamma))$ with $\partial_t g_{3\gamma} \in L^2(-\infty, T; H_\gamma^{-1/2}(\Gamma))$, $g_{5\gamma} \in L^2(-\infty, T; H_\gamma^{1/2}(\Gamma_-))$ with $\partial_t g_{5\gamma} \in L^2(-\infty, T; H_\gamma^{-1/2}(\Gamma_-))$, and the compatibility conditions (66). Then the solution $\mathfrak{H}$ of (63) satisfies $\mathfrak{H}_\gamma \in H_\gamma^1(Q_T^-)$ with*

$$\|\mathfrak{H}_\gamma\|_{H_\gamma^1(Q_T^-)} \leq C(K)\big(\|\chi_\gamma, \Xi_\gamma\|_{L^2(Q_T^-)} + \|\partial_t\chi_\gamma\|_{L^2(-\infty, T; H_{\tau\Gamma}^{-1})} + \|\partial_t\Xi_\gamma\|_{L^2(-\infty, T; H_{0\Gamma_-}^{-1})}$$
$$+ \|g_{3\gamma}\|_{L^2(-\infty, T; H_\gamma^{1/2}(\Gamma))} + \|g_{5\gamma}\|_{L^2(-\infty, T; H_\gamma^{1/2}(\Gamma_-))}$$
$$+ \|\partial_t g_{3\gamma}\|_{L^2(-\infty, T; H_\gamma^{-1/2}(\Gamma))} + \|\partial_t g_{5\gamma}\|_{L^2(-\infty, T; H_\gamma^{-1/2}(\Gamma_-))}\big). \tag{80}$$

*If $\chi_\gamma \in H_\gamma^1(Q_T^-) \cap H_\gamma^2(-\infty, T; H_{\tau\Gamma}^{-1})$, $\Xi_\gamma \in H_\gamma^1(Q_T^-) \cap H_\gamma^2(-\infty, T; H_{0\Gamma_-}^{-1})$, $g_{3\gamma} \in H_\gamma^{3/2}(\omega_T) \cap H_\gamma^2(-\infty, T; H_\gamma^{-1/2}(\Gamma))$, $g_{5\gamma} \in H_\gamma^{3/2}(\omega_T^-) \cap H_\gamma^2(-\infty, T; H_\gamma^{-1/2}(\Gamma_-))$ then $\mathfrak{H}_\gamma \in H_\gamma^2(Q_T^-)$ and*

$$\|\mathfrak{H}_\gamma\|_{H_\gamma^2(Q_T^-)} \leq C(K)\big(\gamma\|\chi_\gamma, \Xi_\gamma\|_{H_\gamma^1(Q_T^-)} + \|\chi_\gamma\|_{H_\gamma^2(-\infty, T; H_{\tau\Gamma}^{-1})} + \|\Xi_\gamma\|_{H_\gamma^2(-\infty, T; H_{0\Gamma_-}^{-1})}$$
$$+ \|g_{3\gamma}\|_{H_\gamma^{3/2}(\omega_T)} + \|g_{3\gamma}\|_{H_\gamma^2(-\infty, T; H_\gamma^{-1/2}(\Gamma))}$$
$$+ \|g_{5\gamma}\|_{H_\gamma^{3/2}(\omega_T^-)} + \|g_{5\gamma}\|_{H_\gamma^2(-\infty, T; H_\gamma^{-1/2}(\Gamma_-))}\big). \tag{81}$$

*Proof.* Multiplying the equations in (63) by $e^{-\gamma t}$, applying (67) and integrating in time yields

$$\|\mathfrak{H}_\gamma\|_{L^2(-\infty, T; H_\gamma^1(\Omega^-))} \leq C\big(\gamma\|\chi_\gamma, \Xi_\gamma\|_{L^2(Q_T^-)} + \|g_{3\gamma}\|_{L^2(-\infty, T; H_\gamma^{1/2}(\Gamma))} + \|g_{5\gamma}\|_{L^2(-\infty, T; H_\gamma^{1/2}(\Gamma_-))}\big). \tag{82}$$

For the estimate of $\partial_t\mathfrak{H}_\gamma$ we consider the decomposition $\partial_t\mathfrak{H}_\gamma = \partial_t\zeta_\gamma + \partial_t\nabla\xi_\gamma$, with $\zeta, \xi$ from (69), (75). By taking the time derivative of (69) we get from (72)

$$\|\partial_t\zeta_\gamma\|_{L^2(Q_T^-)} \leq C\big(\|\partial_t\chi_\gamma\|_{L^2(-\infty, T; H_{\tau\Gamma}^{-1})} + \|\partial_t g_{5\gamma}\|_{L^2(-\infty, T; H_\gamma^{-1/2}(\Gamma_-))}\big). \tag{83}$$

The estimate of $\partial_t\nabla\xi_\gamma$ is more complex. Multiplying (75) by $e^{-\gamma t}$ and differentiating in time gives the system

$$\begin{cases} \operatorname{div}(A\partial_t\nabla\xi_\gamma) = \partial_t\Xi_\gamma - \operatorname{div}\partial_t(A\zeta_\gamma) - \operatorname{div}((\partial_t A)\nabla\xi_\gamma) & \text{in } \Omega^-, \\ (A\partial_t\nabla\xi_\gamma)_1 = \partial_t g_{3\gamma} - \partial_t(A\zeta_\gamma)_1 - ((\partial_t A)\nabla\xi_\gamma)_1 & \text{on } \Gamma, \\ \partial_t\xi_\gamma = 0 & \text{on } \Gamma_-. \end{cases}$$

Multiplying by $\partial_t\xi_\gamma$ and integrating twice by parts yields

$$\int_{\Omega^-} A\nabla\partial_t\xi_\gamma \cdot \nabla\partial_t\xi_\gamma \, dx = -\int_{\Omega^-} \partial_t\Xi_\gamma\,\partial_t\xi_\gamma \, dx$$
$$-\int_{\Omega^-} [\partial_t(A\zeta_\gamma) + (\partial_t A)\nabla\xi_\gamma] \cdot \nabla\partial_t\xi_\gamma \, dx + \int_\Gamma \partial_t g_3\,\partial_t\xi_\gamma \, dx',$$

which gives

$$\|\partial_t \nabla \xi_\gamma\|_{L^2(Q_T^-)} \leq C \big( \|\chi_\gamma\|_{H^1_\gamma(-\infty,T;H^{-1}_{rT})} + \|\Xi_\gamma\|_{H^1_\gamma(-\infty,T;H^{-1}_{0\Gamma_-})}$$
$$+ \|g_{3\gamma}\|_{H^1_\gamma(-\infty,T;H^{-1/2}_\gamma(\Gamma))} + \|g_{5\gamma}\|_{H^1_\gamma(-\infty,T;H^{-1/2}_\gamma(\Gamma_-))} \big). \quad (84)$$

Adding (82), (83), (84) gives (80). The proof of (81) is similar, and so we omit the details. □

## 7. FINAL $H^1$ ESTIMATE FOR THE LINEARIZED PROBLEM (39)

Summarizing the results of Theorem 8 and Theorem 13 we get what follows. Let us recall the linearized problem (39)

$$\begin{cases} \widehat{A}_0 \partial_t \dot{U} + \sum_{j=1}^3 \widehat{A}_j \partial_j \dot{U} + \widehat{C} \dot{U} = f & \text{in } Q_T^+, \\ \nabla \times \dot{\mathfrak{H}} = \chi, \quad \text{div } \dot{\mathfrak{h}} = \Xi & \text{in } Q_T^-, \\ \partial_t \varphi = \dot{v}_N - \hat{v}_2 \partial_2 \varphi - \hat{v}_3 \partial_3 \varphi + \varphi \, \partial_1 \hat{v}_N + g_1, \\ \dot{q} = \widehat{\mathcal{H}} \cdot \dot{\mathcal{H}} - [\partial_1 \hat{q}] \varphi + g_2, \\ \dot{\mathcal{H}}_N = \partial_2 (\widehat{\mathcal{H}}_2 \varphi) + \partial_3 (\widehat{\mathcal{H}}_3 \varphi) + g_3 & \text{on } \omega_T, \\ \dot{v}_1 = g_4 \quad \text{on } \omega_T^+, \qquad \nu \times \dot{\mathcal{H}} = g_5 \quad \text{on } \omega_T^-, \\ (\dot{U}, \dot{\mathcal{H}}, \varphi) = 0 & \text{for } t < 0. \end{cases}$$

For the following analysis it seems more convenient to work in the plasma part with the system analogous to (58) and write the vacuum equations in terms of $\dot{\mathfrak{H}}$, as in (63). We find that $\dot{\mathcal{U}} = (\dot{q}, \dot{u}, \dot{h}, \dot{S})$, where

$$\dot{u} = \hat{\eta} \, \dot{v}, \qquad \dot{h} = \hat{\eta} \, \dot{H},$$

see (56), satisfies with $\dot{\mathfrak{H}}$ the system

$$\begin{cases} \widehat{A}_0 \partial_t \dot{\mathcal{U}} + \sum_{j=1}^3 (\widehat{A}_j + \mathcal{E}_{1j+1}) \partial_j \dot{\mathcal{U}} + \widehat{C}' \dot{\mathcal{U}} = \tilde{f} & \text{in } Q_T^+, \\ \nabla \times \dot{\mathfrak{H}} = \chi, \quad \text{div } (A \dot{\mathfrak{H}}) = \Xi & \text{in } Q_T^-, \\ \partial_t \varphi = \dot{u}_1 - \hat{v}_2 \partial_2 \varphi - \hat{v}_3 \partial_3 \varphi + \varphi \, \partial_1 \hat{v}_N + g_1, \\ \dot{q} = \hat{\mathfrak{h}} \cdot \dot{\mathfrak{H}} - [\partial_1 \hat{q}] \varphi + g_2, \\ (A \dot{\mathfrak{H}})_1 = \partial_2 (\widehat{\mathcal{H}}_2 \varphi) + \partial_3 (\widehat{\mathcal{H}}_3 \varphi) + g_3 & \text{on } \omega_T, \\ \dot{v}_1 = g_4 \quad \text{on } \omega_T^+, \qquad \nu \times \dot{\mathcal{H}} = g_5 \quad \text{on } \omega_T^-, \\ (\dot{\mathcal{U}}, \dot{\mathcal{H}}, \varphi) = 0 & \text{for } t < 0, \end{cases} \quad (85)$$

where we have set $\tilde{f} = \partial_1 \widehat{\Phi}_1 \, \widehat{R} f$. Let us remark that (85) is equivalent to the linearized problem (39). Then we have

**Theorem 15.** *Let $T > 0$. Let the basic state (28) satisfy assumptions (29)–(34), (61). There exists $\gamma_1 \geq 1$ such that for all $\gamma \geq \gamma_1$ and for all $\tilde{f}_\gamma \in H^2_\gamma(Q_T^+)$, $\chi_\gamma \in H^1_\gamma(Q_T^-) \cap H^2_\gamma(-\infty,T;H^{-1}_{rT})$, $\Xi_\gamma \in H^1_\gamma(Q_T^-) \cap H^2_\gamma(-\infty,T;H^{-1}_{0\Gamma_-})$, $(g_{1\gamma},g_{2\gamma}) \in H^{3/2}_\gamma(\omega_T)$, $g_{3\gamma} \in H^{3/2}_\gamma(\omega_T) \cap H^2_\gamma(-\infty,T;H^{-1/2}_\gamma(\Gamma))$, $g_{4\gamma} \in H^{3/2}_\gamma(\omega_T^+)$, $g_{5\gamma} \in H^{3/2}_\gamma(\omega_T^-) \cap H^2_\gamma(-\infty,T;H^{-1/2}_\gamma(\Gamma_-))$, with $(\chi, g_5)$ satisfying the compatibility conditions (66), and all functions vanishing in the past, problem (85) has a unique solution $(\dot{\mathcal{U}}_\gamma, \dot{\mathfrak{H}}_\gamma, \varphi_\gamma) \in H^1_{tan,\gamma}(Q_T^+) \times \dot{H}^1_\gamma(Q_T^-) \times H^{3/2}_\gamma(\omega_T)$ with trace $(\dot{q}_\gamma, \dot{u}_{1\gamma}, \dot{h}_{1\gamma})|_{\omega_T} \in H^{1/2}_\gamma(\omega_T)$. Moreover, the solution obeys the a priori estimates*

$$\gamma \left( \|\dot{\mathcal{U}}_\gamma\|^2_{H^1_{tan,\gamma}(Q_T^+)} + \|\dot{\mathfrak{H}}_\gamma\|^2_{H^1_\gamma(Q_T^-)} + \|(\dot{q}_\gamma, \dot{u}_{1\gamma}, \dot{h}_{1\gamma})|_{\omega_T}\|^2_{H^{1/2}_\gamma(\omega_T)} + \|\varphi_\gamma\|^2_{H^{3/2}_\gamma(\omega_T)} \right)$$
$$\leq \frac{C}{\gamma} \big\{ \|\tilde{f}_\gamma\|^2_{H^2_{tan,\gamma}(Q_T^+)} + \gamma \|\chi_\gamma, \Xi_\gamma\|_{H^1_\gamma(Q_T^-)} + \|\chi_\gamma\|_{H^2_\gamma(-\infty,T;H^{-1}_{rT})} + \|\Xi_\gamma\|_{H^2_\gamma(-\infty,T;H^{-1}_{0\Gamma_-})}$$
$$+ \|g_{1\gamma}, g_{2\gamma}, g_{3\gamma}\|^2_{H^{3/2}_\gamma(\omega_T)} + \|g_{3\gamma}\|_{H^2_\gamma(-\infty,T;H^{-1/2}_\gamma(\Gamma))}$$
$$+ \|g_{4\gamma}\|^2_{H^{3/2}_\gamma(\omega_T^+)} + \|g_{5\gamma}\|^2_{H^{3/2}_\gamma(\omega_T^-)} + \|g_{5\gamma}\|_{H^2_\gamma(-\infty,T;H^{-1/2}_\gamma(\Gamma_-))} \big\}, \quad (86)$$

*where we have set* $\dot{\mathcal{U}}_\gamma = e^{-\gamma t}\,\dot{U}, \dot{\mathfrak{H}}_\gamma = e^{-\gamma t}\,\dot{\mathcal{H}}, \varphi_\gamma = e^{-\gamma t}\,\varphi$ *and so on. Here* $C = C(K, T, \delta_0) > 0$ *is a constant independent of the data* $\tilde{f}, \chi, \Xi, g$ *and* $\gamma$.

*Proof.* The proof follows from Theorem 8 and Lemma 14, with (86) following from the inequalities (48), (51), (62), (80), (81). In particular, $\|\dot{\mathfrak{H}}_\gamma\|_{H^1_\gamma(Q_T^-)} \leq \|\mathfrak{H}'_\gamma\|_{H^1_\gamma(Q_T^-)} + \|\mathfrak{H}''_\gamma\|_{H^1_\gamma(Q_T^-)}$, where the first norm in the right-hand side is estimated by (62), and the second one is less than $(1/\gamma)\|\mathfrak{H}''_\gamma\|_{H^2_\gamma(Q_T^-)}$, which in turn is estimated by (81) (where $\mathfrak{H}_\gamma$ stands for $\mathfrak{H}''_\gamma$). When we estimate the boundary data $g'_2 = g_2 + \widehat{\mathcal{H}} \cdot \mathcal{H}''$ in (48), (51), we use the trace estimate $\|\mathfrak{H}''_{\gamma|\omega_T}\|_{H^{3/2}_\gamma(\omega_T)} \leq C\|\mathfrak{H}''_\gamma\|_{H^2_\gamma(Q_T^-)}$ and (81) again. □

We observe that, differently from what happens in (62), in the above a priori estimates we have a loss of one derivative from the data $\tilde{f}, \chi, \Xi, g$ to the solution $(\dot{\mathcal{U}}, \dot{\mathfrak{H}})$.

## 8. Well-posedness of the linearized problem in anisotropic Sobolev spaces

The aim of this section is to prove the following theorem.

**Theorem 16.** *Let* $T > 0$, $m \in \mathbb{N}, m \geq 1$ *and* $s = \max\{m + 2, 9\}$. *Let the basic state* (28) *satisfy assumptions* (29)–(34), (61) *and*

$$(\widehat{U}, \widehat{\mathcal{H}}, \hat{\varphi}) \in H^s_{*,\gamma}(Q_T^+) \times H^s_\gamma(Q_T^-) \times H^{s+1/2}_\gamma(\omega_T). \tag{87}$$

*There exists* $\gamma_m \geq 1$ *such that for all* $\gamma \geq \gamma_m$ *and for all* $f_\gamma \in H^{m+1}_{*,\gamma}(Q_T^+)$, $(\chi_\gamma, \Xi_\gamma) \in H^{m+1}_\gamma(Q_T^-)$, *satisfying the compatibility conditions* (66), $(g_{1\gamma}, g_{2\gamma}, g_{3\gamma}) \in H^{m+1}_\gamma(\omega_T)$, $g_{4\gamma} \in H^{m+1}_\gamma(\omega_T^+)$, $g_{5\gamma} \in H^{m+1}_\gamma(\omega_T^-)$, *all functions vanishing in the past, problem* (85) *has a solution* $(\dot{\mathcal{U}}_\gamma, \dot{\mathfrak{H}}_\gamma, \varphi_\gamma) \in H^m_{*,\gamma}(Q_T^+) \times H^m_\gamma(Q_T^-) \times H^{m+1/2}_\gamma(\omega_T)$ *with trace* $(\dot{q}_\gamma, \dot{u}_{1\gamma}, \dot{h}_{1\gamma})|_{\omega_T} \in H^{m-1/2}_\gamma(\omega_T)$. *Moreover, the solution obeys the tame estimate*

$$\gamma\Big(\|\dot{\mathcal{U}}_\gamma\|^2_{H^m_{*,\gamma}(Q_T^+)} + \|\dot{\mathfrak{H}}_\gamma\|^2_{H^m_\gamma(Q_T^-)} + \|(\dot{q}_\gamma, \dot{u}_{1\gamma}, \dot{h}_{1\gamma})|_{\omega_T}\|^2_{H^{m-1/2}_\gamma(\omega_T)} + \|\varphi_\gamma\|^2_{H^{m+1/2}_\gamma(\omega_T)}\Big)$$

$$\leq \frac{C}{\gamma}\Big\{\Big(\|f_\gamma\|^2_{H^s_{*,\gamma}(Q_T^+)} + \|\chi_\gamma, \Xi_\gamma\|^2_{H^s_\gamma(Q_T^-)} + \|g_\gamma\|^2_{H^s_\gamma(\omega_T^\pm)}\Big)\times$$

$$\times \Big(\|\widehat{U}\|_{H^{m+2}_{*,\gamma}(Q_T^+)} + \|\widehat{\mathcal{H}}\|_{H^{m+2}_\gamma(Q_T^-)} + \|\hat{\varphi}\|_{H^{m+2.5}_\gamma(\omega_T)}\Big)$$

$$+ \|f_\gamma\|^2_{H^{m+1}_{*,\gamma}(Q_T^+)} + \|\chi_\gamma, \Xi_\gamma\|^2_{H^{m+1}_\gamma(Q_T^-)} + \|g_\gamma\|^2_{H^{m+1}_\gamma(\omega_T^\pm)}\Big\}, \tag{88}$$

*where the constant* $C = C(K, T, \delta_0)$ *is independent of the data* $f, \chi, \Xi, g$ *and* $\gamma$, *and where for the sake of brevity we have set*

$$\|g_\gamma\|^2_{H^{m+1}_\gamma(\omega_T^\pm)} := \|g_{1\gamma}, g_{2\gamma}, g_{3\gamma}\|^2_{H^{m+1}_\gamma(\omega_T)} + \|g_{4\gamma}\|^2_{H^{m+1}_\gamma(\omega_T^+)} + \|g_{5\gamma}\|^2_{H^{m+1}_\gamma(\omega_T^-)}.$$

The proof proceeds by induction. Assume that Theorem 16 holds up to $m-1$. Given the data $(f, \chi, \Xi, g)$ as in Theorem 16, by the inductive hypothesis there exists a solution $(\dot{\mathcal{U}}, \dot{\mathfrak{H}}, \varphi)$ of problem (85) such that $(\dot{\mathcal{U}}_\gamma, \dot{\mathfrak{H}}_\gamma, \varphi_\gamma) \in H^{m-1}_{*,\gamma}(Q_T^+) \times H^{m-1}_\gamma(Q_T^-) \times H^{m-1}_\gamma(\omega_T)$ with trace $(\dot{q}_\gamma, \dot{u}_{1\gamma}, \dot{h}_{1\gamma})|_{\omega_T} \in H^{m-3/2}_\gamma(\omega_T)$. This solution satisfies the corresponding a priori estimate (88) of order $m-1$.

In order to show that $(\dot{\mathcal{U}}_\gamma, \dot{\mathfrak{H}}_\gamma, \varphi_\gamma) \in H^m_{*,\gamma}(Q_T^+) \times H^m_\gamma(Q_T^-) \times H^m_\gamma(\omega_T)$ we have to increase the regularity by one order. For $\dot{\mathcal{U}}$, which has to belong to the anisotropic space $H^m_{*,\gamma}$, we have to increase the regularity by one more tangential derivative and, if $m$ is even, also by one more normal derivative. The idea is the same as in [28, 30], revisited as in [6, 24, 31], with the additional difficulty of the loss of regularity from the source terms, as in [22], and the coupling with the elliptic system for $\dot{\mathfrak{H}}$.

At every step we can estimate some derivatives of $\dot{\mathcal{U}}$ through equations where in the right-hand side we can put other derivatives of $\dot{\mathcal{U}}$ that have already been estimated at previous steps. For the increase of regularity we first consider the system of equations (91) for purely tangential derivatives of $\dot{\mathcal{U}}$, coupled with the elliptic system (93) through the boundary equations (111), where we can use the inductive assumption. The difficulty is that we have to deal with the loss of one derivative in the right-hand side of (91). However the terms in the right-hand side have order $m - 1$; after the loss of one derivative they become essentially of order $m$, and can be absorbed for $\gamma$ large by similar terms in the left-hand side. The regularity of the front is obtained at this step.

Then we consider other systems (119), (122) of equations for mixed tangential and normal derivatives where the boundary matrix vanishes identically, so that no boundary condition is needed and we can apply a standard energy method.

For the sake of brevity, let us denote

$$L = \widehat{\mathcal{A}}_0 \partial_t + \sum_{j=1}^{3} (\widehat{\mathcal{A}}_j + \mathcal{E}_{1j+1}) \partial_j + \widehat{\mathcal{C}}.$$

We decompose $\dot{\mathcal{U}}$ as $\dot{\mathcal{U}} = \begin{pmatrix} \dot{\mathcal{U}}^I \\ \dot{\mathcal{U}}^{II} \end{pmatrix}$, where $\dot{\mathcal{U}}^I = (\dot{q}, \dot{u}_1)$, $\dot{\mathcal{U}}^{II} = (\dot{u}_2, \dot{u}_3, \dot{h}_1, \dot{h}_2, \dot{h}_3, \dot{S})$. A similar decomposition is used for other vectors. We also write the first two rows of $\widehat{\mathcal{A}}_1 + \mathcal{E}_{12}$ as

$$(\widehat{\mathcal{A}}_1^{I,I} + \mathcal{E}_{12}^{I,I} \quad \widehat{\mathcal{A}}_1^{I,II}), \qquad \mathcal{E}_{12}^{I,I} = \begin{pmatrix} 0 & 1 \\ 1 & 0 \end{pmatrix}.$$

The matrix $\mathcal{E}_{12}^{I,I}$ is the invertible part of $\mathcal{E}_{12}$.

### 8.1. Purely tangential regularity in the plasma part.
Let us start by considering all the tangential derivatives $Z^\alpha \dot{\mathcal{U}}$, $|\alpha| = m - 1$. By inverting $\mathcal{E}_{12}^{I,I}$ in the first two rows of (85)$_1$, we can write $\partial_1 \dot{\mathcal{U}}^I$ as the sum of tangential derivatives by

$$\partial_1 \dot{\mathcal{U}}^I = \Lambda Z \dot{\mathcal{U}} + \mathcal{R} \tag{89}$$

where

$$\Lambda Z \dot{\mathcal{U}} = -(\mathcal{E}_{12}^{I,I})^{-1} \big[ (\widehat{\mathcal{A}}_0 Z_0 \dot{\mathcal{U}} + \sum_{j=1}^{2} (\widehat{\mathcal{A}}_j + \mathcal{E}_{1j+1}) Z_j \dot{\mathcal{U}})^I + \widehat{\mathcal{A}}_1^{I,II} \partial_1 \dot{\mathcal{U}}^{II} \big],$$

$$\mathcal{R} = (\mathcal{E}_{12}^{I,I})^{-1} (\tilde{f} - \widehat{\mathcal{C}} \dot{\mathcal{U}})^I.$$

Here and below, everywhere it is needed, we use the fact that, if a given matrix $\mathcal{A}$ vanishes on $\{x_1 = 0\} \cup \{x_1 = 1\}$, we can write $\mathcal{A}\partial_1 \dot{\mathcal{U}} = MZ_1\dot{\mathcal{U}}$, where $M$ is a suitable matrix, and it holds $\|M\|_{H^{s-2}_{*,\gamma}(Q_T^+)} \le c\|\mathcal{A}\|_{H^{s}_{*,\gamma}(Q_T^+)}$, see Lemmata 46 and 47; this trick transforms some normal derivatives into tangential derivatives. We obtain $\Lambda \in H^{s-2}_{*,\gamma}(Q_T^+)$.

Applying the operator $Z^\alpha$ to (85)$_1$, with $\alpha = (\alpha_0, \alpha')$, $\alpha' = (\alpha_1, \alpha_2, \alpha_3)$, and substituting (89) gives

$$L(Z^\alpha \dot{\mathcal{U}}) + \sum_{|\gamma|=|\alpha|-1} \sum_{j \ne 1} Z \widehat{\mathcal{A}}_j Z_j Z^\gamma \dot{\mathcal{U}} + \sum_{|\gamma|=|\alpha|-1} Z \widehat{\mathcal{A}}_1 \begin{pmatrix} \Lambda Z(Z^\gamma \dot{\mathcal{U}}) \\ 0 \end{pmatrix}$$
$$- \alpha_1 (\widehat{\mathcal{A}}_1 + \mathcal{E}_{12}) \begin{pmatrix} \Lambda Z(Z_0^{\alpha_0} Z_1^{\alpha_1-1} Z_2^{\alpha_2} Z_3^{\alpha_3} \dot{\mathcal{U}}) \\ 0 \end{pmatrix} \tag{90}$$
$$+ \bigg( \sum_{|\gamma|=|\alpha|-1} Z \widehat{\mathcal{A}}_1 Z^\gamma - \alpha_1 (\widehat{\mathcal{A}}_1 + \mathcal{E}_{12}) Z_0^{\alpha_0} Z_1^{\alpha_1-1} Z_2^{\alpha_2} Z_3^{\alpha_3} \bigg) \begin{pmatrix} 0 \\ \partial_1 \dot{\mathcal{U}}^{II} \end{pmatrix} = F_\alpha, \qquad \text{in } Q_T^+,$$

where

$$\begin{aligned}
F_\alpha = \ & -\sum_{|\beta|\geq 2, \beta\leq\alpha}\left[\sum_{j\neq 1} Z^\beta \widehat{\mathcal{A}}_j Z_j Z^{\alpha-\beta}\dot{\mathcal{U}} + Z^\beta \widehat{\mathcal{A}}_1 Z^{\alpha-\beta}\left(\begin{array}{c}\Lambda Z\dot{\mathcal{U}}+\mathcal{R}\\ \partial_1\dot{\mathcal{U}}^{II}\end{array}\right)\right]\\
& -\left(\begin{array}{c}\alpha_1\\ 2\end{array}\right)(\widehat{\mathcal{A}}_1+\mathcal{E}_{12})Z_0^{\alpha_0}Z_1^{\alpha_1-2}Z_2^{\alpha_2}Z_3^{\alpha_3}\left(\begin{array}{c}\Lambda Z\dot{\mathcal{U}}+\mathcal{R}\\ \partial_1\dot{\mathcal{U}}^{II}\end{array}\right)\\
& -(\widehat{\mathcal{A}}_1+\mathcal{E}_{12})\partial_1 Z_0^{\alpha_0}\Big[(Z_1-1)^{\alpha_1}-Z_1^{\alpha_1}+\alpha_1(Z_1-1)^{\alpha_1-1}\\
& -\left(\begin{array}{c}\alpha_1\\ 2\end{array}\right)(Z_1-1)^{\alpha_1-2}\Big]Z_2^{\alpha_2}Z_3^{\alpha_3}\dot{\mathcal{U}} -\sum_{|\alpha'|=|\alpha|-1}Z\widehat{\mathcal{A}}_1\left[Z^{\alpha'},\left(\begin{array}{c}\Lambda\\ 0\end{array}\right)\right]Z\dot{\mathcal{U}}\\
& +\alpha_1(\widehat{\mathcal{A}}_1+\mathcal{E}_{12})\left[Z_0^{\alpha_0}Z_1^{\alpha_1-1}Z_2^{\alpha_2}Z_3^{\alpha_3},\left(\begin{array}{c}\Lambda\\ 0\end{array}\right)\right]Z\dot{\mathcal{U}} -\left[Z^\alpha,\widehat{\mathcal{C}'}\right]\dot{\mathcal{U}}\\
& -\left(\sum_{|\alpha'|=|\alpha|-1}Z\widehat{\mathcal{A}}_1 Z^{\alpha'}-\alpha_1(\widehat{\mathcal{A}}_1+\mathcal{E}_{12})Z_0^{\alpha_0}Z_1^{\alpha_1-1}Z_2^{\alpha_2}Z_3^{\alpha_3}\right)\left(\begin{array}{c}\mathcal{R}\\ 0\end{array}\right) + Z^\alpha\tilde{f}.
\end{aligned}$$

$[\,\cdot\,,\,\cdot\,]$ denotes the commutator. Equation (90) takes the form $(L+B)Z^\alpha\dot{\mathcal{U}}=F_\alpha$ with $B\in H^{s-3}_{*,\gamma}(Q_T^+)$.

Then we consider the problem satisfied by the vector of all tangential derivatives $Z^\alpha\dot{\mathcal{U}}$ of order $|\alpha|=m-1$. From (90) this problem takes the form

$$(\mathcal{L}+\mathcal{B})Z^\alpha\dot{\mathcal{U}}=\mathcal{F}_\alpha \quad \text{in } Q_T^+, \tag{91}$$

where

$$\mathcal{L}=\left(\begin{array}{ccc}L & & \\ & \ddots & \\ & & L\end{array}\right),$$

$\mathcal{B}\in H^{s-3}_{*,\gamma}(Q_T^+)$ is a suitable matrix and $\mathcal{F}_\alpha$ is the vector of all right-hand sides $F_\alpha$.

In order to increase by one tangential derivative the regularity of $Z^\alpha\dot{\mathcal{U}}$ we will apply Theorem 15 and in particular the a priori estimate (86). For this, we have to estimate $\mathcal{F}_{\alpha\gamma}$ in $H^2_{tan,\gamma}(Q_T^+)$. Here special care is needed because of the loss of one derivative in (86) from the source term to the solution. However, this is the same calculation of [22], p. 77-79. Proceeding as in [22] gives from Corollary 44 and Theorem 45 the estimate

$$\begin{aligned}
\|\mathcal{F}_{\alpha\gamma}\|_{H^2_{tan,\gamma}(Q_T^+)} \leq C\Big(&\|\dot{\mathcal{U}}_\gamma\|_{H^m_{tan,\gamma}(Q_T^+)}+\|\partial_1\dot{\mathcal{U}}^{II}_\gamma\|_{H^{m-2}_{tan,\gamma}(Q_T^+)}\\
& +\|\dot{\mathcal{U}}_\gamma\|_{H^7_{*,\gamma}(Q_T^+)}\|\widehat{U},\nabla\widehat{\Psi}\|_{H^{m+1}_{*,\gamma}(Q_T^+)}+\|\tilde{f}_\gamma\|_{H^{m+1}_{tan,\gamma}(Q_T^+)}\Big), \tag{92}
\end{aligned}$$

where the constant $C$ depends on $\|\widehat{U},\nabla\widehat{\Psi}\|_{H^9_{*,\gamma}(Q_T^+)}$.

## 8.2. Purely tangential regularity in the vacuum part.

Let us consider now the equations $(85)_2$ for the vacuum magnetic field. Applying the operator $Z^\alpha$ gives

$$\nabla\times Z^\alpha\dot{\mathfrak{H}}=\chi_\alpha, \quad \text{div}\,(A\,Z^\alpha\dot{\mathfrak{H}})=\Xi_\alpha \quad \text{in } Q_T^-, \tag{93}$$

where we have set

$$\chi_\alpha=Z^\alpha\chi-[Z^\alpha,\nabla\times]\dot{\mathfrak{H}}, \qquad \Xi_\alpha=Z^\alpha\Xi-[Z^\alpha,\text{div}\,(A\cdot)]\dot{\mathfrak{H}}. \tag{94}$$

For the estimation of the commutators in $\chi_\alpha,\Xi_\alpha$, it is crucial to have introduced the functional spaces $H^{-1}_{\tau\Gamma},H^{-1}_{0\Gamma_-}$ (in space variables) with negative order, in order to compensate the $H^2$ norm (in time) appearing in the right-hand side of (86).

**Lemma 17.** *The following estimate holds:*

$$\gamma\|\chi_{\alpha\gamma},\Xi_{\alpha\gamma}\|_{H^1_\gamma(Q^-_T)} + \|\chi_{\alpha\gamma},\Xi_{\alpha\gamma}\|_{H^2_\gamma(-\infty,T;H^{-1}_{\gamma T}\times H^{-1}_{0\Gamma_-})}$$

$$\leq C(K)\big\{\|\|\dot{\mathfrak{H}}_\gamma\|_{H^m_{\tan,\gamma}(Q^-_T)} + \|\chi_\gamma,\Xi_\gamma\|_{H^{m+1}_{\tan,\gamma}(Q^-_T)} + \|\varphi_\gamma\|_{H^{m+1/2}_\gamma(\omega_T)}$$

$$+ (\|\dot{\mathfrak{H}}_\gamma\|_{H^4_\gamma(Q^-_T)} + \|\varphi_\gamma\|_{H^5_\gamma(\omega_T)})(\|\varphi_\gamma\|_{H^{m+1/2}_\gamma(\omega_T)} + \|\widehat{\mathcal{H}}\|_{H^{m+1}_\gamma(Q^-_T)})$$

$$+ \|g_{3\gamma}\|_{H^{m-1/2}_\gamma(\omega_T)} + \|g_{5\gamma}\|_{H^{m-1/2}_\gamma(\omega_T)}\big\}. \tag{95}$$

*Proof.* Notice that $[Z^\alpha,\partial_i] = 0$ for $i = 2,3$, and $[Z^\alpha,\partial_1] \neq 0$ only if $\alpha_1 > 0$. In fact, it is

$$[Z_1^{\alpha_1},\partial_1] = -\sigma'\partial_1 Z_1^{\alpha_1-1} - Z_1(\sigma'\partial_1 Z_1^{\alpha_1-2}) - \cdots - Z_1^{\alpha_1-1}(\sigma'\partial_1). \tag{96}$$

It follows from (96) that $[Z^\alpha,\nabla\times]\dot{\mathfrak{H}}$ may be written as[5]

$$[Z^\alpha,\nabla\times]\dot{\mathfrak{H}}_\gamma = P_{|\alpha|-1}(Z)\partial_1\dot{\mathfrak{H}}_\gamma, \tag{97}$$

where $P_{|\alpha|-1}$ is a polynomial in $Z$ of degree $|\alpha|-1$ with $C^\infty$ coefficients only dependent on $x_1$. For the sake of simplicity let us first assume $m = 2$, so that $\alpha = \alpha_1 = 1$ and $P_0 = -\sigma'$. By an integration by parts we have

$$\|\sigma'\partial_1\dot{\mathfrak{H}}_\gamma\|_{H^{-1}_{\gamma T}} = \sup_{\psi\in H^1_{\gamma T}}\frac{|\int_{\Omega^-}\sigma'\partial_1\dot{\mathfrak{H}}_\gamma\cdot\psi\,dx|}{\|\psi\|_{H^1}}$$
$$= \sup_{\psi\in H^1_{\gamma T}}\frac{|\int_\Gamma \sigma'\dot{\mathfrak{H}}_\gamma\cdot\psi\,dx' - \int_{\Gamma_-}\sigma'\dot{\mathfrak{H}}_\gamma\cdot\psi\,dx' - \int_{\Omega^-}\dot{\mathfrak{H}}_\gamma\cdot\partial_1(\sigma'\psi)\,dx|}{\|\psi\|_{H^1}}. \tag{98}$$

Regarding the boundary integrals, $\psi\in H^1_{\gamma\Gamma}$ and the boundary condition on $\Gamma_-$ give

$$\dot{\mathfrak{H}}\cdot\psi = \dot{\mathfrak{H}}_1\psi_1 \qquad \text{on } \Gamma, \qquad \dot{\mathfrak{H}}\cdot\psi = \dot{\mathfrak{H}}_1\psi_1 - g_{5,3}\psi_2 + g_{5,2}\psi_3 \qquad \text{on } \Gamma_-. \tag{99}$$

Since the trace $\psi\in H^{1/2}(\Gamma\cup\Gamma_-)$, for the estimate of the boundary integrals in (98), after (99) we need an estimate of $\dot{\mathfrak{H}}_1$ in $H^{-1/2}(\Gamma\cup\Gamma_-)$. From[6]

$$(A\dot{\mathfrak{H}})_1 = \langle\nabla\hat{\varphi}\rangle^2\dot{\mathfrak{H}}_1 - \dot{\mathfrak{H}}_2\partial_2\hat{\varphi} - \dot{\mathfrak{H}}_3\partial_3\hat{\varphi} \quad \text{on } \Gamma, \qquad (A\dot{\mathfrak{H}})_1 = (\partial_1\widehat{\Phi}_1)^{-1}\dot{\mathfrak{H}}_1 \quad \text{on } \Gamma_-,$$

we get

$$\int_\Gamma \sigma'\dot{\mathfrak{H}}_{1\gamma}\psi_1\,dx' - \int_{\Gamma_-}\sigma'\dot{\mathfrak{H}}_{1\gamma}\psi_1\,dx'$$
$$= \int_\Gamma \sigma'\langle\nabla\hat{\varphi}\rangle^{-2}(A\dot{\mathfrak{H}}_\gamma)_1\psi_1\,dx' - \int_{\Gamma_-}\sigma'\partial_1\widehat{\Phi}_1(A\dot{\mathfrak{H}}_\gamma)_1\psi_1\,dx' + \int_\Gamma \sigma'\langle\nabla\hat{\varphi}\rangle^{-2}(\dot{\mathfrak{H}}_{2\gamma}\partial_2\hat{\varphi} + \dot{\mathfrak{H}}_{3\gamma}\partial_3\hat{\varphi})\psi_1\,dx'$$
$$= \int_{\Omega^-}\text{div}\,(\sigma'\psi_1\langle\nabla\hat{\varphi}\rangle^{-2}A\dot{\mathfrak{H}}_\gamma)\,dx + \int_\Gamma \sigma'\langle\nabla\hat{\varphi}\rangle^{-2}(\dot{\mathfrak{H}}_{2\gamma}\partial_2\hat{\varphi} + \dot{\mathfrak{H}}_{3\gamma}\partial_3\hat{\varphi})\psi_1\,dx'$$
$$= \int_{\Omega^-}\Big(\sigma'\psi_1\langle\nabla\hat{\varphi}\rangle^{-2}\Xi_\gamma + A\dot{\mathfrak{H}}_\gamma\cdot\nabla(\sigma'\psi_1\langle\nabla\hat{\varphi}\rangle^{-2})\Big)\,dx + \int_\Gamma \sigma'\langle\nabla\hat{\varphi}\rangle^{-2}(\dot{\mathfrak{H}}_{2\gamma}\partial_2\hat{\varphi} + \dot{\mathfrak{H}}_{3\gamma}\partial_3\hat{\varphi})\psi_1\,dx'. \tag{100}$$

Thus, from (98), (99), (100) and the well-known trace estimate

$$\|(\dot{\mathfrak{H}}_{2\gamma},\dot{\mathfrak{H}}_{3\gamma})_{|\Gamma}\|_{H^{-1/2}(\Gamma)} \leq C\left(\|\dot{\mathfrak{H}}_\gamma\|_{L^2(\Omega^-)} + \|\nabla\times\dot{\mathfrak{H}}_\gamma\|_{L^2(\Omega^-)}\right) \tag{101}$$

we get

$$\|\sigma'\partial_1\dot{\mathfrak{H}}_\gamma\|_{H^{-1}_{\gamma T}} \leq C\left(\|\dot{\mathfrak{H}}_\gamma\|_{L^2(\Omega^-)} + \|\chi_\gamma,\Xi_\gamma\|_{L^2(\Omega^-)} + \|g_{5\gamma}\|_{H^{-1/2}(\Gamma_-)}\right).$$

For $m$ in the general case as well as for the time derivatives of $\partial_1\dot{\mathfrak{H}}_\gamma$ we follow similar arguments and obtain from (97), the calculus inequality (269) and Sobolev imbeddings

$$\|[Z^\alpha,\nabla\times]\dot{\mathfrak{H}}_\gamma\|_{H^2_\gamma(-\infty,T;H^{-1}_{\gamma T})} \leq C\big\{\|\dot{\mathfrak{H}}_\gamma\|_{H^m_{\tan,\gamma}(Q^-_T)} + \|\chi_\gamma,\Xi_\gamma\|_{H^m_{\tan,\gamma}(Q^-_T)}$$
$$+ \|\dot{\mathfrak{H}}_\gamma\|_{H^4_\gamma(Q^-_T)}\|\hat{\varphi}_\gamma\|_{H^{m+1/2}_\gamma(\omega_T)} + \|g_{5\gamma}\|_{H^{m-1/2}_\gamma(\omega^-_T)}\big\}.$$

---

[5]The equality follows by commuting $\partial_1$ and $Z_1$ till when $\partial_1$ always stands on the right-hand side of each term.

[6]Here we use the notation $\langle\nabla\hat{\varphi}\rangle := (1 + |\partial_2\hat{\varphi}|^2 + |\partial_3\hat{\varphi}|^2)^{1/2}$.

Adding the estimate of $Z^\alpha \chi$ gives

$$\|\chi_{2\gamma}\|_{H^2_\gamma(-\infty,T;H^{-1}_\Gamma)} \leq C\big\{\|\dot{\mathfrak{H}}_\gamma\|_{H^m_{tan,\gamma}(Q^-_T)} + \|\chi_\gamma, \Xi_\gamma\|_{H^{m+1}_{tan,\gamma}(Q^-_T)}$$
$$+ \|\dot{\mathfrak{H}}_\gamma\|_{H^4_\gamma(Q^-_T)}\|\hat{\varphi}_\gamma\|_{H^{m+1/2}_\gamma(\omega_T)} + \|g_{5\gamma}\|_{H^{m-1/2}_\gamma(\omega^-_T)}\big\}. \quad (102)$$

which provides half of (95) for the part of $\chi_\alpha$; in analogous way we prove the estimate of $\|\chi_{2\gamma}\|_{H^1_\gamma(Q^-_T)}$.

Let us consider $\Xi_\alpha$. As to the second commutator in (94) we have

$$[Z^\alpha, \operatorname{div}(A\cdot)]\dot{\mathfrak{H}} = Z^\alpha \operatorname{div}(A\dot{\mathfrak{H}}) - \operatorname{div} Z^\alpha(A\dot{\mathfrak{H}}) + \operatorname{div} Z^\alpha(A\dot{\mathfrak{H}}) - \operatorname{div}(AZ^\alpha\dot{\mathfrak{H}})$$
$$= [Z^\alpha, \partial_1](A\dot{\mathfrak{H}})_1 + \operatorname{div}[Z^\alpha, A\cdot]\dot{\mathfrak{H}} = P_{|\alpha|-1}(Z)\partial_1(A\dot{\mathfrak{H}})_1 + \operatorname{div}[Z^\alpha, A\cdot]\dot{\mathfrak{H}}, \quad (103)$$

from $[Z^\alpha, Z_i] = 0$ if $i = 2, 3$, and (96), where again $P_{|\alpha|-1}$ is a polynomial in $Z$ of degree $|\alpha| - 1$ with $C^\infty$ coefficients dependent on $x_1$. To estimate the first term in the right side assume for the sake of simplicity that $m = 2$, which yields $P_0 = -\sigma'$. Integrating by parts gives

$$\|\sigma'\partial_1(A\dot{\mathfrak{H}})_1\|_{H^{-1}_{0\Gamma_-}} = \sup_{\phi \in H^1_{0\Gamma_-}} \frac{|\int_{\Omega^-} \sigma'\partial_1(A\dot{\mathfrak{H}})_1\phi\, dx|}{\|\phi\|_{H^1_{0\Gamma_-}}}$$
$$= \sup_{\phi \in H^1_{0\Gamma_-}} \frac{|\int_\Gamma \sigma'(A\dot{\mathfrak{H}})_1\phi\, dx' - \int_{\Gamma_-}\sigma'(A\dot{\mathfrak{H}})_1\phi\, dx' - \int_{\Omega^-}(A\dot{\mathfrak{H}})_1\partial_1(\sigma'\phi)\, dx|}{\|\nabla\phi\|_{L^2(\Omega^-)}}$$
$$= \sup_{\phi \in H^1_{0\Gamma_-}} \frac{|\int_{\Omega^-} \operatorname{div}(\sigma'\phi A\dot{\mathfrak{H}}_\gamma)\, dx - \int_{\Omega^-}(A\dot{\mathfrak{H}}_\gamma)_1\partial_1(\sigma'\phi)\, dx|}{\|\nabla\phi\|_{L^2(\Omega^-)}}$$
$$= \sup_{\phi \in H^1_{0\Gamma_-}} \frac{|\int_{\Omega^-}\left(\sigma'\phi\Xi_\gamma + A\dot{\mathfrak{H}}_\gamma \cdot \nabla(\sigma'\phi)\right)\, dx - \int_{\Omega^-}(A\dot{\mathfrak{H}}_\gamma)_1\partial_1(\sigma'\phi)\, dx|}{\|\nabla\phi\|_{L^2(\Omega^-)}},$$

and we readily get

$$\|\sigma'\partial_1(A\dot{\mathfrak{H}}_\gamma)_1\|_{H^{-1}_{0\Gamma_-}} \leq C\left(\|\dot{\mathfrak{H}}_\gamma\|_{L^2(\Omega^-)} + \|\Xi_\gamma\|_{L^2(\Omega^-)}\right).$$

For $m$ in the general case a similar argument gives

$$\|P_{|\alpha|-1}(Z)\partial_1(A\dot{\mathfrak{H}})_1\|_{H^2_\gamma(-\infty,T;H^{-1}_{0\Gamma_-})} \leq C\big\{\|\dot{\mathfrak{H}}_\gamma\|_{H^m_{tan,\gamma}(Q^-_T)} + \|\Xi_\gamma\|_{H^m_{tan,\gamma}(Q^-_T)}$$
$$+ \|\dot{\mathfrak{H}}_\gamma\|_{H^4_\gamma(Q^-_T)}\|\hat{\varphi}_\gamma\|_{H^{m+1/2}_\gamma(\omega_T)}\big\}. \quad (104)$$

Finally, consider the last term in the right-hand side of (103)

$$\operatorname{div}[Z^\alpha, A\cdot]\dot{\mathfrak{H}} = \sum_{i,j}\partial_i[Z^\alpha, A_{ij}\cdot]\dot{\mathfrak{H}}_j.$$

If $i = 1$ we compute

$$\|\partial_1[Z^\alpha, A_{1j}\cdot]\dot{\mathfrak{H}}_\gamma\|_{H^{-1}_{0\Gamma_-}} = \sup_{\phi \in H^1_{0\Gamma_-}} \frac{|\int_\Gamma[Z^\alpha, A_{1j}\cdot]\dot{\mathfrak{H}}_\gamma\phi\, dx' - \int_{\Omega^-}[Z^\alpha, A_{1j}\cdot]\dot{\mathfrak{H}}_\gamma\partial_1\phi\, dx|}{\|\nabla\phi\|_{L^2(\Omega^-)}}. \quad (105)$$

Here the problem is how to estimate the trace of $\dot{\mathfrak{H}}_{1\gamma}$ in the boundary integral, while the traces of $\dot{\mathfrak{H}}_{2\gamma}, \dot{\mathfrak{H}}_{3\gamma}$ can be controlled as in (101). From the boundary condition

$$\dot{\mathfrak{H}}_1 = \hat{\mathcal{H}}_1\partial_1\hat{\Phi}_1 = \sum_{i=2,3}\left(\hat{\mathcal{H}}_i\partial_i\hat{\varphi} + \partial_i(\hat{\mathcal{H}}_i\varphi)\right) + g_3 \quad \text{on } \Gamma, \quad (106)$$

we have

$$\dot{\mathfrak{H}}_2 = \hat{\mathcal{H}}_1\partial_2\hat{\Psi} + \hat{\mathcal{H}}_2 = \left(\sum_i\left(\hat{\mathcal{H}}_i\partial_i\hat{\varphi} + \partial_i(\hat{\mathcal{H}}_i\varphi)\right) + g_3\right)\partial_2\hat{\varphi} + \hat{\mathcal{H}}_2,$$
$$\dot{\mathfrak{H}}_3 = \hat{\mathcal{H}}_1\partial_3\hat{\Psi} + \hat{\mathcal{H}}_3 = \left(\sum_i\left(\hat{\mathcal{H}}_i\partial_i\hat{\varphi} + \partial_i(\hat{\mathcal{H}}_i\varphi)\right) + g_3\right)\partial_3\hat{\varphi} + \hat{\mathcal{H}}_3,$$

and this system can be rewritten as

$$\begin{aligned}
(1 + |\partial_2\hat\varphi|^2)\dot{\mathcal{H}}_2 + \partial_2\hat\varphi\partial_3\hat\varphi\dot{\mathcal{H}}_3 &= \dot{\mathfrak{H}}_2 - \Big(\sum_i \partial_i(\widehat{\mathcal{H}}_i\varphi) + g_3\Big)\partial_2\hat\varphi, \\
\partial_2\hat\varphi\partial_3\hat\varphi\dot{\mathcal{H}}_2 + (1 + |\partial_3\hat\varphi|^2)\dot{\mathcal{H}}_3 &= \dot{\mathfrak{H}}_3 - \Big(\sum_i \partial_i(\widehat{\mathcal{H}}_i\varphi) + g_3\Big)\partial_3\hat\varphi.
\end{aligned} \tag{107}$$

Since system (107) has determinant $1 + |\partial_2\hat\varphi|^2 + |\partial_3\hat\varphi|^2 \neq 0$, we can write $\dot{\mathcal{H}}_2, \dot{\mathcal{H}}_3$ in terms of the right-hand sides. The values of $\dot{\mathcal{H}}_2, \dot{\mathcal{H}}_3$ that are controlled in this way are then used for the estimate of $\dot{\mathfrak{H}}_{1\gamma}$ from (106). From (105)–(107) we get

$$\begin{aligned}
\|\partial_1[Z^\alpha, A_{1j}\cdot]\dot{\mathfrak{H}}_j\|_{H^2_\gamma(-\infty,T;H^{-1}_{0\Gamma_-})} &\leq C\big\{\|\dot{\mathfrak{H}}_\gamma\|_{H^m_{tan,\gamma}(Q^-_T)} + \|\varphi_\gamma\|_{H^{m+1/2}_\gamma(\omega_T)} \\
&\quad + (\|\dot{\mathfrak{H}}_\gamma\|_{H^3_\gamma(Q^-_T)} + \|\varphi_\gamma\|_{H^3_\gamma(\omega_T)})(\|\hat\varphi_\gamma\|_{H^{m+1/2}_\gamma(\omega_T)} + \|\widehat{\mathcal{H}}\|_{H^{m+1}_\gamma(Q^-_T)}) \\
&\quad + \|g_{3\gamma}\|_{H^{m-1/2}_\gamma(\omega_T)} + \|g_{5\gamma}\|_{H^{m-1/2}_\gamma(\omega_T)}\big\}.
\end{aligned} \tag{108}$$

For $i = 2, 3$, applying the last inequality in (65) yields

$$\|\partial_i[Z^\alpha, A_{ij}\cdot]\dot{\mathfrak{H}}_j\|_{H^2_\gamma(-\infty,T;H^{-1}_{0\Gamma_-})} \leq \|[Z^\alpha, A_{ij}\cdot]\dot{\mathfrak{H}}_j\|_{H^2_\gamma(-\infty,T;L^2)}.$$

By Leibniz's formula, we expand $[Z^\alpha, A_{ij}\cdot]\dot{\mathfrak{H}}_j$ (we write it in informal way, only taking account of the order of derivatives and dropping the irrelevant numerical coefficients)

$$[Z^\alpha, A_{ij}\cdot]\dot{\mathfrak{H}}_j = \sum_{|\beta|<|\alpha|} Z^{\alpha-\beta} A_{ij} Z^\beta \dot{\mathfrak{H}}_j. \tag{109}$$

Applying the calculus inequality (265) to each term of the above expansion and a Sobolev imbedding yields

$$\|[Z^\alpha, A_{ij}\cdot]\dot{\mathfrak{H}}_{j\gamma}\|_{H^2_\gamma(-\infty,T;L^2)} \leq C(\|\dot{\mathfrak{H}}_\gamma\|_{H^m_{tan,\gamma}(Q^-_T)} + \|\dot{\mathfrak{H}}_\gamma\|_{H^3_\gamma(Q^-_T)}\|\varphi_\gamma\|_{H^{m+1/2}_\gamma(\omega_T)}), \qquad i = 2, 3. \tag{110}$$

Adding the estimate of $Z^\alpha\Xi$, from (104), (108), (110) we obtain (95) for the part of $\Xi_\alpha$. The proof of $\|\Xi_\gamma\|_{H^1_\gamma(Q^-_T)}$ is similar so we omit it for the sake of brevity. $\qquad\square$

8.3. **Regularity on the boundary.** Applying the operator $Z^\alpha$ to the boundary conditions in (85) gives

$$\begin{aligned}
&(\partial_t + \hat v_2\partial_2 + \hat v_3\partial_3)Z^\alpha\varphi = Z^\alpha\dot v_1 + Z^\alpha\varphi\,\partial_1\hat v_N + g_{1\alpha}, \\
&Z^\alpha\dot q = \hat{\mathfrak{b}}\cdot Z^\alpha\dot{\mathfrak{H}} - [\partial_1\hat q]Z^\alpha\varphi + g_{2\alpha}, \\
&(AZ^\alpha\dot{\mathfrak{H}})_1 = \partial_2(\widehat{\mathcal{H}}_2 Z^\alpha\varphi) + \partial_3(\widehat{\mathcal{H}}_3 Z^\alpha\varphi) + g_{3\alpha} \qquad \text{on } \omega_T, \\
&Z^\alpha\dot v_1 = Z^\alpha g_4 \quad \text{on } \omega^+_T, \qquad \nu \times Z^\alpha\dot{\mathcal{H}} = Z^\alpha g_5 \quad \text{on } \omega^-_T,
\end{aligned} \tag{111}$$

where we have set

$$\begin{aligned}
g_{1\alpha} &= Z^\alpha g_1 + [Z^\alpha, \partial_1\hat v_N]\varphi - [Z^\alpha, \hat v_2\partial_2 + \hat v_3\partial_3]\varphi, \\
g_{2\alpha} &= Z^\alpha g_2 + [Z^\alpha, \hat{\mathfrak{b}}\cdot]\dot{\mathfrak{H}} - [Z^\alpha, [\partial_1\hat q]]\varphi, \\
g_{3\alpha} &= Z^\alpha g_3 + [Z^\alpha, \partial_2(\widehat{\mathcal{H}}_2\cdot) + \partial_3(\widehat{\mathcal{H}}_3\cdot)]\varphi - ([Z^\alpha, A]\dot{\mathfrak{H}})_1.
\end{aligned} \tag{112}$$

It is understood that these terms make sense only for $\alpha_1 = 0$, because of the weight $\sigma$ in $Z_1$, vanishing at $\omega_T$. We immediately get

$$\begin{aligned}
\|Z^\alpha g_{4\gamma}\|_{H^{3/2}_\gamma(\omega^+_T)} + \|Z^\alpha g_{5\gamma}\|_{H^{3/2}_\gamma(\omega^-_T)} + \|Z^\alpha g_{5\gamma}\|_{H^2_\gamma(-\infty,T;H^{-1/2}_\gamma(\Gamma_-))} \\
\leq \|g_{4\gamma}\|_{H^{m+1/2}_\gamma(\omega^+_T)} + \|g_{5\gamma}\|_{H^{m+1}_\gamma(\omega^-_T)}.
\end{aligned} \tag{113}$$

For the other more involved terms we prove:

**Lemma 18.** *The data $g_{i\alpha}, i = 1, 2, 3$, defined in (112), satisfy the estimates*

$$\begin{aligned}
\|g_{1\alpha\gamma}, g_{2\alpha\gamma}, g_{3\alpha\gamma}\|_{H^{3/2}_\gamma(\omega_T)} &\leq \|g_{1\gamma}, g_{2\gamma}, g_{3\gamma}\|_{H^{m+1/2}_\gamma(\omega_T)} + C(K)\Big(\|\dot{\mathfrak{H}}_\gamma\|_{H^m_{tan,\gamma}(Q^-_T)} + \|\varphi_\gamma\|_{H^{m+1/2}_\gamma(\omega_T)} \\
&\quad + \Big(\|\widehat U\|_{H^{m+2}_{*,\gamma}(Q^+_T)} + \|\widehat{\mathcal{H}}\|_{H^{m+2}_\gamma(Q^-_T)} + \|\nabla\widehat\Psi\|_{H^{m+2}_\gamma(Q_T)}\Big)\Big(\|\dot{\mathfrak{H}}_\gamma\|_{H^3_\gamma(Q^-_T)} + \|\varphi_\gamma\|_{H^3_\gamma(\omega_T)}\Big)\Big),
\end{aligned} \tag{114}$$

$$\|g_{3\alpha\gamma}\|_{H^2_\gamma(-\infty,T;H^{-1/2}_\gamma(\Gamma))} \le \|g_{3\gamma}\|_{H^{m+1}_\gamma(\omega_T)}$$
$$+ C(K)\Big(\|\dot{\mathfrak{H}}_\gamma\|_{H^m_{tan,\gamma}(Q^-_T)} + \|\chi_\gamma\|_{H^m_{tan,\gamma}(Q^-_T)} + \|\varphi_\gamma\|_{H^{m+1/2}_\gamma(\omega_T)}$$
$$+ \big(\|\widehat{\mathcal{H}}\|_{H^{m+2}_\gamma(Q^-_T)} + \|\nabla\widehat{\Psi}\|_{H^{m+2}_\gamma(Q_T)}\big)\big(\|\dot{\mathfrak{H}}_\gamma\|_{H^4_\gamma(Q^-_T)} + \|\varphi_\gamma\|_{H^3_\gamma(\omega_T)}\big)\Big). \quad (115)$$

*Proof.* As the commutators in (112) are meaningful only for $\alpha_1 = 0$, they only involve derivatives $Z_0, Z_2, Z_3$, i.e. standard derivatives (not conormal) and the standard analysis of commutators applies. The proof of (114) follows by well-known commutator estimates (see Lemma 48), standard Sobolev imbeddings and Theorem 41. The proof of (115) follows by using part of the arguments needed for the proof of Lemma 17. $\qquad \square$

### 8.4. A priori estimate for purely tangential derivatives.
Now we apply the a priori estimates (86) to the solutions of (91), (93), (111). This is a compound system with the same structure of (85), except for the addition of the zero order terms $\mathcal{B}$ in (91). A check of the proof of Theorem 8 in [35] immediately shows that it works well even with added zero order terms in (85)$_1$, w.r.t. $\dot{\mathcal{U}}$. From (86), (92), (95), (113), (114), (115) we have

$$\gamma\Big(\|\dot{\mathcal{U}}_\gamma\|^2_{H^m_{tan,\gamma}(Q^+_T)} + \|\dot{\mathfrak{H}}_\gamma\|^2_{H^m_{tan,\gamma}(Q^-_T)} + \|\nabla\dot{\mathfrak{H}}_\gamma\|^2_{H^{m-1}_{tan,\gamma}(Q^-_T)}$$
$$+ \|(\dot{q}_\gamma, \dot{u}_{1\gamma}, \dot{h}_{1\gamma})|_{\omega_T}\|^2_{H^{m-1/2}_\gamma(\omega_T)} + \|\varphi_\gamma\|^2_{H^{m+1/2}_\gamma(\omega_T)}\Big)$$
$$\le \frac{C}{\gamma}\Big\{\|\dot{\mathcal{U}}_\gamma\|^2_{H^m_{tan,\gamma}(Q^+_T)} + \|\partial_1\dot{\mathcal{U}}^{II}_\gamma\|^2_{H^{m-2}_{tan,\gamma}(Q^+_T)} + \|\dot{\mathfrak{H}}_\gamma\|_{H^m_{tan,\gamma}(Q^-_T)} + \|\varphi_\gamma\|^2_{H^{m+1/2}_\gamma(\omega_T)}$$
$$+ \Big(\|\dot{\mathcal{U}}_\gamma\|^2_{H^4_{*,\gamma}(Q^+_T)} + \|\dot{\mathfrak{H}}_\gamma\|_{H^4_\gamma(Q^-_T)} + \|\varphi_\gamma\|_{H^3_\gamma(\omega_T)}\Big)\Big(\|\widehat{U}\|_{H^{m+2}_{*,\gamma}(Q^+_T)} + \|\widehat{\mathcal{H}}\|_{H^{m+2}_\gamma(Q^-_T)} + \|\widehat{\Psi}\|_{H^{m+3}_\gamma(Q_T)}\Big)$$
$$+ \|\tilde{f}_\gamma\|^2_{H^{m+1}_{tan,\gamma}(Q^+_T)} + \|\chi_\gamma, \Xi_\gamma\|_{H^{m+1}_{tan,\gamma}(Q^-_T)} + \|g_{1\gamma}, g_{2\gamma}\|_{H^{m+1/2}_\gamma(\omega_T)}$$
$$+ \|g_{3\gamma}\|_{H^{m+1}_\gamma(\omega_T)} + \|g_{4\gamma}\|_{H^{m+1/2}_\gamma(\omega^+_T)} + \|g_{5\gamma}\|_{H^{m+1}_\gamma(\omega^-_T)}\Big\}. \quad (116)$$

The constant $C$ depends on $K$. Taking $\gamma$ sufficiently large yields

$$\gamma\Big(\|\dot{\mathcal{U}}_\gamma\|^2_{H^m_{tan,\gamma}(Q^+_T)} + \|\dot{\mathfrak{H}}_\gamma\|^2_{H^m_{tan,\gamma}(Q^-_T)} + \|\nabla\dot{\mathfrak{H}}_\gamma\|^2_{H^{m-1}_{tan,\gamma}(Q^-_T)}$$
$$+ \|(\dot{q}_\gamma, \dot{u}_{1\gamma}, \dot{h}_{1\gamma})|_{\omega_T}\|^2_{H^{m-1/2}_\gamma(\omega_T)} + \|\varphi_\gamma\|^2_{H^{m+1/2}_\gamma(\omega_T)}\Big) \le \frac{C}{\gamma}\Big\{\|\partial_1\dot{\mathcal{U}}^{II}_\gamma\|^2_{H^{m-2}_{tan,\gamma}(Q^+_T)}$$
$$+ \Big(\|\dot{\mathcal{U}}_\gamma\|^2_{H^4_{*,\gamma}(Q^+_T)} + \|\dot{\mathfrak{H}}_\gamma\|_{H^4_\gamma(Q^-_T)} + \|\varphi_\gamma\|_{H^3_\gamma(\omega_T)}\Big)\Big(\|\widehat{U}\|_{H^{m+2}_{*,\gamma}(Q^+_T)} + \|\widehat{\mathcal{H}}\|_{H^{m+2}_\gamma(Q^-_T)} + \|\widehat{\Psi}\|_{H^{m+3}_\gamma(Q_T)}\Big)$$
$$+ \|\tilde{f}_\gamma\|^2_{H^{m+1}_{tan,\gamma}(Q^+_T)} + \|\chi_\gamma, \Xi_\gamma\|_{H^{m+1}_{tan,\gamma}(Q^-_T)} + \|g_{1\gamma}, g_{2\gamma}\|_{H^{m+1/2}_\gamma(\omega_T)}$$
$$+ \|g_{3\gamma}\|_{H^{m+1}_\gamma(\omega_T)} + \|g_{4\gamma}\|_{H^{m+1/2}_\gamma(\omega^+_T)} + \|g_{5\gamma}\|_{H^{m+1}_\gamma(\omega^-_T)}\Big\}, \quad (117)$$

where the constant $C$ depends on $K$.

### 8.5. Tangential and one normal derivatives.
We apply to the part $II$ of (85)$_1$ (i.e. to the equations for $\dot{\mathcal{U}}^{II} = (\dot{u}_2, \dot{u}_3, \dot{h}_1, \dot{h}_2, \dot{h}_3, \dot{S})$) the operator $Z^\beta\partial_1$, with $|\beta| = m-2$. We obtain equation (28) in [6], that is

$$[(L + \partial_1\widehat{\mathcal{A}}_1)Z^\beta + \sum_{|\gamma|=|\beta|-1}(Z\widehat{\mathcal{A}}_0\partial_t + \sum_{j=1}^n Z\widehat{\mathcal{A}}_j\partial_j)Z^\gamma$$
$$-\beta_1(\widehat{\mathcal{A}}_1 + \mathcal{E}_{12})\partial_1 Z^{\beta_0}_0 Z^{\beta_1-1}_1 Z^{\beta_2}_2 Z^{\beta_3}_3]^{II,II}\partial_1\dot{\mathcal{U}}^{II} = \mathcal{G}, \quad (118)$$

where the exact expression of $\mathcal{G}$, with the lower order terms, may be found in [6]. Using (89) again, we write (118), for variable $|\beta| = m-2$, as

$$(\tilde{\mathcal{L}} + \tilde{\mathcal{C}})Z^\beta\partial_1\dot{\mathcal{U}}^{II} = \mathcal{G}, \quad (119)$$

where
$$\tilde{\mathcal{L}} = \begin{pmatrix} \tilde{L} & & \\ & \ddots & \\ & & \tilde{L} \end{pmatrix}$$

with $\tilde{L} = \widehat{\mathcal{A}}_0^{II,II}\partial_t + \sum_{j=1}^n (\widehat{\mathcal{A}}_j + \mathcal{E}_{1j+1})^{II,II}\partial_j$ and where $\tilde{\mathcal{C}} \in H_{*,\gamma}^{s-2}(Q_T^+)$ is a suitable matrix. Here a crucial point is that (119) is a transport-type equation, because the boundary matrix of $\tilde{\mathcal{L}}$ vanishes at $\{x_1 = 0\} \cup \{x_1 = 1\}$. Thus we do not need any boundary condition. Moreover, a standard energy argument gives an $L^2$ a priori estimate for the solution with no loss of regularity w.r.t. the source term $\mathcal{G}$. For its estimate it is important to observe that the only derivatives of $\dot{\mathcal{U}}$ of order $m$ contained in $\mathcal{G}$ are tangential derivatives, estimated in (117). We get the a priori estimate

$$\gamma\|\partial_1\dot{\mathcal{U}}_\gamma^{II}\|_{H_{tan,\gamma}^{m-2}(Q_T)}^2 \leq \frac{C}{\gamma}\left(\|\tilde{f}_\gamma\|_{H_{*,\gamma}^m(Q_T)}^2 + \|\dot{\mathcal{U}}_\gamma\|_{H_{tan,\gamma}^m(Q_T)}^2 + \gamma\|\dot{\mathcal{U}}_\gamma\|_{H_{*,\gamma}^{m-1}(Q_T)}^2\right), \tag{120}$$

for all $\gamma$ sufficiently large, where the constant $C$ depends on $\|\widehat{U}, \nabla\widehat{\Psi}\|_{W^{1,\infty}(Q_T^+)}^2$.

Combining (117), (120) and applying Theorem 15, we infer $\dot{\mathcal{U}} \in H_{tan,\gamma}^m(Q_T^+)$, with $(\dot{q}_\gamma, \dot{u}_{1\gamma}, \dot{h}_{1\gamma})|_{\omega_T} \in H_\gamma^{m-1/2}(\omega_T), \varphi_\gamma \in H_\gamma^{m+1/2}(\omega_T)$ and $\dot{\mathfrak{H}} \in H_{tan,\gamma}^m(Q_T^-)$ with $\nabla\dot{\mathfrak{H}} \in H_{tan,\gamma}^{m-1}(Q_T^-)$.

We also deduce that equation (119) has a unique solution $Z^\beta\partial_1\dot{\mathcal{U}}^{II} \in L^2(Q_T^+)$, for all $|\beta| = m-2$, i.e. $\partial_1\dot{\mathcal{U}}^{II} \in H_{tan,\gamma}^{m-2}(Q_T^+)$. Using (89) again, we infer $\partial_1\dot{\mathcal{U}} \in H_{tan,\gamma}^{m-2}(Q_T^+)$. Adding (117), (120) and taking $\gamma$ sufficiently large yields

$$\gamma\Big(\|\dot{\mathcal{U}}_\gamma\|_{H_{tan,\gamma}^m(Q_T^+)}^2 + \|\partial_1\dot{\mathcal{U}}_\gamma\|_{H_{tan,\gamma}^{m-2}(Q_T^+)}^2 + \|\dot{\mathfrak{H}}_\gamma\|_{H_{tan,\gamma}^m(Q_T^-)}^2 + \|\nabla\dot{\mathfrak{H}}_\gamma\|_{H_{tan,\gamma}^{m-1}(Q_T^-)}^2$$
$$+ \|(\dot{q}_\gamma, \dot{u}_{1\gamma}, \dot{h}_{1\gamma})|_{\omega_T}\|_{H_\gamma^{m-1/2}(\omega_T)}^2 + \|\varphi_\gamma\|_{H_\gamma^{m+1/2}(\omega_T)}^2\Big) \leq \frac{C}{\gamma}\Big\{\gamma\|\dot{\mathcal{U}}_\gamma\|_{H_{*,\gamma}^{m-1}(Q_T)}^2$$
$$+ \Big(\|\dot{\mathcal{U}}_\gamma\|_{H_{*,\gamma}^4(Q_T^+)}^2 + \|\dot{\mathfrak{H}}_\gamma\|_{H_\gamma^4(Q_T^-)}^2 + \|\varphi_\gamma\|_{H_\gamma^4(\omega_T)}^2\Big)\Big(\|\widehat{U}\|_{H_{*,\gamma}^{m+2}(Q_T^+)}^2 + \|\widehat{\mathcal{H}}\|_{H_\gamma^{m+2}(Q_T^-)}^2 + \|\widehat{\Psi}\|_{H_\gamma^{m+3}(Q_T^+)}^2\Big)$$
$$+ \|\tilde{f}_\gamma\|_{H_{*,\gamma}^{m+1}(Q_T^+)}^2 + \|\chi_\gamma, \Xi_\gamma\|_{H_{tan,\gamma}^{m+1}(Q_T^-)}^2 + \|g_{1\gamma}, g_{2\gamma}\|_{H_\gamma^{m+1/2}(\omega_T)}^2$$
$$+ \|g_{3\gamma}\|_{H_\gamma^{m+1}(\omega_T)}^2 + \|g_{4\gamma}\|_{H_\gamma^{m+1/2}(\omega_T)}^2 + \|g_{5\gamma}\|_{H_\gamma^{m+1}(\omega_T^-)}^2\Big\}, \tag{121}$$

where $C$ depends on $K$.

### 8.6. Normal derivatives.

The last step is again by induction, as in [28], page 867, (ii). For convenience of the reader, we provide a brief sketch of the proof.

Suppose that for some fixed $k$, with $1 \leq k < [m/2]$, it has already been shown that $Z^\alpha\partial_1^h\dot{\mathcal{U}}$ belongs to $L^2(Q_T^+)$, for any $h$ and $\alpha$ such that $h = 1, \cdots, k$, $|\alpha| + 2h \leq m$. From (89) it immediately follows that $Z^\alpha\partial_1^{k+1}\dot{\mathcal{U}}^I \in L^2(Q_T^+)$. It rests to prove that $Z^\alpha\partial_1^{k+1}\dot{\mathcal{U}}^{II} \in L^2(Q_T^+)$.

We apply operator $Z^\alpha\partial_1^{k+1}, |\alpha| + 2k = m-2$, to the part $II$ of (85)$_1$ and obtain an equation similar to (119) of the form

$$(\tilde{\mathcal{L}} + \tilde{\mathcal{C}}_k)Z^\alpha\partial_1^{k+1}\dot{\mathcal{U}}^{II} = \mathcal{G}_k, \tag{122}$$

where $\tilde{\mathcal{C}}_k \in H_{*,\gamma}^{s-3}(Q_T^+)$ is a suitable linear operator. The right-hand side $\mathcal{G}_k$ contains derivatives of $\dot{\mathcal{U}}$ of order $m$ (in $H_{*,\gamma}^m$, i.e. counting 1 for each tangential derivative and 2 for normal derivatives), but contains only normal derivatives that have already been estimated. All products of functions are estimated in spaces $H_{*,\gamma}^m$ by the rules given in Theorem 39 and Lemmata 46 and 47 in Appendix A. We infer $\mathcal{G}_k \in L^2(Q_T)$. Again it is crucial that the boundary matrix of $\tilde{\mathcal{L}}$ vanishes at $\{x_1 = 0\} \cup \{x_1 = 1\}$. We obtain the a priori estimate

$$\gamma\|Z^\alpha\partial_1^{k+1}\dot{\mathcal{U}}_\gamma^{II}\|_{L^2(Q_T)}^2 \leq \frac{C}{\gamma}\left(\|\tilde{f}_\gamma\|_{H_{*,\gamma}^m(Q_T)}^2 + \sum_{|\beta|+2h=m,h\leq k}\|Z^\beta\partial_1^h\dot{\mathcal{U}}_\gamma\|_{H_{tan,\gamma}^m(Q_T)}^2 + \gamma\|\dot{\mathcal{U}}_\gamma\|_{H_{*,\gamma}^{m-1}(Q_T)}^2\right), \tag{123}$$

for all $\gamma$ sufficiently large. The solution $Z^\alpha \partial_1^{k+1} \dot{\mathcal{U}}^{II}$ is in $L^2(Q_T^+)$ for all $\alpha, k$ with $|\alpha| + 2k = m - 2$. By repeating this procedure we obtain the result for any $k \leq [m/2]$, hence $\dot{\mathcal{U}} \in H_{*,\gamma}^m(Q_T^+)$. We refer the reader to [6, 22, 24, 28, 30] for similar details.

From (121) and (123) for varying $k$, plus the direct estimate of the normal derivative of $\dot{\mathcal{U}}^I$ by tangential derivatives via (89) we obtain the full regularity of $\dot{\mathcal{U}}_\gamma$ in $H_{*,\gamma}^m(Q_T^+)$:

$$
\begin{aligned}
\gamma \Big( & \|\dot{\mathcal{U}}_\gamma\|_{H_{*,\gamma}^m(Q_T^+)}^2 + \|\dot{\mathfrak{H}}_\gamma\|_{H_{tan,\gamma}^m(Q_T^-)}^2 + \|\nabla \dot{\mathfrak{H}}_\gamma\|_{H_{tan,\gamma}^{m-1}(Q_T^-)}^2 \\
& + \|(\dot{q}_\gamma, \dot{u}_{1\gamma}, \dot{h}_{1\gamma})|_{\omega_T}\|_{H_\gamma^{m-1/2}(\omega_T)}^2 + \|\varphi_\gamma\|_{H_\gamma^{m+1/2}(\omega_T)}^2 \Big) \\
\leq \frac{C}{\gamma} \Big\{ & \Big( \|\dot{\mathcal{U}}_\gamma\|_{H_{*,\gamma}^7(Q_T^+)}^2 + \|\dot{\mathfrak{H}}_\gamma\|_{H_\gamma^4(Q_T^-)}^2 + \|\varphi_\gamma\|_{H_\gamma^3(\omega_T)} \Big) \Big( \|\widehat{U}\|_{H_{*,\gamma}^{m+2}(Q_T^+)} + \|\widehat{\mathcal{H}}\|_{H_\gamma^{m+2}(Q_T^-)} + \|\widehat{\Psi}\|_{H_\gamma^{m+3}(Q_T)} \Big) \\
& + \|\tilde{f}_\gamma\|_{H_{*,\gamma}^{m+1}(Q_T^+)}^2 + \|\chi_\gamma, \Xi_\gamma\|_{H_{tan,\gamma}^{m+1}(Q_T^-)} + \|g_{1\gamma}, g_{2\gamma}\|_{H_\gamma^{m+1/2}(\omega_T)} \\
& + \|g_{3\gamma}\|_{H_\gamma^{m+1}(\omega_T)} + \|g_{4\gamma}\|_{H_\gamma^{m+1/2}(\omega_T^+)} + \|g_{5\gamma}\|_{H_\gamma^{m+1}(\omega_T^-)} \Big\}, \quad (124)
\end{aligned}
$$

where $C$ depends on $K$. The norms of higher order in the right-hand side are absorbed by taking $\gamma$ sufficiently large.

### 8.7. Regularity of the vacuum magnetic field.
Up to now we have only proved for the vacuum magnetic field $\dot{\mathfrak{H}} \in H_{tan,\gamma}^m(Q_T^-)$ with $\nabla \dot{\mathfrak{H}} \in H_{tan,\gamma}^{m-1}(Q_T^-)$. We write the normal derivatives in terms of the tangential derivatives from

$$
\nabla \times \dot{\mathfrak{H}} = \chi, \quad \operatorname{div}(A\dot{\mathfrak{H}}) = \Xi \quad \text{in } Q_T^-,
$$

where $(\chi_\gamma, \Xi_\gamma) \in H_\gamma^{m+1}(Q_T^-)$. Starting from $\nabla \dot{\mathfrak{H}} \in H_{tan,\gamma}^{m-1}(Q_T^-)$ we may increase the regularity in the normal direction step by step and finally conclude $\dot{\mathfrak{H}} \in H_\gamma^m(Q_T^-)$.

### 8.8. The tame estimate.
From (124), after the previous calculations for the additional regularity of $\dot{\mathfrak{H}}_\gamma$ we obtain, for all $\gamma \geq \gamma_m$ sufficiently large, the following estimate

$$
\begin{aligned}
\gamma \Big( & \|\dot{\mathcal{U}}_\gamma\|_{H_{*,\gamma}^m(Q_T^+)}^2 + \|\dot{\mathfrak{H}}_\gamma\|_{H_\gamma^m(Q_T^-)}^2 + \|(\dot{q}_\gamma, \dot{u}_{1\gamma}, \dot{h}_{1\gamma})|_{\omega_T}\|_{H_\gamma^{m-1/2}(\omega_T)}^2 + \|\varphi_\gamma\|_{H_\gamma^{m+1/2}(\omega_T)}^2 \Big) \\
\leq \frac{C}{\gamma} \Big\{ & \Big( \|\dot{\mathcal{U}}_\gamma\|_{H_{*,\gamma}^7(Q_T^+)}^2 + \|\dot{\mathfrak{H}}_\gamma\|_{H_\gamma^4(Q_T^-)}^2 + \|\varphi_\gamma\|_{H_\gamma^3(\omega_T)} \Big) \Big( \|\widehat{U}\|_{H_{*,\gamma}^{m+2}(Q_T^+)} + \|\widehat{\mathcal{H}}\|_{H_\gamma^{m+2}(Q_T^-)} + \|\widehat{\Psi}\|_{H_\gamma^{m+3}(Q_T)} \Big) \\
& + \|f_\gamma\|_{H_{*,\gamma}^{m+1}(Q_T^+)}^2 + \|\chi_\gamma, \Xi_\gamma\|_{H_\gamma^{m+1}(Q_T^-)} + \|g_{1\gamma}, g_{2\gamma}\|_{H_\gamma^{m+1/2}(\omega_T)} \\
& + \|g_{3\gamma}\|_{H_\gamma^{m+1}(\omega_T)} + \|g_{4\gamma}\|_{H_\gamma^{m+1/2}(\omega_T^+)} + \|g_{5\gamma}\|_{H_\gamma^{m+1}(\omega_T^-)} \Big\}, \quad (125)
\end{aligned}
$$

where $C$ depends on $K$.

The a priori estimate (125) holds for all $m \geq 1$ and $\gamma \geq \gamma_m$ (we may assume that $\gamma_m$ is an increasing function of $m$). From (125) for $m = 7$ and $\gamma \geq \gamma_7$ one gets

$$
\begin{aligned}
\gamma \Big( & \|\dot{\mathcal{U}}_\gamma\|_{H_{*,\gamma}^7(Q_T^+)}^2 + \|\dot{\mathfrak{H}}_\gamma\|_{H_\gamma^7(Q_T^-)}^2 + \|\varphi_\gamma\|_{H_\gamma^{7.5}(\omega_T)} \Big) \\
\leq \frac{C}{\gamma} \Big\{ & \Big( \|\dot{\mathcal{U}}_\gamma\|_{H_{*,\gamma}^7(Q_T^+)}^2 + \|\dot{\mathfrak{H}}_\gamma\|_{H_\gamma^4(Q_T^-)}^2 + \|\varphi_\gamma\|_{H_\gamma^3(\omega_T)} \Big) \Big( \|\widehat{U}\|_{H_{*,\gamma}^9(Q_T^+)} + \|\widehat{\mathcal{H}}\|_{H_\gamma^9(Q_T^-)} + \|\widehat{\Psi}\|_{H_\gamma^{10}(Q_T)} \Big) \\
& + \|f_\gamma\|_{H_{*,\gamma}^8(Q_T^+)}^2 + \|\chi_\gamma, \Xi_\gamma\|_{H_\gamma^8(Q_T^-)} + \|g_{1\gamma}, g_{2\gamma}\|_{H_\gamma^{7.5}(\omega_T)} \\
& + \|g_{3\gamma}\|_{H_\gamma^8(\omega_T)} + \|g_{4\gamma}\|_{H_\gamma^{7.5}(\omega_T^+)} + \|g_{5\gamma}\|_{H_\gamma^8(\omega_T^-)} \Big\}.
\end{aligned}
$$

Then, by taking $\gamma$ larger if needed, one obtains

$$\gamma\Big(\|\dot{\mathcal{U}}_\gamma\|^2_{H^7_{*,\gamma}(Q^+_T)} + \|\dot{\mathfrak{H}}_\gamma\|^2_{H^7_\gamma(Q^-_T)} + \|\varphi_\gamma\|^2_{H^{7.5}_\gamma(\omega_T)}\Big)$$
$$\leq \frac{C}{\gamma}\Big\{\|f_\gamma\|^2_{H^8_{*,\gamma}(Q^+_T)} + \|\chi_\gamma, \Xi_\gamma\|_{H^8_\gamma(Q^-_T)} + \|g_{1\gamma}, g_{2\gamma}\|_{H^{7.5}_\gamma(\omega_T)}$$
$$+ \|g_{3\gamma}\|_{H^8_\gamma(\omega_T)} + \|g_{4\gamma}\|_{H^{7.5}_\gamma(\omega^+_T)} + \|g_{5\gamma}\|_{H^8_\gamma(\omega^-_T)}\Big\}, \quad (126)$$

where $C$ depends on $K$, having used $\|\widehat{\Psi}\|_{H^{10}_\gamma(Q_T)} \leq C\|\hat{\varphi}\|^2_{H^{9.5}_\gamma(\omega_T)}$. At last, substituting (126) in (125) gives the tame estimate

$$\gamma\Big(\|\dot{\mathcal{U}}_\gamma\|^2_{H^m_{*,\gamma}(Q^+_T)} + \|\dot{\mathfrak{H}}_\gamma\|^2_{H^m_\gamma(Q^-_T)} + \|(\dot{q}_\gamma, \dot{u}_{1\gamma}, \dot{h}_{1\gamma})|_{\omega_T}\|^2_{H^{m-1/2}_\gamma(\omega_T)} + \|\varphi_\gamma\|^2_{H^{m+1/2}_\gamma(\omega_T)}\Big)$$
$$\leq \frac{C}{\gamma}\Big\{\Big(\|f_\gamma\|^2_{H^8_{*,\gamma}(Q^+_T)} + \|\chi_\gamma, \Xi_\gamma\|_{H^8_\gamma(Q^-_T)} + \|g_{1\gamma}, g_{2\gamma}\|_{H^{7.5}_\gamma(\omega_T)} + \|g_{3\gamma}\|_{H^8_\gamma(\omega_T)}$$
$$+ \|g_{4\gamma}\|_{H^{7.5}_\gamma(\omega^+_T)} + \|g_{5\gamma}\|_{H^8_\gamma(\omega^-_T)}\Big)\Big(\|\widehat{U}\|_{H^{m+2}_{*,\gamma}(Q^+_T)} + \|\widehat{\mathcal{H}}\|_{H^{m+2}_\gamma(Q^-_T)} + \|\hat{\varphi}\|_{H^{m+2.5}_\gamma(\omega_T)}\Big)$$
$$+ \|f_\gamma\|^2_{H^{m+1}_{*,\gamma}(Q^+_T)} + \|\chi_\gamma, \Xi_\gamma\|_{H^{m+1}_\gamma(Q^-_T)} + \|g_{1\gamma}, g_{2\gamma}\|_{H^{m+1/2}_\gamma(\omega_T)}$$
$$+ \|g_{3\gamma}\|_{H^{m+1}_\gamma(\omega_T)} + \|g_{4\gamma}\|_{H^{m+1/2}_\gamma(\omega^+_T)} + \|g_{5\gamma}\|_{H^{m+1}_\gamma(\omega^-_T)}\Big\}, \quad (127)$$

where $C$ depends on $K$. In the end, for homogeneity of notation we raise the norms of boundary data of fractional order to the higher integer order and get (88). The proof of Theorem 16 is complete.

## 9. Compatibility conditions on the initial data

Assume we are given initial data $U_0 = (q_0, v^0, H^0, S_0)$, $\mathcal{H}^0$ and $\varphi_0$ that satisfy the hyperbolicity condition (8) and the stability condition (22). We also assume

$$\|\varphi_0\|_{H^{2.5}(\Gamma)} \leq \epsilon_1, \quad (128)$$

with $\epsilon_1 := \epsilon_0/2$, for $\epsilon_0$ as in Lemma 3. Let the functions $\Psi_0, \Phi^0$ be defined from $\varphi_0$, as in Lemmata 2 and 3. We assume also that the initial plasma magnetic field $H^0$ satisfies

$$\mathrm{div}\, h^0 = 0 \quad \text{in } \Omega^+,$$
$$H^0_{N_0} = 0 \quad \text{on } \Gamma, \qquad H^0_1 = 0 \quad \text{on } \Gamma_+, \quad (129)$$

where $h^0 = (H^0_{N_0}, H^0_2 \partial_1 \Phi^0, H^0_3 \partial_1 \Phi^0)$, $H^0_{N_0} = H^0_1 - H^0_2 \partial_2 \Psi_0 - H^0_3 \partial_3 \Psi_0$, and the initial vacuum magnetic field $\mathcal{H}^0$ satisfies

$$\nabla \times \mathfrak{H}^0 = 0, \quad \mathrm{div}\, \mathfrak{h}^0 = 0 \quad \text{in } \Omega^-,$$
$$\mathcal{H}^0_{N_0} = 0 \quad \text{on } \Gamma, \qquad \nu \times \mathcal{H}^0 = \mathfrak{J}(0) \quad \text{on } \Gamma_-, \quad (130)$$

where $\mathfrak{H}^0$, $\mathfrak{h}^0$ and $\mathcal{H}^0_{N_0}$ are defined by

$$\mathfrak{H}^j = (\mathcal{H}^j_1 \partial_1 \Phi^1, \mathcal{H}^j_{\tau_2}, \mathcal{H}^j_{\tau_3}), \quad \mathfrak{h}^j = (\mathcal{H}^j_N, \mathcal{H}^j_2 \partial_1 \Phi^1, \mathcal{H}^j_3 \partial_1 \Phi^0),$$
$$\mathcal{H}^j_N = \mathcal{H}^j_1 - \mathcal{H}^j_2 \partial_2 \Psi_0 - \mathcal{H}^j_3 \partial_3 \Psi_0, \quad \mathcal{H}^j_{\tau_i} = \mathcal{H}^j_i \partial_i \Psi_0 + \mathcal{H}^j_i, \quad i = 2, 3, \quad (131)$$

for $j = 0$. Notice that system (130) uniquely determines $\mathcal{H}^0$ from $\varphi_0$ and $\mathfrak{J}(0)$ by Theorem 13, see the comment in Remark 6.

Let us define $U_j = (q_j, v^j, H^j, S_j)$, with $v^j = (v^j_1, v^j_2, v^j_3)$ and $H^j = (H^j_1, H^j_2, H^j_3)$, and $\varphi_j$ by formally taking $j - 1$ time derivatives of (17) and the boundary equation $\partial_t \varphi - v_N = 0$, evaluating at time $t = 0$ and solving for $\partial^j_t U(0), \partial^j_t \varphi(0)$. This procedure inductively determines $\partial^j_t U(0), \partial^j_t \varphi(0)$ in terms of $U_0, \varphi_0$. We denote $U_j = \partial^j_t U(0), \varphi_j = \partial^j_t \varphi(0)$. Corresponding to $\varphi_j$ we compute the functions $\Psi_j, \Phi^j$, as in Lemmata 2 and 3.

Finally, we define the time derivatives at initial time $\mathcal{H}^j$ as the unique solution of the elliptic system

$$\nabla \times \mathfrak{H}^j = \alpha_j, \quad \mathrm{div}\, \mathfrak{h}^j = \beta_j \quad \text{in } \Omega^-,$$
$$\mathcal{H}^j_N = \vartheta_j, \quad \text{on } \Gamma, \qquad \nu \times \mathcal{H}^j = \partial^j_t \mathfrak{J}(0), \quad \text{on } \Gamma_-, \quad (132)$$

where $\mathfrak{H}^j$, $\mathfrak{h}^j$ and $\mathcal{H}_N^j$ are as in (131) and $\alpha_j, \beta_j, \theta_j$ are suitable commutators, for example

$$\vartheta_j = [^{(j)}, \partial_2 \Psi_0] \mathcal{H}_2^0 + [^{(j)}, \partial_3 \Psi_0] \mathcal{H}_3^0$$

with $^{(j)}\Psi_0 := \Psi_j$, $^{(j)}\mathcal{H}_k^0 := \mathcal{H}_k^j$. From the second boundary equation in (18), stating the continuity of the total pressure, we deduce that sufficiently regular solutions should satisfy

$$\partial_t^j (q - \tfrac{1}{2} |\mathcal{H}|^2)\big|_{t=0} = 0 \quad \text{on } \Gamma.$$

These equations yield the compatibility conditions

$$\begin{aligned}
q_0 &= \frac{1}{2} |\mathcal{H}^0|^2 \quad \text{on } \Gamma, \quad j = 0, \\
q_j &= \sum_{i=0}^{j-1} C_{j-1}^i (\mathcal{H}^i, \mathcal{H}^{j-i}) \quad \text{on } \Gamma, \quad j \geq 1.
\end{aligned} \tag{133}$$

Notice that the other boundary conditions in (18) do not give raise to compatibility conditions as these are implicitly included in the above definitions of $\varphi_j, \mathcal{H}^j$.

**Lemma 19.** *Let $k \in \mathbb{N}$, $k \geq 4$, $U_0 \in H^{k-1/2}(\Omega^+)$, $\mathcal{H}^0 \in H^{k-1/2}(\Omega^-)$, $\varphi_0 \in H^k(\Gamma)$ and $\mathfrak{J} \in H^{k-1/2}([0, T_0] \times \Gamma_-)$. Then, the procedure described above determines $U_j \in H^{k-j-1/2}(\Omega^+)$, $\mathcal{H}^j \in H^{k-j-1/2}(\Omega^-)$ and $\varphi_j \in H^{k-j}(\Gamma)$ for $j = 1, \ldots, k-1$. Moreover,*

$$\|\mathcal{H}^0\|_{H^{k-1/2}(\Omega^-)} + \sum_{j=1}^{k-1} \left( \|U_j\|_{H^{k-j-1/2}(\Omega^+)} + \|\mathcal{H}^j\|_{H^{k-j-1/2}(\Omega^-)} + \|\varphi_j\|_{H^{k-j}(\Gamma)} \right) \leq C(M_0), \tag{134}$$

*where the constant $C = C(M_0) > 0$ depends on*

$$M_0 = \|U_0\|_{H^{k-1/2}(\Omega^+)} + \|\varphi_0\|_{H^k(\Gamma)} + \sum_{j=0}^{k-1} \|\partial_t^j \mathfrak{J}(0)\|_{H^{k-j-1}(\Gamma_-)}. \tag{135}$$

*Proof.* As follows from the construction of the functions $U_j$ and $\varphi_j$, we can estimate them separately from $\mathcal{H}^j$:

$$\sum_{j=1}^{k-1} \left( \|U_j\|_{H^{k-j-1/2}(\Omega^+)} + \|\varphi_j\|_{H^{k-j}(\Gamma)} \right) \leq CM_0, \tag{136}$$

where the constant $C > 0$ depends only on $k$ and the norms $\|U_0\|_{W^{1,\infty}(\Omega^+)}$ and $\|\varphi_0\|_{W^{1,\infty}(\Gamma)}$. For the proof of (136) we refer to [27, 21].

Problem (130) has the form of problem (63) with $\chi = 0$, $\Xi = 0$, $g_3 = 0$, and Theorem 13 shows its solvability with $H^2$ regularity of the solution. By classical results one can improve the regularity to

$$\|\mathfrak{H}^0\|_{H^{k-1/2}(\Omega^-)} \leq C(K_1) \left( \|\varphi_0\|_{H^k(\Gamma)} + \|\mathfrak{J}(0)\|_{H^{k-1}(\Gamma_-)} \right),$$

with $K_1 = \|\varphi_0\|_{W^{2,\infty}(\Gamma)}$. Since $H^k(\Gamma) \hookrightarrow W^{2,\infty}(\Gamma)$, we get

$$\|\mathcal{H}^0\|_{H^{k-1/2}(\Omega^-)} \leq C(M_0). \tag{137}$$

For estimating $\mathcal{H}^j$ we can use as well regularity results for elliptic systems. Indeed, the elliptic problem (132) has the form of problem (63) with

$$\chi = \alpha_j, \quad \Xi = \beta_j, \quad g_3 = \vartheta_j, \tag{138}$$

and $\hat{\varphi}$, $\widehat{\Psi}$ and $\widehat{\Phi}_1$ replaced by $\varphi_0$, $\Psi_0$ and $\Phi_1^0$ respectively. The proof of estimate (134) follows then by finite induction with respect to the upper limit of the sum. Applying the Moser-type calculus inequality (269) for estimating the commutators $\alpha_1$, $\beta_1$ and $\vartheta_1$, taking into account Lemmata 2 and 3 and using estimates (136), (137) and Sobolev's imbeddings, for problem (132) with $j = 1$ we derive the a priori estimate

$$\|\mathcal{H}^1\|_{H^{k-3/2}(\Omega^-)} \leq C(M_0)$$

justifying the basis for the induction. Exploiting similar arguments, from the inductive hypothesis

$$\|\mathcal{H}^0\|_{H^{k-1/2}(\Omega^-)} + \sum_{j=1}^{k-2} \left( \|U_j\|_{H^{k-j-1/2}(\Omega^+)} + \|\mathcal{H}^j\|_{H^{k-j-1/2}(\Omega^-)} + \|\varphi_j\|_{H^{k-j}(\Gamma)} \right) \leq C(M_0)$$

we derive the desired estimate (134). But, in this step some terms appearing in the commutators are treated in a different way. For example, for the term $\mathcal{H}_2^0|_\Gamma \partial_2 \varphi_{k-1}$ appearing in $\vartheta_{k-1}$ we do not use the Moser-type inequality and estimate it as follows:

$$\|\mathcal{H}_2^0 \partial_2 \varphi_{k-1}\|_{H^{1/2}(\Gamma)} \leq \|\mathcal{H}_2^0\|_{L^\infty(\Gamma)} \|\varphi_{k-1}\|_{H^{3/2}(\Gamma)} \leq \|\mathcal{H}^0\|_{H^k(\Gamma)} \|\varphi_{k-1}\|_{H^{3/2}(\Gamma)} \leq C(M_0),$$

where we have used Sobolev's imbedding and estimates (136) and (137).      $\square$

**Definition 20.** *Let $k \in \mathbb{N}$, $k \geq 4$. Consider initial data $U_0 \in H^{k-1/2}(\Omega^+)$, $\mathcal{H}^0 \in H^{k-1/2}(\Omega^-)$, and $\varphi_0 \in H^k(\Gamma)$ that satisfy* (8), (22), (128), (129) *and* (130). *The initial data $(U_0, \mathcal{H}^0, \varphi_0)$ are said to be compatible up to order $k-1$ if they satisfy* (133) *on $\Gamma$, $v_1^j = 0$ on $\Gamma_+$, for $j = 0, \ldots, k-2$, and*

$$\int_\Gamma \left| q_{k-1} - \sum_{i=0}^{k-2} C_{k-2}^i (\mathcal{H}^i, \mathcal{H}^{k-1-i}) \right|^2 \frac{dx_1}{x_1} dx' + \int_{\Gamma_+} |v_1^{k-1}|^2 \frac{dx_1}{x_1} dx' < +\infty. \tag{139}$$

Observe that $U_j \in H^{3/2}(\Omega^+)$, $\mathcal{H}^j \in H^{3/2}(\Omega^-)$ for $j = 0, \ldots, k-2$, so it is legitimate to take the traces at $\{x_1 = 0\}$.

## 10. CONSTRUCTION OF AN APPROXIMATE SOLUTION

We now introduce the following "approximate" solution. As regards the plasma equations, these are solutions in the sense of Taylor's series at $t = 0$. Let us set

$$Q = \mathbb{R} \times \Omega, \quad Q^\pm = \mathbb{R} \times \Omega^\pm, \quad \omega = \mathbb{R} \times \Gamma, \quad \omega^\pm = \mathbb{R} \times \Gamma_\pm.$$

First we extend the density current $\mathfrak{J} \in H^{k-1/2}([0, T_0] \times \Gamma_-)$ to the whole real line of times by preserving the same regularity. Thus from now on we assume that $\mathfrak{J} \in H^{k-1/2}(\omega^-)$.

**Lemma 21.** *Let $k \in \mathbb{N}$, $k \geq 4$ and $\mathfrak{J} \in H^{k-1/2}(\omega^-)$. Consider initial data $U_0 \in H^{k-1/2}(\omega^+)$, $\mathcal{H}^0 \in H^{k-1/2}(\Omega^-)$, and $\varphi_0 \in H^k(\Gamma)$ that satisfy* (8), (22), (128), (129) *and* (130) *and are compatible up to order $k-1$ in the sense of Definition 20. Then there exist functions $(U^a, \mathcal{H}^a, \varphi^a)$ such that $U^a \in H^k(Q^+)$, $\mathcal{H}^a \in H^k(Q^-)$, $\varphi^a \in H^{k+1/2}(\omega)$, and such that*

$$\partial_t^j \mathbb{P}(U^a, \Psi^a)_{|t=0} = 0 \qquad \text{for } j = 0, \ldots, k-1, \tag{140}$$

$$\mathbb{V}(\mathcal{H}^a, \Psi^a) = 0 \qquad \text{in } Q^-, \tag{141}$$

$$\mathbb{B}(U^a, \mathcal{H}^a, \varphi^a) = \bar{\mathfrak{J}} \qquad \text{on } \omega^3 \times \omega^+ \times \omega^-, \tag{142}$$

*where $\Psi^a$ is constructed from $\varphi^a$ by Lemma 2, and $\bar{\mathfrak{J}} = (0, 0, 0, 0, \mathfrak{J})^T$. Moreover $(U^a, \mathcal{H}^a, \varphi^a)$ satisfy* $(29)_2$, (30), (61) *with a strict inequality, the constraint* (31) *and the estimate*

$$\|U^a\|_{H^k(Q^+)} + \|\mathcal{H}^a\|_{H^k(Q^-)} + \|\varphi^a\|_{H^{k+1/2}(\omega)} \leq C(M_0), \tag{143}$$

*with $C = C(M_0) > 0$ and $M_0$ defined in* (135).

*Proof.* Given the initial data, let us take $U_j = (q_j, v^j, H^j, S_j)$ and $\varphi_j$, with $v^j = (v_1^j, v_2^j, v_3^j)$ and $H^j = (H_1^j, H_2^j, H_3^j)$, as in Lemma 19. We first take $(v^a, S^a) \in H^k(Q^+)$ such that

$$\partial_t^j (v^a, S^a)_{|t=0} = (v^j, S_j) \quad \text{in } \Omega^+, \, j = 0, \ldots, k-1, \qquad v_1^a = 0 \quad \text{on } \omega^+. \tag{144}$$

Given $v^a$, we find $\varphi^a$ from

$$\begin{cases} \partial_t \varphi^a = v_{N^a}^a := v_1^a - v_2^a \partial_2 \varphi^a - v_3^a \partial_3 \varphi^a & \text{on } \omega, \\ \varphi^a_{|t=0} = \varphi_0. \end{cases} \tag{145}$$

As $v^a|_\omega \in H^{k-1/2}(\omega)$ we get $\varphi^a \in H^{k-1/2}(\omega)$. From (144), deriving (145) in time, it follows

$$\partial_t^j \varphi^a_{|t=0} = \varphi_j, \qquad j = 0, \ldots, k-1.$$

Since $\varphi^a$ satisfies (128) at the initial time $t = 0$, by a cut-off argument we can choose $\varphi^a$ that satisfies $(29)_2$ for all times in strict sense. From $\varphi^a$ we compute $\Psi^a$, $\Phi^a$ as in Lemmata 2, 3. Now we solve

$$\partial_t H^a + \frac{1}{\partial_1 \Phi_1^a} \left\{ (w^a \cdot \nabla) H^a - (h^a \cdot \nabla) v^a + H^a \mathrm{div}\, u^a \right\} = 0 \qquad \text{in } Q_T^+, \tag{146}$$
$$H^a|_{t=0} = H_0,$$

where

$$h^a = (H_1^a - H_2^a \partial_2 \Psi^a - H_3^a \partial_3 \Psi^a, H_2^a \partial_1 \Phi_1^a, H_3^a \partial_1 \Phi_1^a),$$

$$u^a = (v_{N^a}^a, v_2^a \partial_1 \Phi_1^a, v_3^a \partial_1 \Phi_1^a), \quad w^a = u^a - (\partial_t \Psi^a, 0, 0).$$

We have $w_1^a = 0$ on $\Gamma \cup \Gamma_+$, so that (146) doesn't need any boundary condition.

With the usual calculations we find from (146) and the initial constraints (129), that $H^a$ satisfies for $t > 0$

$$\mathrm{div}\, h(t) = 0 \quad \text{in } \Omega^+, \tag{147}$$
$$H_N(t) = 0 \quad \text{on } \Gamma, \qquad H_1(t) = 0 \quad \text{on } \Gamma_+.$$

From (146) we also obtain

$$\partial_t^j H^a|_{t=0} = H^j, \qquad j = 0, \ldots, k-1.$$

Then we compute the vacuum magnetic field $\mathcal{H}^a$, for each fixed $t$, from the problem

$$\nabla \times \mathfrak{H}^a = 0, \quad \mathrm{div}\, \mathfrak{h}^a = 0 \quad \text{in } \Omega^-, \tag{148}$$
$$\mathcal{H}_{N^a}^a = 0 \quad \text{on } \Gamma, \qquad \nu \times \mathcal{H}^a = \mathfrak{J} \quad \text{on } \Gamma_-,$$

where $\mathcal{H}_{N^a}^a = \mathcal{H}_1^a - \mathcal{H}_2^a \partial_2 \varphi^a - \mathcal{H}_3^a \partial_3 \varphi^a$, and where $\mathfrak{H}^a, \mathfrak{h}^a$ are defined as usual from $\mathcal{H}^a, \Psi^a, \Phi^a$. When $t = 0$ we get

$$\nabla \times \mathfrak{H}^a(0) = 0, \quad \mathrm{div}\, \mathfrak{h}^a(0) = 0 \quad \text{in } \Omega^-, \tag{149}$$
$$\mathcal{H}_{N^a}^a(0) = 0 \quad \text{on } \Gamma, \qquad \nu \times \mathcal{H}^a(0) = \mathfrak{J}(0) \quad \text{on } \Gamma_-.$$

On the other hand we are prescribing for $\mathcal{H}_0$ the initial constraints (130). As the right hand side of (130) and (149) is the same, by uniqueness of Theorem 13 we get $\mathfrak{H}^a(0) = \mathfrak{H}^0$, that is $\mathcal{H}^a(0) = \mathcal{H}^0$. By taking the time derivatives of (148) we obtain

$$\partial_t^j \mathcal{H}^a|_{t=0} = \mathcal{H}^j, \qquad j = 0, \ldots, k-1.$$

To conclude, we define $q^a \in H^k(Q^+)$ by requiring

$$\partial_t^j q^a|_{t=0} = q_j \quad \text{in } \Omega^+, \quad j = 0, \ldots, k-1, \tag{150}$$
$$q^a = \frac{1}{2}|\mathcal{H}^a|^2 \quad \text{on } \omega.$$

Such a lifting is possible because of the compatibility conditions (133) for $j = 0, \ldots, k-2$ and (139), see [20, Theorem 2.3].

The regularity $\varphi^a \in H^{k+1/2}(\omega)$ follows from $(147)_2$, $(148)_2$ and the stability condition (22). We observe that $(148)_1$ may be equivalently written as (141), and $(144)_2$, $(145)_1$, $(148)_2$, $(150)_2$ as (142).

At last, since $U^a, \mathcal{H}^a$ satisfy (30), (61) at the initial time $t = 0$, by a cut-off argument in the above procedure we can choose $U^a, \mathcal{H}^a$ that satisfy (30), (61) for all times in strict sense. $\qquad \square$

The approximate solution $(U^a, \mathcal{H}^a, \varphi^a)$ enables us to reformulate the original problem as a nonlinear problem with zero initial data. Let us take $k = m + 10$ in Lemmata 19, 21, where $m \in \mathbb{N}$. Introduce:

$$\begin{cases} f^a := -\mathbb{P}(U^a, \Psi^a), & t > 0, \\ f^a := 0, & t < 0. \end{cases} \tag{151}$$

Because $U^a \in H^{m+10}(Q^+)$, $\varphi^a \in H^{m+10.5}(\omega)$, (140) yields $f^a \in H^{m+9}(Q^+)$. From (143), we also get the estimate:

$$\|f^a\|_{H^{m+9}(Q^+)} \leq C(M_0). \tag{152}$$

Given the approximate solution $(U^a, \mathcal{H}^a, \varphi^a)$ of Lemma 21, and $f^a$ defined by (151), then $(U, \mathcal{H}, \varphi) = (U^a, \mathcal{H}^a, \varphi^a) + (V, \mathcal{K}, \psi)$ is a solution of (17)–(19) on $Q_T^+ \times Q_T^-$, if $(V, \mathcal{K}, \psi)$ satisfies the following system:

$$
\begin{cases}
\mathcal{L}(V, \Psi) = f^a, & \text{in } Q_T^+, \\
\mathcal{E}(\mathcal{K}, \Psi) = 0, & \text{in } Q_T^-, \\
\mathcal{B}(V, \mathcal{K}, \psi) = 0, & \text{on } \omega_T \times \omega_T^\pm, \\
(V, \mathcal{K}, \psi) = 0, & \text{for } t < 0,
\end{cases}
\tag{153}
$$

where

$$
\begin{aligned}
\mathcal{L}(V, \Psi) &:= \mathbb{P}(U^a + V, \Psi^a + \Psi) - \mathbb{P}(U^a, \Psi^a), \\
\mathcal{E}(\mathcal{K}, \Psi) &:= \mathbb{V}(\mathcal{H}^a + \mathcal{K}, \Psi^a + \Psi), \\
\mathcal{B}(V, \mathcal{K}, \psi) &:= \mathbb{B}(U^a + V, \mathcal{H}^a + \mathcal{K}, \varphi^a + \psi) - \tilde{\mathfrak{I}},
\end{aligned}
\tag{154}
$$

with $\Psi$ denoting the extension constructed from $\psi$ by Lemma 2. We note that the properties of the approximate solution imply that $(V, \mathcal{K}, \psi) = 0$ satisfy (153) for $t < 0$. Therefore the initial nonlinear problem on $[0, T] \times \Omega^\pm$ is now substituted for a problem on $Q_T^+ \times Q_T^-$. The initial data (19) are absorbed into the source term, and the problem has to be solved in the class of functions vanishing in the past, which is exactly the class of functions in which we have a well-posedness result for the linearized problem. Thanks to (143), we see that $(U^a, \mathcal{H}^a, \varphi^a)$ satisfies the first inequality in (29) if we choose $K \geq C(M_0)$.

## 11. Description of the iterative scheme

We solve problem (153) by a Nash-Moser type iteration. (We refer to [2, 15, 16] for a general description of the method). This method requires a family of smoothing operators whose construction is inspired from [1, 13], see also [8].

### 11.1. The smoothing operators.

We begin with a few notations. In what follows, $\Omega_T$ stands alternatively for $Q_T^+, Q_T^-, Q_T$. For $T > 0$, $s \geq 0$, and $\gamma \geq 1$, we let

$$
\mathcal{F}_\gamma^s(\Omega_T) := \left\{ u \in H_\gamma^s(\Omega_T), \quad u = 0 \text{ for } t < 0 \right\}.
$$

This is a closed subspace of $H_\gamma^s(\Omega_T)$, that we equip with the induced norm. In case of $Q_T^+$, in the definition of $\mathcal{F}_\gamma^s$ the space $H_\gamma^s(Q_T^+)$ is substituted by $H_{*,\gamma}^s(Q_T^+)$. The definition of $\mathcal{F}_\gamma^s(\omega_T)$ is entirely similar.

**Lemma 22.** *There exists a family* $\{S_\theta\}_{\theta \geq 1}$ *of operators* $S_\theta : \mathcal{F}_\gamma^0(\Omega_T) \longrightarrow \bigcap_{\beta \geq 0} \mathcal{F}_\gamma^\beta(\Omega_T)$, *such that*

$$
\|S_\theta u\|_{H_\gamma^\beta(\Omega_T)} \leq C \, \theta^{(\beta - \alpha)_+} \|u\|_{H_\gamma^\alpha(\Omega_T)} \qquad \forall \alpha, \beta \geq 0,
\tag{155a}
$$

$$
\|S_\theta u - u\|_{H_\gamma^\beta(\Omega_T)} \leq C \, \theta^{\beta - \alpha} \|u\|_{H_\gamma^\alpha(\Omega_T)} \qquad 0 \leq \beta \leq \alpha,
\tag{155b}
$$

$$
\left\| \frac{d}{d\theta} S_\theta u \right\|_{H_\gamma^\beta(\Omega_T)} \leq C \, \theta^{\beta - \alpha - 1} \|u\|_{H_\gamma^\alpha(\Omega_T)} \qquad \forall \alpha, \beta \geq 0.
\tag{155c}
$$

*Here we use the classical notation* $(\beta - \alpha)_+ := \max(0, \beta - \alpha)$. *The constants in the inequalities are uniform with respect to* $\alpha, \beta$ *when* $\alpha, \beta$ *belong to some bounded interval. In case of* $Q_T^+$, *in* (155) *the norm of* $H_{*,\gamma}^s(Q_T^+)$ *substitutes the norm of* $H_\gamma^s(Q_T^+)$, $s = \alpha, \beta$.

*Moreover, there is another family of operators, still denoted* $S_\theta$, *that acts on functions that are defined on the boundary* $\omega_T$, *and that enjoy the properties* (155), *with the norms* $\| \cdot \|_{H_\gamma^s(\omega_T)}$.

### 11.2. Description of the iterative scheme.

Let us describe the iterative scheme. The scheme starts from $V_0 = 0, \mathcal{K}_0 = 0, \Psi_0 = 0, \psi_0 = 0$. Assume that $V_k, \mathcal{K}_k, \Psi_k, \psi_k$ are already given for $k = 0, \ldots, n$ and verify

$$
\begin{cases}
V_{k,2} = 0 \quad \text{on } \omega_T^+, \qquad \nu \times \mathcal{K}_k = 0 \quad \text{on } \omega_T^-, \\
(V_k, \mathcal{K}_k, \Psi_k, \psi_k) = 0 \quad \text{for } t < 0.
\end{cases}
\tag{156}
$$

As in [2], we consider

$$
\begin{aligned}
V_{n+1} &= V_n + \delta V_n, & \mathcal{K}_{n+1} &= \mathcal{K}_n + \delta \mathcal{K}_n, \\
\Psi_{n+1} &= \Psi_n + \delta \Psi_n, & \psi_{n+1} &= \psi_n + \delta \psi_n,
\end{aligned}
\tag{157}
$$

where the differences $\delta V_n, \delta \mathcal{K}_n, \delta \Psi_n, \delta \psi_n$ will be specified later on. Given $\theta_0 \geq 1$, let us set $\theta_n := (\theta_0^2 + n)^{1/2}$, and consider the smoothing operators $S_{\theta_n}$. We decompose

$$
\begin{aligned}
\mathcal{L}(V_{n+1}, \Psi_{n+1}) - \mathcal{L}(V_n, \Psi_n) &= \mathbb{P}(U^a + V_{n+1}, \Psi^a + \Psi_{n+1}) - \mathbb{P}(U^a + V_n, \Psi^a + \Psi_n) \\
&= \mathbb{P}'(U^a + V_n, \Psi^a + \Psi_n)(\delta V_n, \delta \Psi_n) + e'_n \\
&= \mathbb{P}'(U^a + S_{\theta_n} V_n, \Psi^a + S_{\theta_n} \Psi_n)(\delta V_n, \delta \Psi_n) + e'_n + e''_n,
\end{aligned}
$$

where $e'_n$ denotes the usual "quadratic" error of Newton's scheme, and $e''_n$ the "first substitution" error. The operator $\mathbb{P}'$ is given explicitly in (36).

Similarly, in the vacuum side we have

$$
\mathcal{E}(\mathcal{K}_n, \Psi_{n+1}) - \mathcal{E}(\mathcal{K}_n, \Psi_n) = \mathbb{V}'(\mathcal{H}^a + S_{\theta_n} \mathcal{K}_n, \Psi^a + S_{\theta_n} \Psi_n)(\delta \mathcal{K}_n, \delta \Psi_n) + \hat{e}'_n + \hat{e}''_n,
$$

where $\hat{e}'_n$ denotes the "quadratic" error, and $\hat{e}''_n$ the "first substitution" error. The operator $\mathbb{V}'$ is given in (36) as well.

On the boundary $\Gamma$ we have

$$
\begin{aligned}
\mathcal{B}((V_{n+1}, \mathcal{K}_{n+1})_{|_{x_1=0}}, \psi_{n+1}) &- \mathcal{B}((V_n, \mathcal{K}_n)_{|_{x_1=0}}, \psi_n) \\
&= \mathbb{B}'\big((U^a + V_n, \mathcal{H}^a + \mathcal{K}_n)_{|_{x_1=0}}, \varphi^a + \psi_n\big)((\delta V_n, \delta \mathcal{K}_n)_{|_{x_1=0}}, \delta \psi_n) + \tilde{e}'_n \\
&= \mathbb{B}'\big((U^a + S_{\theta_n} V_n, \mathcal{H}^a + S_{\theta_n} \mathcal{K}_n)_{|_{x_1=0}}, \varphi^a + S_{\theta_n} \psi_n\big)((\delta V_n, \delta \mathcal{K}_n)_{|_{x_1=0}}, \delta \psi_n) + \tilde{e}'_n + \tilde{e}''_n,
\end{aligned}
$$

where $\tilde{e}'_n$ denotes the "quadratic" error, and $\tilde{e}''_n$ the "first substitution" error. This decomposition is meaningful only for the first three components, defined on $\omega_T$, as the last two components of $\mathbb{B}$ are linear.

The inversion of the operator $(\mathbb{L}', \mathbb{V}', \mathbb{B}')$ requires the linearization around a state satisfying the constraints (29)–(34), (61), that is the constraints of the basic state in Section 4. We thus need to introduce a smooth modified state, denoted $V_{n+1/2}, \mathcal{K}^{n+1/2}, \Psi_{n+1/2}, \psi_{n+1/2}$, that satisfies the above mentioned constraints. (The exact definition of this intermediate state is detailed in subsection 12.4.) A similar difficulty was found in [8, 36]. Accordingly, we introduce the decompositions

$$
\mathcal{L}(V_{n+1}, \Psi_{n+1}) - \mathcal{L}(V_n, \Psi_n) = \mathbb{P}'(U^a + V_{n+1/2}, \Psi^a + \Psi_{n+1/2})(\delta V_n, \delta \Psi_n) + e'_n + e''_n + e'''_n,
$$

$$
\mathcal{E}(\mathcal{K}_{n+1}, \Psi_{n+1}) - \mathcal{E}(\mathcal{K}_n, \Psi_n) = \mathbb{V}'(\mathcal{H}^a + \mathcal{K}^{n+1/2}, \Psi^a + \Psi_{n+1/2})(\delta \mathcal{K}_n, \delta \Psi_n) + \hat{e}'_n + \hat{e}''_n + \hat{e}'''_n,
$$

$$
\begin{aligned}
\mathcal{B}((V_{n+1}, \mathcal{K}_{n+1})_{|_{x_1=0}}, \psi_{n+1}) - \mathcal{B}((V_n, \mathcal{K}_n)_{|_{x_1=0}}, \psi_n) = \\
\mathbb{B}'\big((U^a + V_{n+1/2}, \mathcal{H}^a + \mathcal{K}^{n+1/2})_{|_{x_1=0}}, \varphi^a + \psi_{n+1/2}\big)((\delta V_n, \delta \mathcal{K}_n)_{|_{x_1=0}}, \delta \psi_n) + \tilde{e}'_n + \tilde{e}''_n + \tilde{e}'''_n,
\end{aligned}
$$

where $e'''_n, \hat{e}'''_n, \tilde{e}'''_n$ denote the "second substitution" errors.

The final step is the introduction of the "good unknown" (compare with (35)):

$$
\delta \dot{V}_n := \delta V_n - \delta \Psi_n \frac{\partial_1 (U^a + V_{n+1/2})}{\partial_1 (\Phi_1^a + \Psi_{n+1/2})}, \qquad \delta \dot{\mathcal{K}}_n := \delta \mathcal{K}_n - \delta \Psi_n \frac{\partial_1 (\mathcal{H}^a + \mathcal{K}^{n+1/2})}{\partial_1 (\Phi_1^a + \Psi_{n+1/2})}. \tag{158}
$$

For the interior equations this leads to

$$
\begin{aligned}
\mathcal{L}(V_{n+1}, \Psi_{n+1}) - \mathcal{L}(V_n, \Psi_n) &= \mathbb{P}'_e(U^a + V_{n+1/2}, \Psi^a + \Psi_{n+1/2}) \delta \dot{V}_n \\
&\quad + e'_n + e''_n + e'''_n + \frac{\delta \Psi_n}{\partial_1 (\Phi_1^a + \Psi_{n+1/2})} \partial_1 \big\{ \mathbb{P}(U^a + V_{n+1/2}, \Psi^a + \Psi_{n+1/2}) \big\}, \tag{159}
\end{aligned}
$$

$$
\begin{aligned}
\mathcal{E}(\mathcal{K}_{n+1}, \Psi_{n+1}) - \mathcal{E}(\mathcal{K}_n, \Psi_n) &= \mathbb{V}'_e(\mathcal{H}^a + \mathcal{K}^{n+1/2}, \Psi^a + \Psi_{n+1/2}) \delta \dot{\mathcal{K}}_n \\
&\quad + \hat{e}'_n + \hat{e}''_n + \hat{e}'''_n + \frac{\delta \Psi_n}{\partial_1 (\Phi_1^a + \Psi_{n+1/2})} \partial_1 \big\{ \mathbb{V}(\mathcal{H}^a + \mathcal{K}^{n+1/2}, \Psi^a + \Psi_{n+1/2}) \big\}, \tag{160}
\end{aligned}
$$

recalling (36), (38). For the boundary terms we obtain

$$
\begin{aligned}
\mathcal{B}((V_{n+1}, \mathcal{K}_{n+1})_{|_{x_1=0}}, \psi_{n+1}) &- \mathcal{B}((V_n, \mathcal{K}_n)_{|_{x_1=0}}, \psi_n) \\
&= \mathbb{B}'_e((U^a + V_{n+1/2}, \mathcal{H}^a + \mathcal{K}^{n+1/2})_{|_{x_1=0}}, \varphi^a + \psi_{n+1/2})((\delta \dot{V}_n, \delta \dot{\mathcal{K}}_n)_{|_{x_1=0}}, \delta \psi_n) + \tilde{e}'_n + \tilde{e}''_n + \tilde{e}'''_n, \tag{161}
\end{aligned}
$$

where $\mathbb{B}'_e$ is defined in (37). For the sake of brevity we set

$$D_{n+1/2} := \frac{1}{\partial_1(\Phi_1^a + \Psi_{n+1/2})} \, \partial_1 \Big\{ \mathbb{P}(U^a + V_{n+1/2}, \Psi^a + \Psi_{n+1/2}) \Big\},$$

$$\hat{D}_{n+1/2} := \frac{1}{\partial_1(\Phi_1^a + \Psi_{n+1/2})} \, \partial_1 \Big\{ \mathbb{V}(\mathcal{H}^a + \mathcal{K}^{n+1/2}, \Psi^a + \Psi_{n+1/2}) \Big\},$$

$$\mathbb{B}'_{n+1/2} := \mathbb{B}'_e\big(U^a + V_{n+1/2}, \mathcal{H}^a + \mathcal{K}^{n+1/2}, \varphi^a + \psi_{n+1/2}\big).$$

Let us also set

$$e_n := e'_n + e''_n + e'''_n + D_{n+1/2} \, \delta\Psi_n \, , \quad \hat{e}_n := \hat{e}'_n + \hat{e}''_n + \hat{e}'''_n + \hat{D}_{n+1/2} \, \delta\Psi_n \, , \quad \tilde{e}_n := \tilde{e}'_n + \tilde{e}''_n + \tilde{e}'''_n. \quad (162)$$

The iteration proceeds as follows. Given

$$V_0 := 0 \, , \quad \mathcal{K}_0 := 0 \, , \quad \Psi_0 := 0 \, , \quad \psi_0 := 0 \, ,$$
$$f_0 := S_{\theta_0} f^a \, , \quad \hat{f}_0 := 0 \, , \quad \tilde{f}_0 := 0 \, , \quad E_0 := 0 \, , \quad \hat{E}_0 := 0 \, , \quad \tilde{E}_0 := 0 \, ,$$
$$V_1, \ldots, V_n \, , \quad \mathcal{K}_1, \ldots, \mathcal{K}_n \, , \quad \Psi_1, \ldots, \Psi_n \, , \quad \psi_1, \ldots, \psi_n \, ,$$
$$f_1, \ldots, f_{n-1} \, , \quad \hat{f}_1, \ldots, \hat{f}_{n-1} \, , \quad \tilde{f}_1, \ldots, \tilde{f}_{n-1} \, , \quad e_0, \ldots, e_{n-1} \, , \quad \hat{e}_0, \ldots, \hat{e}_{n-1} \, , \quad \tilde{e}_0, \ldots, \tilde{e}_{n-1} \, ,$$

we first compute for $n \geq 1$

$$E_n := \sum_{k=0}^{n-1} e_k \, , \qquad \hat{E}_n := \sum_{k=0}^{n-1} \hat{e}_k \, , \qquad \tilde{E}_n := \sum_{k=0}^{n-1} \tilde{e}_k \, . \quad (163)$$

These are the accumulated errors at the step $n$. Then we compute $f_n, \hat{f}_n$, and $\tilde{f}_n$ from the equations:

$$\sum_{k=0}^{n} f_k + S_{\theta_n} E_n = S_{\theta_n} f^a \, , \qquad \sum_{k=0}^{n} \hat{f}_k + S_{\theta_n} \hat{E}_n = 0 \, , \qquad \sum_{k=0}^{n} \tilde{f}_k + S_{\theta_n} \tilde{E}_n = 0 \, , \quad (164)$$

and we solve the linear problem

$$\begin{cases} \mathbb{P}'_e(U^a + V_{n+1/2}, \Psi^a + \Psi_{n+1/2}) \, \delta\dot{V}_n = f_n & \text{in } Q_T^+ \, , \\ \mathbb{V}'_e(\mathcal{H}^a + \mathcal{K}^{n+1/2}, \Psi^a + \Psi_{n+1/2}) \, \delta\dot{\mathcal{K}}_n = \hat{f}_n & \text{in } Q_T^- \, , \\ \mathbb{B}'_{n+1/2}(\delta\dot{V}_n, \delta\dot{\mathcal{K}}_n, \delta\psi_n) = \tilde{f}_n & \text{on } \omega_T^3 \times \omega_T^\pm \, , \\ \delta\dot{V}_n = 0, \quad \delta\dot{\mathcal{K}}_n = 0, \quad \delta\psi_n = 0 & \text{for } t < 0 \, , \end{cases} \quad (165)$$

finding $(\delta\dot{V}_n, \delta\dot{\mathcal{K}}_n, \delta\psi_n)$. Then we construct $\delta\Psi_n$ from $\delta\psi_n$ by Lemma 2, the functions $\delta V_n, \delta\mathcal{K}_n$ are obtained from (158), and the functions $V_{n+1}, \mathcal{K}_{n+1}, \Psi_{n+1}, \psi_{n+1}$ are obtained from (157). Finally, we compute $e_n, \hat{e}_n, \tilde{e}_n$ from

$$\mathcal{L}(V_{n+1}, \Psi_{n+1}) - \mathcal{L}(V_n, \Psi_n) = f_n + e_n \, ,$$

$$\mathcal{E}(V_{n+1}, \Psi_{n+1}) - \mathcal{E}(V_n, \Psi_n) = \hat{f}_n + \hat{e}_n \, , \quad (166)$$

$$\mathcal{B}(V_{n+1}, \mathcal{K}_{n+1}, \psi_{n+1}) - \mathcal{B}(V_n, \mathcal{K}_n, \psi_n) = \tilde{f}_n + \tilde{e}_n \, .$$

Adding (166) from 0 to $N$, and combining with (164) gives

$$\mathcal{L}(V_{N+1}, \Psi_{N+1}) - f^a = (S_{\theta_N} - I) f^a + (I - S_{\theta_N}) E_N + e_N \, ,$$

$$\mathcal{E}(\mathcal{K}_{N+1}, \Psi_{N+1}) = (I - S_{\theta_N}) \hat{E}_N + \tilde{e}_N \, ,$$

$$\mathcal{B}(V_{N+1}, \mathcal{K}_{N+1}, \psi_{N+1}) = (I - S_{\theta_N}) \tilde{E}_N + \tilde{e}_N \, .$$

Because $S_{\theta_N} \to I$ as $N \to +\infty$, and since we expect $(e_N, \hat{e}_N, \tilde{e}_N) \to 0$, we will formally obtain the solution of the problem (153) from

$$\mathcal{L}(V_{N+1}, \Psi_{N+1}) \to f^a \, , \qquad \mathcal{E}(\mathcal{K}_{N+1}, \Psi_{N+1}) \to 0 \, , \qquad \mathcal{B}(V_{N+1}, \mathcal{K}_{N+1}, \psi_{N+1}) \to 0 \, .$$

11.3. **Tame estimate for the second derivatives.** For the control of the errors in the iteration scheme, we need to estimate the second derivative of the operators $\mathbb{P}$, $\mathbb{V}$, and $\mathbb{B}$. Let us first define the spaces

$$W_*^{1,\infty}(Q_T^+) = \{u \in L^\infty(Q_T^+) : Z_i u \in L^\infty(Q_T^+), i = 0, \ldots, 3\},$$
$$W_*^{2,\infty}(Q_T^+) = \{u \in W_*^{1,\infty}(Q_T^+) : \nabla u \in W_*^{1,\infty}(Q_T^+)\}, \tag{167}$$

equipped with its natural norms. From Theorem 41 we have $H_*^6(Q_T^+) \hookrightarrow W_*^{2,\infty}(Q_T^+)$. We consider a fixed time $T > 0$, and we take $(\hat{U}, \hat{\mathcal{H}}, \hat{\varphi})$ such that

$$\hat{U} \in W_*^{2,\infty}(Q_T^+), \quad \hat{\mathcal{H}} \in W^{1,\infty}(Q_T^-), \quad \hat{\Psi} \in W^{2,\infty}(Q_T),$$
$$\|\hat{U}\|_{W_*^{2,\infty}(Q_T^+)} + \|\hat{\mathcal{H}}\|_{W^{1,\infty}(Q_T^-)} + \|\hat{\Psi}\|_{W^{2,\infty}(Q_T)} \leq \tilde{K}, \qquad \|\hat{\varphi}\|_{\mathcal{C}([0,T];H^{2.5}(\Gamma))} \leq \epsilon_0, \tag{168}$$

where $\tilde{K}$ is a positive constant. As usual, corresponding to the given $\hat{\varphi}$ we construct $\hat{\Psi}$ and the diffeomorphism $\hat{\Phi}$ as in Lemmata 2 and 3. Let $\tilde{\alpha}$ be a sufficiently large integer that will be chosen later on. We have the following result:

**Proposition 23.** *Let* $m \in \mathbb{N}, m \in [6, \tilde{\alpha} - 2]$, *and let* $T > 0$. *Assume that* $(\hat{U}, \hat{\mathcal{H}}, \hat{\varphi})$ *satisfy* (168), *and*

$$(\hat{U}, \hat{\mathcal{H}}, \hat{\Psi}) \in H_{*,\gamma}^{\tilde{\alpha}}(Q_T^+) \times H_\gamma^{\tilde{\alpha}}(Q_T^-) \times H_\gamma^{\tilde{\alpha}}(Q_T).$$

*Then there exist two constants* $\tilde{K}_0 > 0$, *and* $C > 0$, *dependent on* $\tilde{K}_0$ *but independent of* $T$, *such that, if* $\tilde{K} \leq \tilde{K}_0$, *and if* $(V', \Psi'), (V'', \Psi'') \in H_{*,\gamma}^{m+2}(Q_T^+)$, *then one has*

$$\|\mathbb{P}''(\hat{U}, \hat{\Psi})\big((V', \Psi'), (V'', \Psi'')\big)\|_{H_{*,\gamma}^m(Q_T^+)}$$
$$\leq C \Big\{ \|(\hat{U}, \hat{\Psi})\|_{H_{*,\gamma}^{m+2}(Q_T^+)} \|(V', \Psi')\|_{W_*^{2,\infty}(Q_T^+)} \|(V'', \Psi'')\|_{W_*^{2,\infty}(Q_T^+)}$$
$$+ \|(V', \Psi')\|_{H_{*,\gamma}^{m+2}(Q_T^+)} \|(V'', \Psi'')\|_{W_*^{2,\infty}(Q_T^+)} + \|(V'', \Psi'')\|_{H_{*,\gamma}^{m+2}(Q_T^+)} \|(V', \Psi')\|_{W_*^{2,\infty}(Q_T^+)} \Big\}. \tag{169}$$

*If* $(\mathcal{K}', \nabla\Psi'), (\mathcal{K}'', \nabla\Psi'') \in H_\gamma^{m+1}(Q_T^-) \times H_\gamma^{m+1}(Q_T^-)$ *one has*

$$\|\mathbb{V}''(\hat{\mathcal{H}}, \hat{\Psi})\big((\mathcal{K}', \Psi'), (\mathcal{K}'', \Psi'')\big)\|_{H_\gamma^m(Q_T^-)}$$
$$\leq C \Big\{ \|(\hat{\mathcal{H}}, \nabla\hat{\Psi})\|_{H_\gamma^{m+1}(Q_T^-)} \|(\mathcal{K}', \nabla\Psi')\|_{W^{1,\infty}(Q_T^-)} \|(\mathcal{K}'', \nabla\Psi'')\|_{W^{1,\infty}(Q_T^-)}$$
$$+ \|(\mathcal{K}', \nabla\Psi')\|_{H_\gamma^{m+1}(Q_T^-)} \|(\mathcal{K}'', \nabla\Psi'')\|_{W^{1,\infty}(Q_T^-)}$$
$$+ \|(\mathcal{K}'', \nabla\Psi'')\|_{H_\gamma^{m+1}(Q_T^-)} \|(\mathcal{K}', \nabla\Psi')\|_{W^{1,\infty}(Q_T^-)} \Big\}. \tag{170}$$

*If* $(V', \mathcal{K}', \psi'), (V'', \mathcal{K}'', \psi'') \in H_{*,\gamma}^{m+1}(Q_T^+) \times H_\gamma^{m+1}(Q_T^-) \times H_\gamma^{m+1}(\omega_T)$, *then one has*

$$\|\mathbb{B}''\big((V', \mathcal{K}', \psi'), (V'', \mathcal{K}'', \psi'')\big)_{1,2,3}\|_{H_\gamma^m(\omega_T)} \leq C \Big\{ \|V'\|_{H_{*,\gamma}^{m+1}(Q_T^+)} \|\psi''\|_{W^{1,\infty}(\omega_T)}$$
$$+ \|V'\|_{L^\infty(Q_T^+)} \|\psi''\|_{H_\gamma^{m+1}(\omega_T)} + \|V''\|_{H_{*,\gamma}^{m+1}(Q_T^+)} \|\psi'\|_{W^{1,\infty}(\omega_T)} + \|V''\|_{L^\infty(Q_T^+)} \|\psi'\|_{H_\gamma^{m+1}(\omega_T)}$$
$$+ \|\mathcal{K}'\|_{H_\gamma^3(Q_T^-)} \|\mathcal{K}''\|_{H_\gamma^{m+1}(Q_T^-)} + \|\mathcal{K}'\|_{H_\gamma^{m+1}(Q_T^-)} \|\mathcal{K}''\|_{H_\gamma^3(Q_T^-)}$$
$$+ \|\mathcal{K}'\|_{H_\gamma^3(Q_T^-)} \|\psi''\|_{H_\gamma^{m+1}(\omega_T)} + \|\mathcal{K}'\|_{H_\gamma^{m+1}(Q_T^-)} \|\psi''\|_{W^{1,\infty}(\omega_T)}$$
$$+ \|\mathcal{K}''\|_{H_\gamma^3(Q_T^-)} \|\psi'\|_{H_\gamma^{m+1}(\omega_T)} + \|\mathcal{K}''\|_{H_\gamma^{m+1}(Q_T^-)} \|\psi'\|_{W^{1,\infty}(\omega_T)} \Big\}. \tag{171}$$

*Proof.* The proof follows from the long, but straightforward calculation of the explicit expression of $\mathbb{P}''$, $\mathbb{V}''$, $\mathbb{B}''$, from Moser-type inequalities in standard Sobolev spaces, see Lemma 49, and from Theorem 40 when we argue in $H_*^m$ spaces. Also for later use, it is useful to observe that $H_*^3(Q_T^+) \hookrightarrow L^\infty(Q_T^+)$, $H_*^4(Q_T^+) \hookrightarrow W_*^{1,\infty}(Q_T^+)$, see Theorem 41. (171) regards only the first three components, defined on $\omega_T$, as the last two components of $\mathbb{B}$ are linear and therefore the second order derivative is zero. For its proof we use the trace estimate $\|u_{|\omega_T}\|_{H_\gamma^m(\omega_T)} \leq C\|u\|_{H_{*,\gamma}^{m+1}(Q_T^+)}$, see [25]. The constant $\tilde{K}_0$ is fixed so that under the constraint $\tilde{K} \leq \tilde{K}_0$, $U$ takes its values in a fixed compact domain of the hyperbolicity region. $\square$

The estimates (169), (170), (171) hold for every $m$, with a constant $C$ that may depend on $m$. Since in Proposition 23, $m$ is taken in a bounded interval, the constant $C$ may be assumed to be independent of $m$.

Without loss of generality, we assume that the constant $\tilde{K}_0 = 2C(M_0)$, where $C(M_0)$ is the constant in (143).

## 12. Proof of the existence of smooth solutions

We recall that the sequence $(\theta_n)$ is defined by $\theta_0 \geq 1$, $\theta_n := (\theta_0^2 + n)^{1/2}$, and that we denote $\Delta_n := \theta_{n+1} - \theta_n$. In particular, the sequence $(\Delta_n)$ is decreasing, and tends to zero. Moreover, one has

$$\forall n \in \mathbb{N}, \quad \frac{1}{3\theta_n} \leq \Delta_n = \sqrt{\theta_n^2 + 1} - \theta_n \leq \frac{1}{2\theta_n}.$$

12.1. **Introduction of the iterative scheme.** Given an integer $\tilde{\alpha}$ that will be chosen later on, let us assume that the following estimate holds:

$$\|U^a\|_{H_*^{\tilde{\alpha}+2}(Q_T^+)} + \|\mathcal{H}^a\|_{H_\gamma^{\tilde{\alpha}+2}(Q_T^-)} + \|\Psi^a\|_{H_\gamma^{\tilde{\alpha}+3}(Q_T)} + \|\varphi^a\|_{H_\gamma^{\tilde{\alpha}+5/2}(\omega_T)} + \|f^a\|_{H_*^{\tilde{\alpha}+1}(Q_T^+)} \leq \delta'(T), \quad (172)$$

where $\delta'(T) \to 0$ as $T \to 0$. Given the integer $\alpha := m + 1$ and a small number $\delta > 0$, our inductive assumption reads:

$$(H_{n-1}) \begin{cases} \text{a) } \forall k = 0, \ldots, n-1, \quad \forall s \in [6, \tilde{\alpha}] \cap \mathbb{N}, \\ \|\delta V_k\|_{H_{*,\gamma}^s(Q_T^+)} + \|\delta \mathcal{K}_k\|_{H_\gamma^s(Q_T^-)} + \|\delta \Psi_k\|_{H_\gamma^{s+1}(Q_T)} + \|\delta \psi_k\|_{H_\gamma^{s+1/2}(\omega_T)} \leq \delta \, \theta_k^{s-\alpha-1} \, \Delta_k \,, \\[2mm] \text{b) } \forall k = 0, \ldots, n-1, \quad \forall s \in [6, \tilde{\alpha}-2] \cap \mathbb{N}, \\ \|\mathcal{L}(V_k, \Psi_k) - f^a\|_{H_{*,\gamma}^s(Q_T^+)} \leq 2 \, \delta \, \theta_k^{s-\alpha-1} \,, \qquad \|\mathcal{E}(\mathcal{K}_k, \Psi_k)\|_{H_\gamma^s(Q_T^-)} \leq 2 \, \delta \, \theta_k^{s-\alpha-1} \,, \\[2mm] \text{c) } \forall k = 0, \ldots, n-1, \quad \forall s \in [6, \tilde{\alpha}-2] \cap \mathbb{N}, \\ \|\mathcal{B}(V_k, \mathcal{K}_k, \psi_k)_{1,2,3}\|_{H_\gamma^s(\omega_T)} \leq \delta \, \theta_k^{s-\alpha-1} \,. \end{cases}$$

For $k = 0, \ldots, n$, the functions $V_k, \mathcal{K}_k, \Psi_k, \psi_k$ are also assumed to satisfy (156).

The first task is to prove that for a suitable choice of the parameters $\theta_0 \geq 1$, and $\delta > 0$, and for $T > 0$ small enough, $(H_{n-1})$ implies $(H_n)$. In the end, we shall prove that $(H_0)$ holds for $T > 0$ sufficiently small.

From now on, we assume that $(H_{n-1})$ holds. Let us show some basic consequences:

**Lemma 24.** *If $\theta_0$ is big enough, then for every $k = 0, \ldots, n$, and for every integer $s \in [6, \tilde{\alpha}]$, we have*

$$\|V_k\|_{H_{*,\gamma}^s(Q_T^+)} + \|\mathcal{K}_k\|_{H_\gamma^s(Q_T^-)} + \|\Psi_k\|_{H_\gamma^{s+1}(Q_T)} + \|\psi_k\|_{H_\gamma^{s+1/2}(\omega_T)} \leq \delta \, \theta_k^{(s-\alpha)_+}, \qquad \alpha \neq s, \quad (173a)$$

$$\|V_k\|_{H_{*,\gamma}^\alpha(Q_T^+)} + \|\mathcal{K}_k\|_{H_\gamma^\alpha(Q_T^-)} + \|\Psi_k\|_{H_\gamma^{\alpha+1}(Q_T)} + \|\psi_k\|_{H_\gamma^{\alpha+1/2}(\omega_T)} \leq \delta \log \theta_k \,. \quad (173b)$$

The proof follows from the triangle inequality, and from the classical comparisons between series and integrals.

**Lemma 25.** *If $\theta_0$ is big enough, then for every $k = 0, \ldots, n$, and for every integer $s \in [6, \tilde{\alpha}+4]$, we have*

$$\|S_{\theta_k} V_k\|_{H_{*,\gamma}^s(Q_T^+)} + \|S_{\theta_k} \mathcal{K}_k\|_{H_\gamma^s(Q_T^-)} + \|S_{\theta_k} \Psi_k\|_{H_\gamma^{s+1}(Q_T)} + \|S_{\theta_k} \psi_k\|_{H_\gamma^{s+1/2}(\omega_T)} \leq C \, \delta \, \theta_k^{(s-\alpha)_+}, \; s \neq \alpha, \quad (174a)$$

$$\|S_{\theta_k} V_k\|_{H_{*,\gamma}^\alpha(Q_T^+)} + \|S_{\theta_k} \mathcal{K}_k\|_{H_\gamma^\alpha(Q_T^-)} + \|S_{\theta_k} \Psi_k\|_{H_\gamma^{\alpha+1}(Q_T)} + \|S_{\theta_k} \psi_k\|_{H_\gamma^{\alpha+1/2}(\omega_T)} \leq C \, \delta \, \log \theta_k \,. \quad (174b)$$

*For every $k = 0, \ldots, n$, and for every integer $s \in [6, \tilde{\alpha}]$, we have*

$$\|(I - S_{\theta_k})V_k\|_{H^s_{*,\gamma}(Q_T^+)} + \|(I - S_{\theta_k})\mathcal{K}_k\|_{H^s_\gamma(Q_T^-)}$$
$$+ \|(I - S_{\theta_k})\Psi_k\|_{H^{s+1}_\gamma(Q_T)} + \|(I - S_{\theta_k})\psi_k\|_{H^{s+1/2}_\gamma(\omega_T)} \leq C\,\delta\,\theta_k^{s-\alpha}. \quad (175)$$

The proof follows from Lemma 24 and the properties (155) of the smoothing operators. The estimates (174), (175) actually hold for every $s$, with a constant $C$ that may depend on $s$. Taking $s$ in a bounded interval, the constant $C$ may be assumed to be independent of $s$.

12.2. **Estimate of the quadratic errors.** We start by proving an estimate for the quadratic errors $e'_k$, $\hat{e}'_k$, $\tilde{e}'_k$ of the iterative scheme. Recall that these errors are defined by[7]

$$e'_k := \mathcal{L}(V_{k+1}, \Psi_{k+1}) - \mathcal{L}(V_k, \Psi_k) - \mathcal{L}'(V_k, \Psi_k)(\delta V_k, \delta \Psi_k), \quad (176)$$

$$\hat{e}'_k := \mathcal{E}(\mathcal{K}_{k+1}, \Psi_{k+1}) - \mathcal{E}(\mathcal{K}_k, \Psi_k) - \mathcal{E}'(\mathcal{K}_k, \Psi_k)(\delta \mathcal{K}_k, \delta \Psi_k), \quad (177)$$

$$\tilde{e}'_k := \mathcal{B}\big((V_{k+1}, \mathcal{K}_{k+1})_{|_{x_1=0}}, \psi_{k+1}\big)_{1,2,3} - \mathcal{B}\big((V_k, \mathcal{K}_k)_{|_{x_1=0}}, \psi_k\big)_{1,2,3}$$
$$- \mathcal{B}'\big((V_k, \mathcal{K}_k)_{|_{x_1=0}}, \psi_k\big)_{1,2,3}\big((\delta V_k, \delta \mathcal{K}_k)_{|_{x_1=0}}, \delta \psi_k\big), \quad (178)$$

where $\mathcal{L}$, $\mathcal{E}$, and $\mathcal{B}$ are defined by (154).

**Lemma 26.** *Let $\alpha \geq 7$. There exist $\delta > 0$ sufficiently small, and $\theta_0 \geq 1$ sufficiently large, such that for all $k = 0, \ldots, n-1$, and for all integer $s \in [6, \tilde{\alpha} - 2]$, one has*

$$\|e'_k\|_{H^s_{*,\gamma}(Q_T^+)} \leq C\,\delta^2\,\theta_k^{L_1(s)-1}\,\Delta_k, \quad (179a)$$

$$\|\hat{e}'_k\|_{H^s_\gamma(Q_T^-)} \leq C\,\delta^2\,\theta_k^{L_1(s)-1}\,\Delta_k, \quad (179b)$$

$$\|\tilde{e}'_k\|_{H^s_\gamma(\omega_T)} \leq C\,\delta^2\,\theta_k^{s+5-2\alpha}\,\Delta_k, \quad (179c)$$

*where $L_1(s) := \max\{(s+2-\alpha)_+ + 10 - 2\alpha; s + 6 - 2\alpha\}$.*

*Proof.* The quadratic error given in (176) may be written as

$$e'_k = \int_0^1 (1-\tau)\,\mathbb{P}''(U^a + V_k + \tau\,\delta V_k, \Psi^a + \Psi_k + \tau\,\delta \Psi_k)\Big((\delta V_k, \delta \Psi_k),(\delta V_k, \delta \Psi_k)\Big)\,d\tau.$$

From Theorem 41, (172), and (173a), we have

$$\sup_{\tau \in [0,1]} \Big(\|U^a + V_k + \tau\,\delta V_k\|_{W^{2,\infty}(Q_T^+)} + \|\Psi^a + \Psi_k + \tau\,\delta \Psi_k\|_{W^{2,\infty}(Q_T)}\Big)$$
$$\leq C \sup_{\tau \in [0,1]} \Big(\|U^a + V_k + \tau\,\delta V_k\|_{H^6_{*,\gamma}(Q_T^+)} + \|\Psi^a + \Psi_k + \tau\,\delta \Psi_k\|_{H^6_{*,\gamma}(Q_T)}\Big) \leq C(M_0) + C\,\delta,$$

and so for $\delta$ sufficiently small, we can apply Proposition 23. Using $(H_{n-1})$, (143) and (173) we obtain (179a). The estimate (179b) of $\hat{e}'_k$ is similar, and follows from (170). The quadratic error $\tilde{e}'_k$ is estimated by means of (171), a classical trace estimate and the Sobolev imbedding Theorem. □

12.3. **Estimate of the first substitution errors.** Now we estimate the first substitution errors $e''_k, \hat{e}''_k, \tilde{e}''_k$ of the iterative scheme, defined by

$$e''_k := \mathcal{L}'(V_k, \Psi_k)(\delta V_k, \delta \Psi_k) - \mathcal{L}'(S_{\theta_k}V_k, S_{\theta_k}\Psi_k)(\delta V_k, \delta \Psi_k), \quad (180)$$

$$\hat{e}''_k := \mathcal{E}'(\mathcal{K}_k, \Psi_k)(\delta \mathcal{K}_k, \delta \Psi_k) - \mathcal{E}'(S_{\theta_k}\mathcal{K}_k, S_{\theta_k}\Psi_k)(\delta \mathcal{K}_k, \delta \Psi_k), \quad (181)$$

$$\tilde{e}''_k := \mathcal{B}'\big((V_k, \mathcal{K}_k)_{|_{x_1=0}}, \psi_k\big)_{1,2,3}\big((\delta V_k, \delta \mathcal{K}_k)_{|_{x_1=0}}, \delta \psi_k\big)$$
$$- \mathcal{B}'\big((S_{\theta_k}V_k, S_{\theta_k}\mathcal{K}_k)_{|_{x_1=0}}, S_{\theta_k}\psi_k\big)_{1,2,3}\big((\delta V_k, \delta \mathcal{K}_k)_{|_{x_1=0}}, \delta \psi_k\big) \quad (182)$$

---

[7] With abuse of notation with respect to Section 11, we identify the boundary errors terms $\tilde{e}'_k, \tilde{e}''_k, \tilde{e}'''_k$ with the only meaningful first three components.

**Lemma 27.** *Let $\alpha \geq 7$. There exist $\delta > 0$ sufficiently small, and $\theta_0 \geq 1$ sufficiently large, such that for all $k = 0, \ldots, n-1$, and for all integer $s \in [6, \tilde{\alpha} - 2]$, one has*

$$\|e_k''\|_{H_{*,\gamma}^s(Q_T^+)} \leq C\,\delta^2\,\theta_k^{L_2(s)-1}\,\Delta_k\,, \tag{183a}$$

$$\|\hat{e}_k''\|_{H_\gamma^s(Q_T^-)} \leq C\,\delta^2\,\theta_k^{L_2(s)-1}\,\Delta_k\,, \tag{183b}$$

$$\|\tilde{e}_k''\|_{H_\gamma^s(\omega_T)} \leq C\,\delta^2\,\theta_k^{s+7-2\alpha}\,\Delta_k\,, \tag{183c}$$

*where $L_2(s) := \max\{(s+2-\alpha)_+ + 12 - 2\alpha; s+8-2\alpha\}$.*

*Proof.* The substitution error given in (180) may be written as

$$e_k'' = \int_0^1 \mathbb{P}''\Big(U^a + S_{\theta_k}V_k + \tau(I - S_{\theta_k})V_k, \Psi^a + S_{\theta_k}\Psi_k + \tau(I - S_{\theta_k})\Psi_k\Big)$$
$$\Big((\delta V_k, \delta \Psi_k), ((I - S_{\theta_k})V_k, (I - S_{\theta_k})\Psi_k)\Big)\,d\tau\,.$$

We first show that we can apply Proposition 23 for $\delta$ sufficiently small, as in the previous proof. For $s + 2 \neq \alpha$, and $s + 2 \leq \tilde{\alpha}$, the estimate (183a) follows from (143), $(H_{n-1})$, (174a), and (175). For $s + 2 = \alpha$, the proof requires the use of (174b). (183b) follows in the same way. The substitution error given in (182) is estimated by using (171), $(H_{n-1})$, and (175). $\qquad\square$

**12.4. Construction and estimate of the modified state.** The next step requires the construction of the smooth modified state $V_{n+1/2}, \mathcal{K}^{n+1/2}, \Psi_{n+1/2}, \psi_{n+1/2}$ satisfying the constraints on the basic state stated in Section 4. We will focus especially on (31)–(33), because the additional constraints (29), (30), (61) will be simply obtained by choosing $T > 0$ small enough.

**Proposition 28.** *Let $\alpha \geq 8$. If $T > 0$ is sufficiently small, there exist some functions $V_{n+1/2} = (q_{n+1/2}, v^{n+1/2}, H^{n+1/2}, S_{n+1/2})$, $\mathcal{K}^{n+1/2}$, $\Psi_{n+1/2}$, $\psi_{n+1/2}$, that vanish in the past, and such that $U^a + V_{n+1/2}$, $\mathcal{H}^a + \mathcal{K}^{n+1/2}$, $\Phi^a + \Psi_{n+1/2}$, $\varphi^a + \psi_{n+1/2}$ satisfy the constraints (29)–(33), (61) and*

$$\Psi_{n+1/2} = S_{\theta_n}\Psi_n\,, \quad \psi_{n+1/2} := S_{\theta_n}\psi_n \tag{184a}$$

$$q_{n+1/2} = S_{\theta_n}q_n\,, \quad v_i^{n+1/2} = S_{\theta_n}v_{n,i} \text{ for } i = 2,3\,, \quad S_{n+1/2} = S_{\theta_n}S_n\,. \tag{184b}$$

*Moreover, for $\delta > 0$ sufficiently small, and $\theta_0 \geq 1$ sufficiently large, these functions satisfy:*

$$\|V_{n+1/2} - S_{\theta_n}V_n\|_{H_{*,\gamma}^s(Q_T^+)} \leq C\,\delta\,\theta_n^{s+2-\alpha}\,, \quad \text{for } s \in [6, \tilde{\alpha} + 3]\,, \tag{185a}$$

$$\|\mathcal{K}^{n+1/2} - S_{\theta_n}\mathcal{K}_n\|_{H_\gamma^s(Q_T^-)} \leq C\,\delta\,\theta_n^{s-\alpha}\,, \quad \text{for } s \in [6, \tilde{\alpha} + 3]\,. \tag{185b}$$

*Proof.* Let us define $\Psi_{n+1/2}, \psi_{n+1/2}, q_{n+1/2}, S_{n+1/2}$, and the components $v_2^{n+1/2}, v_3^{n+1/2}$ as in (184). It is easily checked that all these functions vanish in the past. Then we define the error $\varepsilon^n$ and a function $\mathcal{G}$ by

$$\varepsilon^n := \mathcal{B}(V_n, \mathcal{K}_n, \psi_n)_1 = \partial_t\psi_n + (v_i^a + v_{n,i})|_{x_1=0}\partial_i\psi_n + (v_{n,i})|_{x_1=0}\partial_i\varphi^a - (v_{n,1})|_{x_1=0}\,, \tag{186}$$

$$\mathcal{G} := \partial_t\psi_{n+1/2} + (v_i^a + v_i^{n+1/2})|_{x_1=0}\partial_i\psi_{n+1/2} + (v_i^{n+1/2})|_{x_1=0}\partial_i\varphi^a - (S_{\theta_n}v_{n,1})|_{x_1=0}\,, \tag{187}$$

where summation over $i = 2,3$ is understood, and the normal component of the velocity $v_1^{n+1/2}$ by

$$v_1^{n+1/2} := S_{\theta_n}v_{n,1} + \mathcal{R}_T\mathcal{G}\,, \tag{188}$$

where $\mathcal{R}_T$ is a lifting operator $H^{s-1}(\omega_T) \mapsto H_*^s(Q_T^+), s > 1$, see [25]. It is easily checked that $v_1^{n+1/2}$ vanishes in the past. We now prove the estimate (185a) for the part regarding $v^{n+1/2}$. We have

$$\|v^{n+1/2} - S_{\theta_n}v_n\|_{H_{*,\gamma}^s(Q_T^+)} \leq C\|\mathcal{G}\|_{H_\gamma^{s-1}(\omega_T)}.$$

Using the definitions (184b), (186), (187) we get

$$\mathcal{G} = S_{\theta_n}\varepsilon^n + [\partial_t, S_{\theta_n}]\psi_n + (v_i^a + S_{\theta_n}v_{n,i})\partial_iS_{\theta_n}\psi_n - S_{\theta_n}\big((v_i^a + v_{n,i})\partial_i\psi_n\big)$$
$$+ (S_{\theta_n}v_{n,i})\,\partial_i\varphi^a - S_{\theta_n}\big(v_{n,i}\partial_i\varphi^a\big)\,. \tag{189}$$

To estimate the first term $S_{\theta_n}\varepsilon^n$ on the right-hand side we use the decomposition:

$$\varepsilon^n = \mathcal{B}(V_{n-1}, \mathcal{K}_{n-1}, \psi_{n-1})_1 + \partial_t(\delta\psi_{n-1}) + (v_i^a + v_{n-1,i})\partial_i(\delta\psi_{n-1}) + \delta v_{n-1,i}\partial_i(\varphi^a + \psi_{n-1}) - \delta v_{n-1,1},$$

then exploit point c) of $(H_{n-1})$ and the properties of the smoothing operators. We get

$$\|S_{\theta_n}\varepsilon^n\|_{H_\gamma^{s-1}(\omega_T)} \le C\,\delta\,\theta_n^{s-\alpha-1}. \tag{190}$$

For the estimate of the commutators in (189) we proceed as in [8] and obtain

$$\|v^{n+1/2} - S_{\theta_n}v_n\|_{H_{*,\gamma}^s(Q_T^+)} \le C\,\delta\,\theta_n^{s-\alpha}, \quad \text{for } s \in [6, \tilde{\alpha}+3]. \tag{191}$$

12.4.1. *The modified plasma magnetic field.* Let us see now how to define the modified magnetic field $H^{n+1/2}$, following somehow [36]. This field should be such that $H^a + H^{n+1/2}$ satisfies (31), together with $v^a + v^{n+1/2}, \Psi^a + \Psi_{n+1/2}$, that is, denoting the equations for the magnetic field in (17) by $\mathbb{P}_H(v, H, \Psi) = 0$, we require

$$\mathbb{P}_H(v^a + v^{n+1/2}, H^a + H^{n+1/2}, \Psi^a + \Psi_{n+1/2}) = 0. \tag{192}$$

Recalling (146), i.e. $\mathbb{P}_H(v^a, H^a, \Psi^a) = 0$, and the definition (154), we write (192) as

$$\mathcal{L}_H(v^{n+1/2}, H^{n+1/2}, \Psi_{n+1/2}) = 0. \tag{193}$$

We notice that (193) does not need to be supplemented with any boundary condition. In fact, the coefficient of $\partial_1 H^{n+1/2}$ in (193) is

$$v_1^a + v_1^{n+1/2} - (v_i^a + v_i^{n+1/2})\partial_i(\varphi^a + \psi_{n+1/2}) - \partial_t(\varphi^a + \psi_{n+1/2}),$$

which vanishes at the boundary because of (145) and (188). Given $v^{n+1/2}, \Psi_{n+1/2}$ as above, we define $H^{n+1/2}$ as the unique solution vanishing in the past of (193).

In order to estimate $H^{n+1/2} - S_{\theta_n}H_n$ we first observe that (192) yields

$$\mathbb{P}_H(v^a + v^{n+1/2}, H^{n+1/2} - S_{\theta_n}H_n, \Psi^a + \Psi_{n+1/2}) =$$
$$\mathbb{P}_H(v^a + v^{n+1/2}, H^a + H^{n+1/2} - S_{\theta_n}H_n, \Psi^a + \Psi_{n+1/2}) - \mathbb{P}_H(v^a + v^{n+1/2}, H^a, \Psi^a + \Psi_{n+1/2})$$
$$= -\mathbb{P}_H(v^a + v^{n+1/2}, H^a + S_{\theta_n}H_n, \Psi^a + S_{\theta_n}\Psi_n).$$

Then $H^{n+1/2} - S_{\theta_n}H_n$ solves the equation

$$\mathbb{P}_H(v^a + v^{n+1/2}, H^{n+1/2} - S_{\theta_n}H_n, \Psi^a + \Psi_{n+1/2}) = F_H^{n+1/2}, \tag{194}$$

where

$$F_H^{n+1/2} := \Delta_1 + \Delta_2 - S_{\theta_n}\mathbb{P}_H(v^a + v_n, H^a + H_n, \Psi^a + \Psi_n),$$
$$\Delta_1 := S_{\theta_n}\mathbb{P}_H(v^a + v_n, H^a + H_n, \Psi^a + \Psi_n) - \mathbb{P}_H(v^a + S_{\theta_n}v_n, H^a + S_{\theta_n}H_n, \Psi^a + S_{\theta_n}\Psi_n),$$
$$\Delta_2 := \mathbb{P}_H(v^a + S_{\theta_n}v_n, H^a + S_{\theta_n}H_n, \Psi^a + S_{\theta_n}\Psi_n) - \mathbb{P}_H(v^a + v^{n+1/2}, H^a + S_{\theta_n}H_n, \Psi^a + S_{\theta_n}\Psi_n).$$

For $T > 0$ small enough and $\alpha \ge 7$ we obtain from (172) with $T$ sufficiently small, (173)–(175), (260)

$$\|\Delta_1\|_{H_{*,\gamma}^s(Q_T^+)} \le C\,\delta\,\theta_k^{s+2-\alpha} \quad \text{for } s \in [6, \tilde{\alpha}+3], \tag{195}$$

and from (191) we also get

$$\|\Delta_2\|_{H_{*,\gamma}^s(Q_T^+)} \le C\,\delta\,\theta_k^{s+2-\alpha} \quad \text{for } s \in [6, \tilde{\alpha}+3]. \tag{196}$$

To estimate the last term of $F_H^{n+1/2}$, we write

$$S_{\theta_n}\mathbb{P}_H(v^a + v_n, H^a + H_n, \Psi^a + \Psi_n) = S_{\theta_n}\mathcal{L}_H(V_n, \Psi_n) =$$
$$= S_{\theta_n}\mathcal{L}_H(V_{n-1}, \Psi_{n-1}) + S_{\theta_n}\big(\mathcal{L}_H(V_{n-1} + \delta V_{n-1}, \Psi_{n-1} + \delta\Psi_{n-1}) - \mathcal{L}_H(V_{n-1}, \Psi_{n-1})\big).$$

From (155) and point b) of $(H_{n-1})$ we have:

$$\|S_{\theta_n}\mathcal{L}_H(V_{n-1}, \Psi_{n-1})\|_{H_{*,\gamma}^s(Q_T^+)} \le C\,\theta_n^{s-6}\,\|\mathcal{L}_H(V_{n-1}, \Psi_{n-1})\|_{H_{*,\gamma}^6(Q_T^+)} \le C\,\delta\,\theta_n^{s-\alpha-1}, \tag{197}$$

for all integer $s \geq 6$. Similarly, from (155), point a) of $(H_{n-1})$, (173) and (259) we obtain

$$\|S_{\theta_n}\big(\mathcal{L}_H(V_{n-1} + \delta V_{n-1}, \Psi_{n-1} + \delta\Psi_{n-1}) - \mathcal{L}_H(V_{n-1}, \Psi_{n-1})\big)\|_{H^s_{*,\gamma}(Q^+_T)} \leq C\,\delta\,\theta_n^{s-\alpha-1}, \tag{198}$$

for all integer $s \geq 6$. Collecting (195)-(198) yields

$$\|F_H^{n+1/2}\|_{H^s_{*,\gamma}(Q^+_T)} \leq C\,\delta\,\theta_n^{s+2-\alpha} \quad \text{for } s \in [6, \tilde{\alpha}+3]. \tag{199}$$

Now, equation (194) solved by $H^{n+1/2} - S_{\theta_n}H_n$ has the form

$$\partial_t Y + \textstyle\sum_{j=1}^3 \mathcal{D}_j(b)\partial_j Y + Q(b)Y = F_H^{n+1/2}, \tag{200}$$

for $Y = H^{n+1/2} - S_{\theta_n}H_n$, $b = (v^a + v^{n+1/2}, \Psi^a + \Psi_{n+1/2})$, and where $\mathcal{D}_j$ and $Q$ are some matrices. The matrices $\mathcal{D}_j$ are diagonal and, more important, $\mathcal{D}_1$ vanishes at $\{x_1 = 0\} \cup \{x_1 = 1\}$. This yields that system (200) does not need any boundary condition. A standard energy argument applied to (200) and (174), (191), (259), (260) give for $\gamma$ large and $T > 0$ small the a priori estimate

$$\gamma\|Y_\gamma\|_{H^s_{*,\gamma}(Q^+_T)} \leq C\left(\|(F_H^{n+1/2})_\gamma\|_{H^s_{*,\gamma}(Q^+_T)} + \delta^2\,\theta_n^{(s+1-\alpha)_+}\|Y_\gamma\|_{H^6_{*,\gamma}(Q^+_T)}\right). \tag{201}$$

Choosing $s = 6, \alpha \geq 8$, and taking $\delta$ small in (201) yields

$$\gamma\|Y_\gamma\|_{H^6_{*,\gamma}(Q^+_T)} \leq C\|(F_H^{n+1/2})_\gamma\|_{H^6_{*,\gamma}(Q^+_T)},$$

and substituting in (201) gives

$$\gamma\|Y_\gamma\|_{H^s_{*,\gamma}(Q^+_T)} \leq C\left(\|(F_H^{n+1/2})_\gamma\|_{H^s_{*,\gamma}(Q^+_T)} + \delta^2\,\theta_n^{(s+1-\alpha)_+}\|(F_H^{n+1/2})_\gamma\|_{H^6_{*,\gamma}(Q^+_T)}\right). \tag{202}$$

Finally, plugging (199) in (202) yields

$$\gamma\|Y_\gamma\|_{H^s_{*,\gamma}(Q^+_T)} = \gamma\|(H^{n+1/2} - S_{\theta_n}H_n)_\gamma\|_{H^s_{*,\gamma}(Q^+_T)} \leq C\delta\,\theta_n^{s+2-\alpha} \quad \text{for } s \in [6, \tilde{\alpha}+3]. \tag{203}$$

This completes the proof of (185a).

From now on $\gamma$ is assumed fixed, satisfying the requirements for Theorem 16 and the proof of (203).

12.4.2. *The modified vacuum magnetic field.* The modified vacuum magnetic field $\mathcal{K}^{n+1/2}$ is supposed to satisfy with $\mathcal{H}^a$ and $\Psi^a + \Psi_{n+1/2}$ the constraints (32), (33)$_{2,4}$. This means[8]

$$\begin{aligned}
&\partial_1\{\mathcal{H}_1^a + \mathcal{K}_1^{n+1/2} - (\mathcal{H}_i^a + \mathcal{K}_i^{n+1/2})\partial_i(\Psi^a + \Psi_{n+1/2})\} \\
&\qquad\qquad + \partial_i\{(\mathcal{H}_i^a + \mathcal{K}_i^{n+1/2})\partial_1(\Phi_1^a + \Psi_{n+1/2})\} = 0 \quad \text{in } Q_T^-, \\
&\mathcal{H}_1^a + \mathcal{K}_1^{n+1/2} - (\mathcal{H}_i^a + \mathcal{K}_i^{n+1/2})\partial_i(\varphi^a + \psi_{n+1/2}) = 0 \qquad \text{on } \omega_T, \\
&\nu \times (\mathcal{H}^a + \mathcal{K}^{n+1/2}) = \mathfrak{J} \qquad\qquad\qquad\qquad\qquad\qquad \text{on } \omega_T^-.
\end{aligned}$$

Taking account of (148), this can be rewritten as

$$\begin{aligned}
&\partial_1\{\mathcal{K}_1^{n+1/2} - \mathcal{K}_i^{n+1/2}\partial_i(\Psi^a + \Psi_{n+1/2})\} + \partial_i\{\mathcal{K}_i^{n+1/2}\partial_1(\Phi_1^a + \Psi_{n+1/2})\} \\
&\qquad = \partial_1\{\mathcal{H}_i^a\partial_i\Psi_{n+1/2}\} - \partial_i\{\mathcal{H}_i^a\partial_1\Psi_{n+1/2}\} \qquad\qquad \text{in } Q_T^-, \\
&\mathcal{K}_1^{n+1/2} - \mathcal{K}_i^{n+1/2}\partial_i(\varphi^a + \psi_{n+1/2}) = \mathcal{H}_i^a\partial_i\psi_{n+1/2} \qquad\qquad \text{on } \omega_T, \\
&\nu \times \mathcal{K}^{n+1/2} = 0 \qquad\qquad\qquad\qquad\qquad\qquad\qquad\qquad \text{on } \omega_T^-.
\end{aligned} \tag{204}$$

Denoting

$$\begin{aligned}
&\mathfrak{h}^{n+1/2} = (\mathcal{K}_N^{n+1/2}, \mathcal{K}_2^{n+1/2}\partial_1(\Phi_1^a + \Psi_{n+1/2}), \mathcal{K}_3^{n+1/2}\partial_1(\Phi_1^a + \Psi_{n+1/2})), \\
&\mathcal{K}_N^{n+1/2} = \mathcal{K}_1^{n+1/2} - \mathcal{K}_i^{n+1/2}\partial_i(\Psi^a + \Psi_{n+1/2}),
\end{aligned}$$

(204) is rephrased as

$$\begin{aligned}
&\text{div}\,\mathfrak{h}^{n+1/2} = \partial_1\{\mathcal{H}_i^a\partial_i\Psi_{n+1/2}\} - \partial_i\{\mathcal{H}_i^a\partial_1\Psi_{n+1/2}\} \quad \text{in } Q_T^-, \\
&\mathfrak{h}_1^{n+1/2} = \mathcal{H}_i^a\partial_i\psi_{n+1/2} \qquad\qquad\qquad\qquad\qquad \text{on } \omega_T \\
&\nu \times \mathcal{K}^{n+1/2} = 0 \qquad\qquad\qquad\qquad\qquad\qquad\quad \text{on } \omega_T^-.
\end{aligned} \tag{205}$$

---

[8]Here summation over $i = 2, 3$ is understood.

For given $\Psi_{n+1/2}$, we define $\mathcal{K}^{n+1/2}$ as a solution of (205)[9]. $\mathcal{K}^{n+1/2}$ is not unique, but it is defined up to an arbitrary $\nabla \times \mathfrak{h}^{n+1/2}$. For instance, we can take $\mathcal{K}^{n+1/2}$ such that $\mathfrak{h}^{n+1/2}$ solves (205) and $\nabla \times \mathfrak{h}^{n+1/2} = 0$ in $Q_T^-$.

In order to estimate $\mathcal{K}^{n+1/2} - S_{\theta_n}\mathcal{K}_n$ we consider the problem (compare with (204))

$$
\begin{aligned}
\partial_1\{\mathcal{K}_1^{n+1/2} - S_{\theta_n}\mathcal{K}_{n,1} - (\mathcal{K}_i^{n+1/2} - S_{\theta_n}\mathcal{K}_{n,i})\partial_i(\Psi^a + \Psi_{n+1/2})\} & \\
+ \partial_i\{(\mathcal{K}_i^{n+1/2} - S_{\theta_n}\mathcal{K}_{n,i})\partial_1(\Phi_1^a + \Psi_{n+1/2})\} = G^{n+1/2} & \quad \text{in } Q_T^-, \\
\mathcal{K}_1^{n+1/2} - S_{\theta_n}\mathcal{K}_{n,1} - (\mathcal{K}_i^{n+1/2} - S_{\theta_n}\mathcal{K}_{n,i})\partial_i(\varphi^a + \psi_{n+1/2}) = g^{n+1/2} & \quad \text{on } \omega_T, \\
\nu \times (\mathcal{K}^{n+1/2} - S_{\theta_n}\mathcal{K}_n) = -\nu \times S_{\theta_n}\mathcal{K}_n & \quad \text{on } \omega_T^-,
\end{aligned}
\tag{206}
$$

where we have set

$$
\begin{aligned}
G^{n+1/2} := \partial_1\{\mathcal{H}_i^a\partial_i\Psi_{n+1/2}\} - \partial_i\{\mathcal{H}_i^a\partial_1\Psi_{n+1/2}\} & \\
- \partial_1\{S_{\theta_n}\mathcal{K}_{n,1} - (S_{\theta_n}\mathcal{K}_{n,i})\partial_i(\Psi^a + \Psi_{n+1/2})\} - \partial_i\{(S_{\theta_n}\mathcal{K}_{n,i})\partial_1(\Phi_1^a + \Psi_{n+1/2})\}, & \\
g^{n+1/2} := \mathcal{H}_i^a\partial_i\psi_{n+1/2} - S_{\theta_n}\mathcal{K}_{n,1} + (S_{\theta_n}\mathcal{K}_{n,i})\partial_i(\varphi^a + \psi_{n+1/2}). &
\end{aligned}
$$

As for $G^{n+1/2}$, we decompose it as

$$
G^{n+1/2} = \Delta_3 + \Delta_4,
$$

where

$$
\Delta_3 := -\partial_1\{\mathcal{K}_{n,1} - (\mathcal{H}_i^a + \mathcal{K}_{n,i})\partial_i\Psi_n - \mathcal{K}_{n,i}\partial_i\Psi^a\} - \partial_i\{(\mathcal{H}_i^a + \mathcal{K}_{n,i})\partial_1\Psi_n + \mathcal{K}_{n,i}\partial_1\Phi_1^a\},
$$

$$
\begin{aligned}
\Delta_4 := \partial_1\{\mathcal{H}_i^a\partial_i(S_{\theta_n} - I)\Psi_n\} - \partial_i\{\mathcal{H}_i^a\partial_1(S_{\theta_n} - I)\Psi_n\} & \\
- \partial_1\{(S_{\theta_n} - I)\mathcal{K}_{n,1} - (S_{\theta_n}\mathcal{K}_{n,i})\partial_i(\Psi^a + S_{\theta_n}\Psi_n) + \mathcal{K}_{n,i}\partial_i(\Psi^a + \Psi_n)\} & \\
- \partial_i\{(S_{\theta_n}\mathcal{K}_{n,i})\partial_1(\Phi_1^a + S_{\theta_n}\Psi_n) - \mathcal{K}_{n,i}\partial_1(\Phi_1^a + \Psi_n)\}. &
\end{aligned}
$$

$\Delta_3$ is decomposed as

$$
\Delta_3 = \Delta_3' + \Delta_3'',
$$

where we have set

$$
\begin{aligned}
\Delta_3' := -\partial_1\{\mathcal{K}_{n-1,1} - (\mathcal{H}_i^a + \mathcal{K}_{n-1,i})\partial_i\Psi_{n-1} - \mathcal{K}_{n-1,i}\partial_i\Psi^a\} & \\
- \partial_i\{(\mathcal{H}_i^a + \mathcal{K}_{n-1,i})\partial_1\Psi_{n-1} + \mathcal{K}_{n-1,i}\partial_1\Phi_1^a\}, &
\end{aligned}
$$

$$
\begin{aligned}
\Delta_3'' := -\partial_1\{\delta\mathcal{K}_{n-1,1} - (\mathcal{H}_i^a + \mathcal{K}_{n-1,i})\partial_i(\delta\Psi_{n-1}) - \delta\mathcal{K}_{n-1,i}\partial_i(\Psi^a + \Psi_{n-1})\} & \\
- \partial_i\{(\mathcal{H}_i^a + \mathcal{K}_{n-1,i})\partial_1(\delta\Psi_{n-1}) + \delta\mathcal{K}_{n-1,i}\partial_1(\Phi_1^a + \Psi_{n-1})\}. &
\end{aligned}
$$

Notice that

$$
\Delta_3' = -\mathbb{V}(\mathcal{H}^a + \mathcal{K}_{n-1}, \Psi^a + \Psi_{n-1}) = -\mathcal{E}(\mathcal{K}_{n-1}, \Psi_{n-1}).
$$

Then from point b) of $(H_{n-1})$ we get

$$
\|\Delta_3'\|_{H_\gamma^s(Q_T^-)} \leq 2\,\delta\,\theta_{n-1}^{s-\alpha-1} \leq C\,\delta\,\theta_n^{s-\alpha-1}.
\tag{207}
$$

We also obtain from $(H_{n-1})$ and (173)

$$
\|\Delta_3''\|_{H_\gamma^s(Q_T^-)} \leq C\,\delta\,\theta_n^{s-\alpha},
$$

which gives with (207)

$$
\|\Delta_3\|_{H_\gamma^s(Q_T^-)} \leq C\,\delta\,\theta_n^{s-\alpha}.
\tag{208}
$$

Moreover, from (173)–(175) we get the estimate

$$
\|\Delta_4\|_{H_\gamma^s(Q_T^-)} \leq C\,\delta\,\theta_n^{s-\alpha+1}.
\tag{209}
$$

Thus from (208), (209) we obtain

$$
\|G^{n+1/2}\|_{H_\gamma^s(Q_T^-)} \leq C\,\delta\,\theta_n^{s-\alpha+1}.
\tag{210}
$$

---

[9]Once $\mathfrak{h}^{n+1/2}$ is found from (205), the vector $\mathcal{K}^{n+1/2}$ is immediately obtained from the defining formula for $\mathfrak{h}^{n+1/2}$.

For the boundary term $g^{n+1/2}$ we use similar decompositions. We write

$$g^{n+1/2} = \tilde{\Delta}_3 + \tilde{\Delta}_4,$$

where

$$\tilde{\Delta}_3 := -\mathcal{K}_{n,1} + (\mathcal{H}_i^a + \mathcal{K}_{n,i})\partial_i\psi_n + \mathcal{K}_{n,i}\partial_i\varphi^a,$$

$$\tilde{\Delta}_4 := (I - S_{\theta_n})\mathcal{K}_{n,1} + \mathcal{H}_i^a\partial_i(S_{\theta_n} - I)\psi_n + \mathcal{K}_{n,i}\partial_i(S_{\theta_n} - I)\psi_n + (S_{\theta_n} - I)\mathcal{K}_{n,i}\partial_i(\varphi^a + S_{\theta_n}\psi_n),$$

and $\tilde{\Delta}_3$ is decomposed as

$$\tilde{\Delta}_3 = \tilde{\Delta}_3' + \tilde{\Delta}_3'',$$

where we have set

$$\tilde{\Delta}_3' := -\mathcal{K}_{n-1,1} + (\mathcal{H}_i^a + \mathcal{K}_{n-1,i})\partial_i\psi_{n-1} + \mathcal{K}_{n-1,i}\partial_i\varphi^a,$$

$$\tilde{\Delta}_3'' := -\delta\mathcal{K}_{n-1,1} + (\mathcal{H}_i^a + \mathcal{K}_{n-1,i})\partial_i(\delta\psi_{n-1}) + \delta\mathcal{K}_{n-1,i}\partial_i(\varphi^a + \psi_{n-1}).$$

From point c) of $(H_{n-1})$ we have

$$\|\tilde{\Delta}_3'\|_{H_\gamma^{s-1/2}(\omega_T)} \leq \|\mathcal{B}(V_{n-1}, \mathcal{K}_{n-1}, \psi_{n-1})\|_{H_\gamma^s(\omega_T)} \leq \delta\,\theta_{n-1}^{s-\alpha-1} \leq C\,\delta\,\theta_n^{s-\alpha-1}. \tag{211}$$

We also obtain

$$\|\tilde{\Delta}_3''\|_{H_\gamma^{s-1/2}(\omega_T)} \leq C\,\delta\,\theta_n^{s-\alpha-1},$$

which gives with (211)

$$\|\tilde{\Delta}_3\|_{H_\gamma^{s-1/2}(\omega_T)} \leq C\,\delta\,\theta_n^{s-\alpha-1}. \tag{212}$$

Moreover, we have the estimate

$$\|\tilde{\Delta}_4\|_{H_\gamma^{s-1/2}(\omega_T)} \leq C\,\delta\,\theta_n^{s-\alpha}. \tag{213}$$

Thus from (212), (213) we obtain

$$\|g^{n+1/2}\|_{H_\gamma^{s-1/2}(\omega_T)} \leq C\,\delta\,\theta_n^{s-\alpha}. \tag{214}$$

Finally, from (206), (210), (214) we have

$$\|\mathcal{K}^{n+1/2} - S_{\theta_n}\mathcal{K}_n\|_{H_\gamma^s(Q_T^-)} \leq C\left(\|G^{n+1/2}\|_{H_\gamma^{s-1}(Q_T^-)} + \|g^{n+1/2}\|_{H_\gamma^{s-1/2}(\omega_T)} + \|\nu \times S_{\theta_n}\mathcal{K}_n\|_{H_\gamma^{s-1/2}(\omega_T^-)}\right)$$
$$\leq C\,\delta\,\theta_n^{s-\alpha}. \tag{215}$$

which completes the proof of (185b).

12.4.3. *Conclusion of the proof.* Since the approximate solutions $U^a, \mathcal{H}^a$ satisfy (30), (61) for all times with a strict inequality, and the modified states $V_{n+1/2}, \mathcal{K}^{n+1/2}$ vanish in the past, then $U^a + V_{n+1/2}, \mathcal{H}^a + \mathcal{K}^{n+1/2}$ will satisfy (30), (61) for a sufficiently short time $T > 0$. $\qquad\square$

## 12.5. Estimate of the second substitution errors.

Now we may estimate the second substitution errors $e_k''', \hat{e}_k'''$ and $\tilde{e}_k'''$ of the iterative scheme, that are defined by

$$e_k''' := \mathcal{L}'(S_{\theta_k}V_k, S_{\theta_k}\Psi_k)(\delta V_k, \delta\Psi_k) - \mathcal{L}'(V_{k+1/2}, \Psi_{k+1/2})(\delta V_k, \delta\Psi_k), \tag{216}$$

$$\hat{e}_k''' := \mathcal{E}'(S_{\theta_k}\mathcal{K}_k, S_{\theta_k}\Psi_k)(\delta\mathcal{K}_k, \delta\Psi_k) - \mathcal{E}'(\mathcal{K}_{k+1/2}, \Psi_{k+1/2})(\delta\mathcal{K}_k, \delta\Psi_k), \tag{217}$$

$$\tilde{e}_k''' := \mathcal{B}'\big((S_{\theta_k}V_k, S_{\theta_k}\mathcal{K}_k)|_{x_1=0}, \psi_k\big)_{1,2,3}((\delta V_k, \delta\mathcal{K}_k)|_{x_1=0}, \delta\psi_k)$$
$$- \mathcal{B}'\big((V_{k+1/2}, \mathcal{K}_{k+1/2})|_{x_1=0}, \psi_{k+1/2}\big)_{1,2,3}((\delta V_k, \delta\mathcal{K}_k)|_{x_1=0}, \delta\psi_k). \tag{218}$$

**Lemma 29.** *Let $\alpha \geq 8$. There exist $T > 0$ and $\delta > 0$ sufficiently small, and $\theta_0 \geq 1$ sufficiently large such that for all $k = 0, \ldots, n-1$, and for all integer $s \in [6, \tilde{\alpha} - 2]$, one has*

$$\|e_k'''\|_{H_{*,\gamma}^s(Q_T^+)} \leq C \, \delta^2 \, \theta_k^{L_3(s)-1} \, \Delta_k \,, \tag{219a}$$

$$\|\hat{e}_k'''\|_{H_\gamma^s(Q_T^-)} \leq C \, \delta^2 \, \theta_k^{L_3(s)-1} \, \Delta_k \,, \tag{219b}$$

$$\|\tilde{e}_k'''\|_{H_\gamma^s(\omega_T)} \leq C \, \delta^2 \, \theta_k^{s+6-2\alpha} \, \Delta_k \,, \tag{219c}$$

*where $L_3(s) := \max\{(s+2-\alpha)_+ + 16 - 2\alpha; s + 10 - 2\alpha\}$.*

*Proof.* Using (184a), the substitution error given in (216) may be written as

$$e_k''' = \int_0^1 \mathbb{P}''\big(U^a + V_{k+1/2} + \tau(S_{\theta_k}V_k - V_{k+1/2}), \Psi^a + S_{\theta_k}\Psi_k\big)\big((\delta V_k, \delta \Psi_k), (S_{\theta_k}V_k - V_{k+1/2}, 0)\big)\, d\tau \,.$$

By Lemma 25 and Proposition 28 we first derive the bound

$$\|\big(U^a + V_{k+1/2} + \tau(S_{\theta_k}V_k - V_{k+1/2}), \Psi^a + S_{\theta_k}\Psi_k\big)\|_{H_{*,\gamma}^{s+2}(Q_T^+)}$$
$$\leq C \, \delta \, (\theta_k^{(s+2-\alpha)_+} + \theta_k^{s+4-\alpha}) \,, \quad s \in [6, \tilde{\alpha} - 2] \,.$$

Then (219a) follows by applying the Theorems 40 and 41, Proposition 23, $(H_{n-1})$ and Proposition 28, provided that $T > 0$ and $\delta > 0$ are small enough. A similar argument applies to $\hat{e}_k'''$ and yields (219b). We write the substitution error given in (218) as

$$\tilde{e}_k''' = \mathbb{B}''\Big(\big((\delta V_k, \delta \mathcal{K}_k)|_{x_1=0}, \delta \psi_k\big), \big((V_{k+1/2} - S_{\theta_k}V_k, \mathcal{K}_{k+1/2} - S_{\theta_k}\mathcal{K}_k)|_{x_1=0}, 0\big)\Big)_{1,2,3} \,.$$

Using the exact expression of $\mathbb{B}''$ and (184) gives

$$\tilde{e}_k''' = \begin{pmatrix} 0 \\ \delta \mathcal{K}_k \cdot (S_{\theta_k}\mathcal{K}_k - \mathcal{K}_{k+1/2}) \\ (S_{\theta_k}\mathcal{K}_{k,i} - \mathcal{K}_i^{k+1/2})\partial_i(\delta \psi_k) \end{pmatrix} \,.$$

Then (219c) follows by applying $(H_{n-1})$ and Proposition 28. $\qquad \square$

12.6. **Estimate of the last error terms.** In our iterative scheme we have two last error terms to be estimated, namely

$$D_{k+1/2}\,\delta \Psi_k := \frac{\delta \Psi_k}{\partial_1(\Phi_1^a + \Psi_{k+1/2})}\,\partial_1\left\{\mathbb{P}(U^a + V_{k+1/2}, \Psi^a + \Psi_{k+1/2})\right\} \,,$$

$$\hat{D}_{k+1/2}\,\delta \Psi_k := \frac{\delta \Psi_k}{\partial_1(\Phi_1^a + \Psi_{k+1/2})}\,\partial_1\left\{\mathbb{V}(\mathcal{H}^a + \mathcal{K}_{k+1/2}, \Psi^a + \Psi_{k+1/2})\right\} \,,$$

which result from the introduction of the good unknown in the decomposition of the linearized equations, see (159), (160). Let us set

$$R_k := \partial_1\left\{\mathbb{P}(U^a + V_{k+1/2}, \Psi^a + \Psi_{k+1/2})\right\} \,.$$

Since $V_{k+1/2}$, and $\Psi_{k+1/2}$ vanish in the past, $R_k$ does not vanish in the past. However, $\delta \Psi_k$ vanishes in the past, so the error term $D_{k+1/2}\,\delta \Psi_k$ also vanishes in the past. Moreover, Theorem 40 enables us to obtain:

$$\|D_{k+1/2}\,\delta \Psi_k\|_{H_{*,\gamma}^s(Q_T^+)} \leq C \left\{\|\delta \Psi_k\|_{H_\gamma^s(Q_T^+)}\,\|R_k\|_{W_*^{1,\infty}(Q_T^+)} \right.$$
$$\left. + \|\delta \Psi_k\|_{W_*^{1,\infty}(Q_T^+)}\big(\|R_k\|_{H_{*,\gamma}^s(Q_T^+)} + \|R_k\|_{W_*^{1,\infty}(Q_T^+)}\|\dot{\Phi}^a + \Psi_{k+1/2}\|_{H_\gamma^{s+1}(Q_T^+)}\big) \right\} \,. \tag{220}$$

**Lemma 30.** *Let $\alpha \geq 8$, $\tilde{\alpha} \geq \alpha + 2$. For $\delta > 0$, $T > 0$ sufficiently small, $\theta_0 \geq 1$ sufficiently large, for all $k = 0, \ldots, n-1$, and for all integer $s \in [4, \tilde{\alpha} - 4]$, one has*

$$\|R_k\|_{H_{*,\gamma}^s(Q_T^+)} \leq C \, \delta \left(\theta_k^{s+6-\alpha} + \theta_k^{(s+4-\alpha)_+ + 6-\alpha}\right) \,. \tag{221}$$

*Proof.* We proceed as in [1, 8]. We introduce the following decomposition:

$$\mathbb{P}(U^a + V_{k+1/2}, \Psi^a + \Psi_{k+1/2}) = \mathbb{P}(U^a + V_{k+1/2}, \Psi^a + \Psi_{k+1/2}) - \mathbb{P}(U^a + V_k, \Psi^a + \Psi_k)$$
$$+ \mathcal{L}(V_k, \Psi_k) - f^a. \quad (222)$$

Then the estimate follows from the induction assumption $(H_{n-1})$, Lemma 25 and Proposition 28. $\quad\square$

We are ready to prove the following estimate:

**Lemma 31.** *Let $\alpha \geq 8$, $\tilde{\alpha} \geq \alpha + 2$. There exist $\delta > 0$, $T > 0$ sufficiently small, and $\theta_0 \geq 1$ sufficiently large such that for all $k = 0, \ldots, n-1$, and for all integer $s \in [6, \tilde{\alpha} - 4]$, one has*

$$\|D_{k+1/2} \, \delta\Psi_k\|_{H^s_{*,\gamma}(Q^+_T)} \leq C \, \delta^2 \, \theta_k^{L_4(s)-1} \, \Delta_k \,, \quad (223)$$

*where $L_4(s) := \max\{(s - \alpha)_+ + 16 - 2\alpha; s + 12 - 2\alpha\}$.*

*Proof.* We first use Lemma 30 to derive the bound $\|R_k\|_{W^{1,\infty}(Q^+_T)} \leq C \|R_k\|_{H^4_{*,\gamma}(Q^+_T)} \leq C \, \delta \, \theta_k^{10-\alpha}$. We combine this bound and (221) in (220). The terms in $\delta\Psi_k$ are estimated by the induction assumption $(H_{n-1})$, and the terms in $\Psi_{k+1/2} = S_{\theta_k}\Psi_k$ are estimated by Lemma 25. Putting all these estimates together yields (223). $\quad\square$

A similar argument gives

**Lemma 32.** *Let $\alpha \geq 8$, $\tilde{\alpha} \geq \alpha + 2$. There exist $T > 0$ sufficiently small, and $\theta_0 \geq 1$ sufficiently large such that for all $k = 0, \ldots, n-1$, and for all integer $s \in [6, \tilde{\alpha} - 2]$, one has*

$$\|\hat{D}_{k+1/2} \, \delta\Psi_k\|_{H^s_\gamma(Q^-_T)} \leq C \, \delta^2 \, \theta_k^{L_4(s)-1} \, \Delta_k \,. \quad (224)$$

**12.7. Convergence of the iteration scheme.** We first estimate the errors $e_k$, $\hat{e}_k$, and $\tilde{e}_k$:

**Lemma 33.** *Let $\alpha \geq 8$. There exist $\delta > 0$, $T > 0$ sufficiently small, and $\theta_0 \geq 1$ sufficiently large, such that for all $k = 0, \ldots, n-1$ and all integer $s \in [6, \tilde{\alpha} - 4]$, one has*

$$\|e_k\|_{H^s_{*,\gamma}(Q^+_T)} + \|\hat{e}_k\|_{H^s_\gamma(Q^-_T)} \leq C \, \delta^2 \, \theta_k^{L(s)-1} \, \Delta_k \,,$$
$$\|\tilde{e}_k\|_{H^s_\gamma(\omega_T)} \leq C \, \delta^2 \, \theta_k^{s+7-2\alpha} \, \Delta_k \,, \quad (225)$$

*where $L(s) := \max\{(s + 2 - \alpha)_+ + 16 - 2\alpha; s + 12 - 2\alpha\}$.*

*Proof.* We recall that $e_k, \hat{e}_k, \tilde{e}_k$ are defined in (162) as the sum of all the error terms of the $k$-th step. Adding the estimates (179), (183), (219), (223) and (224) we obtain (225). $\quad\square$

The preceeding Lemma immediately yields the estimate of the accumulated errors $E_n$, and $\hat{E}_n$:

**Lemma 34.** *Let $\alpha \geq 14$, $\tilde{\alpha} = \alpha + 7$. There exist $\delta > 0$, $T > 0$ sufficiently small, $\theta_0 \geq 1$ sufficiently large, such that*

$$\|E_n\|_{H^{\alpha+3}_{*,\gamma}(Q^+_T)} + \|\hat{E}_n\|_{H^{\alpha+3}_\gamma(Q^-_T)} + \|\tilde{E}_n\|_{H^{\alpha+3}_\gamma(\omega_T)} \leq C \, \delta^2 \, \theta_n \,. \quad (226)$$

*Proof.* One can check that $L(\alpha + 3) \leq 1$ if $\alpha \geq 14$. Moreover, in order to apply (225) for $s = \alpha + 3$ one needs $\alpha + 3 \leq \tilde{\alpha} - 4$; the best choice is $\alpha + 3 = \tilde{\alpha} - 4$, which explains why $\tilde{\alpha} = \alpha + 7$. It follows from (225) that

$$\|E_n\|_{H^{\alpha+3}_{*,\gamma}(Q^+_T)} + \|\hat{E}_n\|_{H^{\alpha+3}_\gamma(Q^-_T)} + \|\tilde{E}_n\|_{H^{\alpha+3}_\gamma(\omega_T)} \leq C \, \delta^2 \sum_{k=0}^{n-1} \theta_k^{L(\alpha+3)-1} \, \Delta_k \leq C \, \delta^2 \, \theta_n.$$
$$\square$$

Going on with the iteration scheme, the next Lemma gives the estimates of the source terms $f_n, \hat{f}_n, \tilde{f}_n$, defined by equations (164). Notice that only the first three components of $\tilde{f}_n$ may be different from zero.

**Lemma 35.** *Let $\alpha \geq 14$, and let $\tilde{\alpha}$ be given as in Lemma 34. There exist $\delta > 0, T > 0$ sufficiently small and $\theta_0 \geq 1$ sufficiently large, such that for all integer $s \in [6, \tilde{\alpha} + 1]$, one has*

$$\|f_n\|_{H^s_{*,\gamma}(Q_T^+)} \leq C \, \Delta_n \{\theta_n^{s-\alpha-3}\|f^a\|_{H^{\alpha+2}_{*,\gamma}(Q_T^+)} + \delta^2 \, (\theta_n^{s-\alpha-3} + \theta_n^{L(s)-1})\} \,, \tag{227a}$$

$$\|\hat{f}_n\|_{H^s_\gamma(Q_T^-)} + \|\tilde{f}_n\|_{H^s_\gamma(\omega_T)} \leq C \, \delta^2 \, \Delta_n \, (\theta_n^{s-\alpha-3} + \theta_n^{L(s)-1}) \,. \tag{227b}$$

*Proof.* From (164) we have

$$f_n = (S_{\theta_n} - S_{\theta_{n-1}})f^a - (S_{\theta_n} - S_{\theta_{n-1}})E_{n-1} - S_{\theta_n}e_{n-1} \,.$$

Using (155), (225) and (226) gives (227a), with $\Delta_{n-1}, \theta_{n-1}$ instead of $\Delta_n, \theta_n$. Using $\theta_{n-1} \leq \theta_n \leq \sqrt{2}\theta_{n-1}$, and $\Delta_{n-1} \leq 3\,\Delta_n$, yields (227a). Estimate (227b) follows in the same way. $\qquad\square$

We now consider problem (165), that gives the solution $(\delta\dot{V}_n, \delta\dot{\mathcal{K}}_n, \delta\psi_n)$. Then we find $\Psi_{n+1}$, and consequently $(\delta V_n, \delta\mathcal{K}_n, \delta\psi_n)$:

**Lemma 36.** *Assume $\alpha \geq 14$. If $\delta > 0$ and $T > 0$ are sufficiently small, $\theta_0 \geq 1$ is sufficiently large, then for all $6 \leq s \leq \tilde{\alpha}$, one has*

$$\|\delta V_n\|_{H^s_{*,\gamma}(Q_T^+)} + \|\delta\mathcal{K}_n\|_{H^s_\gamma(Q_T^-)} + \|\delta\Psi_n\|_{H^{s+1}_\gamma(Q_T)} + \|\delta\psi_n\|_{H^{s+1/2}_\gamma(\omega_T)} \leq \delta \, \theta_n^{s-\alpha-1} \, \Delta_n \,. \tag{228}$$

*Proof.* Let us consider problem (165). This problem has the form (37), (38), i.e. in explicit form (39); thus it is equivalent to (85) and will be solved by applying Theorem 16. We first notice that $U^a + V_{n+1/2}, \mathcal{H}^a + \mathcal{K}^{n+1/2}, \varphi^a + \psi_{n+1/2}$ satisfy the required constraints (29)–(34), (61). In order to apply Theorem 16, we verify

$$\|U^a + V_{n+1/2}\|_{H^9_{*,\gamma}(Q_T^+)} + \|\mathcal{H}^a + \mathcal{K}^{n+1/2}\|^2_{H^9_\gamma(Q_T^-)} + \|\varphi^a + \psi_{n+1/2}\|_{H^{9.5}_\gamma(\omega_T)} \leq K,$$
$$\|\varphi^a + \psi_{n+1/2}\|_{\mathcal{C}([0,T];H^{2.5}(\Gamma))} \leq \epsilon_0,$$

by means of (172), (174a), (184), (185) and taking $\delta > 0$ and $T > 0$ sufficiently small (here we only use $\alpha \geq 11$). Thus we may apply our tame estimate (88) and obtain

$$\|\delta\dot{V}_n\|_{H^s_{*,\gamma}(Q_T^+)} + \|\delta\dot{\mathcal{K}}_n\|_{H^s_\gamma(Q_T^-)} + \|\delta\psi_n\|_{H^{s+1/2}_\gamma(\omega_T)} \leq C \left\{ \|f_n\|_{H^{s+1}_{*,\gamma}(Q_T^+)} + \|\hat{f}_n\|_{H^{s+1}_\gamma(Q_T^-)} + \|\tilde{f}_n\|_{H^{s+1}_\gamma(\omega_T)} \right.$$
$$+ \left(\|f_n\|_{H^8_{*,\gamma}(Q_T^+)} + \|\hat{f}_n\|_{H^8_\gamma(Q_T^-)} + \|\tilde{f}_n\|_{H^8_\gamma(\omega_T)}\right) \times$$
$$\left. \times \left(\|U^a + V_{n+1/2}\|_{H^{s+2}_{*,\gamma}(Q_T^+)} + \|\mathcal{H}^a + \mathcal{K}^{n+1/2}\|_{H^{s+2}_\gamma(Q_T^-)} + \|\varphi^a + \psi_{n+1/2}\|_{H^{s+2.5}_\gamma(\omega_T)}\right) \right\} \,. \tag{229}$$

On the other hand, from (158) it follows

$$\|\delta V_n\|_{H^s_{*,\gamma}(Q_T^+)} + \|\delta\mathcal{K}_n\|_{H^s_\gamma(Q_T^-)} \leq \|\delta\dot{V}_n\|_{H^s_{*,\gamma}(Q_T^+)} + \|\delta\dot{\mathcal{K}}_n\|_{H^s_\gamma(Q_T^-)}$$
$$+ C\{\|\delta\Psi_n\|_{H^s_\gamma(Q_T)} + \|\delta\Psi_n\|_{H^6_\gamma(Q_T)}(\|U^a + V_{n+1/2}\|_{H^{s+2}_{*,\gamma}(Q_T^+)}$$
$$+ \|\mathcal{H}^a + \mathcal{K}^{n+1/2}\|_{H^{s+1}_\gamma(Q_T^-)} + \|\Psi^a + \Psi_{n+1/2}\|_{H^{s+1}_\gamma(Q_T)})\} \,. \tag{230}$$

From Lemma 2 we have

$$\|\delta\Psi_n\|_{H^s_\gamma(Q_T)} \leq C\|\delta\psi_n\|_{H^{s-1/2}_\gamma(\omega_T)}. \tag{231}$$

Choosing $s = 5$ in (229) gives from (231)

$$\|\delta\Psi_n\|_{H^6_\gamma(Q_T)} \leq C \left(\|f_n\|_{H^8_{*,\gamma}(Q_T^+)} + \|\hat{f}_n\|_{H^8_\gamma(Q_T^-)} + \|\tilde{f}_n\|_{H^8_\gamma(\omega_T)}\right) \times$$
$$\times \left(1 + \|U^a + V_{n+1/2}\|_{H^7_{*,\gamma}(Q_T^+)} + \|\mathcal{H}^a + \mathcal{K}^{n+1/2}\|_{H^7_\gamma(Q_T^-)} + \|\varphi^a + \psi_{n+1/2}\|_{H^{7.5}_\gamma(\omega_T)}\right) \,. \tag{232}$$

Therefore, we can combine (229)–(232) and eventually obtain

$$\|\delta V_n\|_{H^s_{*,\gamma}(Q_T^+)} + \|\delta\mathcal{K}_n\|_{H^s_\gamma(Q_T^-)} + \|\delta\psi_n\|_{H^{s+1/2}_\gamma(\omega_T)} \leq C\left\{\|f_n\|_{H^{s+1}_{*,\gamma}(Q_T^+)} + \|\hat{f}_n\|_{H^{s+1}_\gamma(Q_T^-)} + \|\tilde{f}_n\|_{H^{s+1}_\gamma(\omega_T)}\right.$$
$$+ \left(\|f_n\|_{H^s_{*,\gamma}(Q_T^+)} + \|\hat{f}_n\|_{H^s_\gamma(Q_T^-)} + \|\tilde{f}_n\|_{H^s_\gamma(\omega_T)}\right) \times$$
$$\left. \times \left(\|U^a + V_{n+1/2}\|_{H^{s+2}_{*,\gamma}(Q_T^+)} + \|\mathcal{H}^a + \mathcal{K}^{n+1/2}\|_{H^{s+2}_\gamma(Q_T^-)} + \|\varphi^a + \psi_{n+1/2}\|_{H^{s+2.5}_\gamma(\omega_T)}\right) \right\}. \quad (233)$$

for all integer $s \in [6, \tilde{\alpha}]$. The remaining part of the work is to estimate the right-hand side of (233). Using Lemma 35, (174a) and Proposition 28, (233) becomes

$$\|\delta V_n\|_{H^s_{*,\gamma}(Q_T^+)} + \|\delta\mathcal{K}_n\|_{H^s_\gamma(Q_T^-)} + \|\delta\psi_n\|_{H^{s+1/2}_\gamma(\omega_T)}$$
$$\leq C\left\{\theta_n^{s-\alpha-2}\left(\|f^a\|_{H^{s+2}_{*,\gamma}(Q_T^+)} + \delta^2\right) + \delta^2\,\theta_n^{L(s+1)-1}\right\}\Delta_n$$
$$+ C\,\delta\,\Delta_n\left\{\theta_n^{5-\alpha}\left(\|f^a\|_{H^{s+2}_{*,\gamma}(Q_T^+)} + \delta^2\right) + \delta^2\,\theta_n^{19-2\alpha}\right\}\left(\theta_n^{(s+2-\alpha)_+} + \theta_n^{s+4-\alpha}\right). \quad (234)$$

One checks that, for $\alpha \geq 14$, and $s \in [6, \tilde{\alpha}]$, the following inequalities hold true:

$$L(s+1) \leq s - \alpha, \qquad\qquad (s+2-\alpha)_+ + 5 - \alpha \leq s - \alpha - 1,$$
$$s + 9 - 2\alpha \leq s - \alpha - 1, \qquad\qquad (s+2-\alpha)_+ + 19 - 2\alpha \leq s - \alpha - 1,$$
$$s + 23 - 3\alpha \leq s - \alpha - 1.$$

From (234), we thus obtain

$$\|\delta V_n\|_{H^s_{*,\gamma}(Q_T^+)} + \|\delta\mathcal{K}_n\|_{H^s_\gamma(Q_T^-)} + \|\delta\psi_n\|_{H^{s+1/2}_\gamma(\omega_T)}$$
$$\leq C\left\{\|f^a\|_{H^{s+2}_{*,\gamma}(Q_T^+)} + \delta^2\right\}\theta_n^{s-\alpha-1}\,\Delta_n \leq C\left\{\delta'(T) + \delta^2\right\}\theta_n^{s-\alpha-1}\,\Delta_n,$$

and (228) follows for $\delta > 0$ and $T > 0$ sufficiently small. $\qquad\square$

From (163), (164) the components of $\tilde{f}_n$ defined on $\omega_T^\pm$ are zero. Then the resolution of (165) gives $\delta\dot{V}_{n,2} = 0$ on $\omega_T^+$, $n \times \delta\dot{\mathcal{K}}_n = 0$ on $\omega_T^-$. From (158) the same is true for $\delta V_n, \delta\mathcal{K}_n$ because by construction $\delta\Psi_n = 0$ at $\omega_T^\pm$. This shows that (156) holds for $k = n + 1$ as well.

We now check the three remaining inequalities in $(H_n)$.

**Lemma 37.** *Assume $\alpha \geq 14$. If $\delta > 0$ and $T > 0$ are sufficiently small, and $\theta_0 \geq 1$ is sufficiently large, then for all $6 \leq s \leq \tilde{\alpha} - 2$, one has*

$$\|\mathcal{L}(V_n, \Psi_n) - f^a\|_{H^s_{*,\gamma}(Q_T^+)} \leq 2\,\delta\,\theta_n^{s-\alpha-1}, \quad (235a)$$

$$\|\mathcal{E}(\mathcal{K}_n, \Psi_n)\|_{H^s_\gamma(Q_T^-)} \leq 2\,\delta\,\theta_n^{s-\alpha-1}, \quad (235b)$$

$$\|\mathcal{B}(V_n, \mathcal{K}_n, \psi_n)_{1,2,3}\|_{H^s_\gamma(\omega_T)} \leq \delta\,\theta_n^{s-\alpha-1}. \quad (235c)$$

*Proof.* Recall that, by summing the relations (166), we have

$$\mathcal{L}(V_n, \Psi_n) - f^a = (S_{\theta_{n-1}} - I)f^a + (I - S_{\theta_{n-1}})E_{n-1} + e_{n-1}.$$

The proof of (235a) then follows by applying (155), (225) and (226), provided that $\delta > 0$ and $T > 0$ are taken sufficiently small. The proof of (235b) and (235c) is similar. $\qquad\square$

Lemmata 36 and 37 show that $(H_{n-1})$ implies $(H_n)$ provided that $\alpha \geq 14$, $\tilde{\alpha} = \alpha + 7$, (172) holds, $\delta > 0$ is small enough, $T > 0$ is small enough, and $\theta_0 \geq 1$ is large enough. We fix $\alpha$, $\tilde{\alpha}$, $\delta > 0$ and $\theta_0 \geq 1$, and we finally prove $(H_0)$.

**Lemma 38.** *If $T > 0$ is sufficiently small, then property $(H_0)$ holds.*

*Proof.* Recall that $V_0 = \mathcal{K}_0 = \Psi_0 = \psi_0 = 0$. Thanks to the properties of the approximate solution (see Lemma 21), we see that $U^a + V_0, \mathcal{H}^a + \mathcal{K}_0, \Psi^a + \Psi_0, \varphi^a + \psi_0$ satisfy the constraints (29)–(33) and (61). Consequently, the contruction of Proposition 28 yields $V_{1/2} = \mathcal{K}_{1/2} = \Psi_{1/2} = \psi_{1/2} = 0$. Consider the problem

$$\begin{aligned}
&\mathbb{P}'_e(U^a, \Psi^a)\, \delta\dot{V}_0 = S_{\theta_0} f^a && \text{in } Q_T^+, \\
&\mathbb{V}'_e(\mathcal{H}^a, \Psi^a)\, \delta\dot{\mathcal{K}}_0 = 0 && \text{in } Q_T^-, \\
&\mathbb{B}_{1/2}(\delta\dot{V}_0, \delta\dot{\mathcal{K}}_0, \delta\psi_0) = 0 && \text{on } \omega_T^3 \times \omega_T^{\pm}, \\
&\delta\dot{V}_0 = 0, \quad \delta\dot{\mathcal{K}}_0 = 0, \quad \delta\psi_0 = 0 && \text{for } t < 0.
\end{aligned}$$

Because of (143), we may apply (88) and obtain

$$\begin{aligned}
\|\delta\dot{V}_0\|_{H^s_{*,\gamma}(Q_T^+)} + \|\delta\dot{\mathcal{K}}_0\|_{H^s_\gamma(Q_T^-)} + \|\delta\psi_0\|_{H^{s+1/2}_\gamma(\omega_T)} \leq C\Big\{ \|S_{\theta_0}f^a\|_{H^{s+1}_{*,\gamma}(Q_T^+)} \\
+ \|S_{\theta_0}f^a\|_{H^8_{*,\gamma}(Q_T^+)} \left( \|U^a\|_{H^{s+2}_{*,\gamma}(Q_T^+)} + \|\mathcal{H}^a\|_{H^{s+2}_\gamma(Q_T^+)} + \|\varphi^a\|_{H^{s+2.5}_\gamma(\omega_T)} \right) \Big\}. \quad (236)
\end{aligned}$$

Then we find $\delta\Psi_0$ from $\delta\psi_0$ by Lemma 2. From (158) we finally obtain:

$$\begin{aligned}
\|\delta V_0\|_{H^s_{*,\gamma}(Q_T^+)} + \|\delta\mathcal{K}_0\|_{H^s_\gamma(Q_T^-)} + \|\delta\psi_0\|_{H^{s+1/2}_\gamma(\omega_T)} \leq C\Big\{ \|S_{\theta_0}f^a\|_{H^{s+1}_{*,\gamma}(Q_T^+)} \\
+ \|S_{\theta_0}f^a\|_{H^8_{*,\gamma}(Q_T^+)} \left( \|U^a\|_{H^{s+2}_{*,\gamma}(Q_T^+)} + \|\mathcal{H}^a\|_{H^{s+2}_\gamma(Q_T^+)} + \|\varphi^a\|_{H^{s+2.5}_\gamma(Q_T)} \right) \Big\}, \quad (237)
\end{aligned}$$

for all integer $s \in [6, \tilde{\alpha}]$. If $T > 0$ is sufficiently small, then

$$\|\delta V_0\|_{H^s_{*,\gamma}(Q_T^+)} + \|\delta\mathcal{K}_0\|_{H^s_\gamma(Q_T^-)} + \|\delta\psi_0\|_{H^{s+1/2}_\gamma(\omega_T)} \leq \delta\, \theta_0^{s-\alpha-1}\Delta_0, \qquad 6 \leq s \leq \tilde{\alpha}.$$

The other inequalities in $(H_0)$ are readily satisfied by taking $T > 0$ small enough. $\qquad\square$

From Lemmata $36 - 38$, we get that $(H_n)$ holds for every $n \geq 0$, provided that the parameters are well-chosen.

*Conclusion of the proof of the existence of smooth solutions in Theorem 5.*

Given an integer $\alpha \geq 14$, in agreement with the requirements of Lemma 34, we set $\tilde{\alpha} = \alpha + 7$. Let $m = \alpha - 1 \geq 13$. Let us consider initial data $U_0 \in H^{m+9.5}(\Omega^+)$, $\mathcal{H}^0 \in H^{m+9.5}(\Omega^-)$, and $\varphi_0 \in H^{m+10}(\Gamma)$ that satisfy (8), (22), (128)–(130) and are compatible up to order $m + 9$ in the sense of Definition 20. Then, by Lemma 21 (with $k = \alpha + 9 = m + 10$) we may find an approximate solution $(U^a, \mathcal{H}^a, \varphi^a)$ such that $U^a \in H^{m+10}(Q^+)$, $\mathcal{H}^a \in H^{m+10}(Q^-)$, $\varphi^a \in H^{m+10.5}(\omega)$, with the properties listed there and (172). If $\delta > 0$ and $T > 0$ are small enough, and $\theta_0 \geq 1$ is large enough, one gets Lemmata 36, 37, 38. Hence the property $(H_n)$ holds true for all $n$. In particular, it follows that

$$\sum_{n\geq 0} \|\delta V_n\|_{H^m_{*,\gamma}(Q_T^+)} + \|\delta\mathcal{K}_n\|_{H^m_\gamma(Q_T^-)} + \|\delta\psi_n\|_{H^{m+1/2}_\gamma(\omega_T)} < +\infty,$$

so the sequences $(V_n)$, $(\mathcal{K}_n)$, $(\psi_n)$, converge in $H^m_{*,\gamma}(Q_T^+) \times H^m_\gamma(Q_T^-) \times H^{m+1/2}_\gamma(\omega_T)$ towards some limits $V, \mathcal{K}, \psi$. Passing to the limit in (235), for $s = m$, we obtain (153). Therefore $U = U^a + V, \mathcal{H} = \mathcal{H}^a + \mathcal{K}, \varphi = \varphi^a + \psi$ is a solution on $Q_T^+ \times Q_T^-$ of (17)–(19).

The proof of the existence part of Theorem 5 is complete.

## 13. Proof of the uniqueness of a smooth solution

Having in hand the existence of a smooth solution $(U, \mathcal{H}, \varphi)$ from Theorem 5, our goal now is to prove its uniqueness. Let, on the contrary, the exist one more solution $(U', \mathcal{H}', \varphi')$ of problem (17)–(19). Omitting calculations, for the differences

$$\widetilde{U} = U - U', \quad \widetilde{\mathcal{H}} = \mathcal{H} - \mathcal{H}', \quad \tilde{\varphi} = \varphi - \varphi'$$

we obtain the following initial boundary value problem:

$$P(U, \Psi)\widetilde{U} - \big\{ P(U, \Psi)\widetilde{\Psi} \big\} \frac{\partial_1 U'}{\partial_1 \Phi'_1} + \mathcal{R} = 0 \qquad \text{in } Q_T^+, \qquad (238)$$

$$\mathbb{V}(\widetilde{\mathcal{H}}, \Psi) + \begin{pmatrix} \nabla \mathcal{H}'_1 \times \nabla \widetilde{\Psi} \\ \nabla \times \begin{pmatrix} 0 \\ -\mathcal{H}'_3 \\ \mathcal{H}'_2 \end{pmatrix} \cdot \nabla \widetilde{\Psi} \end{pmatrix} = 0 \qquad \text{in } Q^-_T, \tag{239}$$

$$\begin{aligned} &\partial_t \tilde{\varphi} + v'_2 \partial_2 \tilde{\varphi} + v'_3 \partial_3 \tilde{\varphi} - \tilde{v}_N = 0, \\ &\tilde{q} - \mathcal{H}' \cdot \widetilde{\mathcal{H}} = \mathfrak{R}, \\ &\widetilde{\mathcal{H}}_N - \mathcal{H}'_2 \partial_2 \tilde{\varphi} - \mathcal{H}'_3 \partial_3 \tilde{\varphi} = 0 \qquad \text{on } \omega_T, \end{aligned} \tag{240}$$

$$\tilde{v}_1 = 0 \quad \text{on } \omega^+_T, \qquad \nu \times \widetilde{\mathcal{H}} = 0 \quad \text{on } \omega^-_T, \tag{241}$$

and we may assume that

$$(\widetilde{U}, \widetilde{\mathcal{H}}, \tilde{\varphi}) \quad \text{for } t < 0 \tag{242}$$

because we have the trivial initial data $(\widetilde{U}, \widetilde{\mathcal{H}}, \tilde{\varphi})|_{t=0} = 0$. Here

$$\mathcal{R} = (P(U, \Psi') - P(U', \Psi')) U', \qquad \mathfrak{R} = \frac{1}{2} |\widetilde{\mathcal{H}}|^2,$$

$$\tilde{v}_N = \tilde{v}_1 - \tilde{v}_2 \partial_2 \Psi - \tilde{v}_3 \partial_3 \Psi, \quad \widetilde{H}_N = \widetilde{H}_1 - \widetilde{H}_2 \partial_2 \Psi - \widetilde{H}_3 \partial_3 \Psi, \quad \widetilde{\mathcal{H}}_N = \widetilde{\mathcal{H}}_1 - \widetilde{\mathcal{H}}_2 \partial_2 \Psi - \widetilde{\mathcal{H}}_3 \partial_3 \Psi,$$

$$\widetilde{\Psi} = \Psi - \Psi', \quad \tilde{q} = q - q', \quad \tilde{v}_j = v_j - v'_j, \quad \widetilde{H}_j = H_j - H'_j, \quad \text{etc.}$$

Note also that the functions $\Psi'$ and $\Phi'$ are defined through the function $\varphi'$ exactly in the same manner as the functions $\Psi$ and $\Phi$ entering (17)–(19) are defined through $\varphi$.

Since all the terms entering the differential operator $P(U, \Psi') - P(U', \Psi')$ contain the differences of matrices $A_\alpha(U) - A_\alpha(U')$ (for a certain $\alpha = \overline{0,3}$), by using the mean value theorem, we can represent the rest term $\mathcal{R}$ in (238) in the form

$$\mathcal{R} = \widehat{\mathcal{C}} \widetilde{U} \tag{243}$$

where $\widehat{\mathcal{C}} = \mathcal{C}(U^*, U', \Psi')$ and the matrix $\mathcal{C}$ depends on the space-time gradients of $U'$, $\Psi'$ as well as on the vector of "mean values" $U^*$ whose components can be estimated through the norms of $U$ and $U'$. It is worth noting that the both solutions $(U, \mathcal{H}, \varphi)$ and $(U', \mathcal{H}', \varphi')$ satisfy constraints (20) and (21). This gives the following equations for the differences:

$$\operatorname{div} \tilde{h} + \nabla \times \begin{pmatrix} 0 \\ -H'_3 \\ H'_2 \end{pmatrix} \cdot \nabla \widetilde{\Psi} = 0 \quad \text{in } Q^+_T, \tag{244}$$

$$\widetilde{H}_N - H'_2 \partial_2 \tilde{\varphi} - H'_3 \partial_3 \tilde{\varphi} = 0 \quad \text{on } \omega_T, \qquad \widetilde{H}_1 = 0 \quad \text{on } \omega^+_T, \tag{245}$$

where $\tilde{h} = (\widetilde{H}_N, \widetilde{H}_2 \partial_1 \Phi_1, \widetilde{H}_3 \partial_1 \Phi_1)$.

Now, as for the linearized problem in Section 4, we pass in (238)–(245) to the "good unknown"

$$\dot{U} := \widetilde{U} - \frac{\widetilde{\Psi}}{\partial_1 \Phi'_1} \partial_1 U', \quad \dot{\mathcal{H}} := \widetilde{\mathcal{H}} - \frac{\widetilde{\Psi}}{\partial_1 \Phi'_1} \partial_1 \mathcal{H}' \tag{246}$$

for the differences of solutions (cf. (35)). Taking into account (243) and omitting detailed calculations, we rewrite (238)–(242) as follows (cf. (39)):

$$\widehat{A}_0 \partial_t \dot{U} + \sum_{j=1}^{3} \widehat{A}_j \partial_j \dot{U} + \widehat{\mathcal{C}} \dot{U} = f \qquad \text{in } Q_T^+, \tag{247a}$$

$$\nabla \times \dot{\mathfrak{H}} = \chi, \quad \text{div } \dot{\mathfrak{h}} = \Xi \qquad \text{in } Q_T^-, \tag{247b}$$

$$\partial_t \tilde{\varphi} = \dot{v}_N - v_2' \partial_2 \tilde{\varphi} - v_3' \partial_3 \tilde{\varphi} + \tilde{\varphi} \, \partial_1 \dot{v}_N, \tag{247c}$$

$$\dot{q} = \mathcal{H}' \cdot \dot{\mathcal{H}} - [\partial_1 q'] \tilde{\varphi} + g_2, \tag{247d}$$

$$\dot{\mathcal{H}}_N = \partial_2 \big( \mathcal{H}_2' \tilde{\varphi} \big) + \partial_3 \big( \mathcal{H}_3' \tilde{\varphi} \big) + g_3 \qquad \text{on } \omega_T, \tag{247e}$$

$$\dot{v}_1 = 0 \qquad \text{on } \omega_T^+, \qquad \nu \times \dot{\mathcal{H}} = 0 \qquad \text{on } \omega_T^-, \tag{247f}$$

$$(\dot{U}, \dot{\mathcal{H}}, \tilde{\varphi}) = 0 \qquad \text{for } t < 0, \tag{247g}$$

where

$$\widehat{A}_\alpha =: A_\alpha(U), \quad \alpha = 0, 2, 3, \quad \widehat{A}_1 =: \widetilde{A}_1(U, \Psi), \quad f = \hat{a} \widetilde{\Psi},$$

$$\dot{\mathfrak{H}} = (\dot{\mathcal{H}}_1 \partial_1 \Phi_1, \dot{\mathcal{H}}_{\tau_2}, \dot{\mathcal{H}}_{\tau_3}), \quad \dot{\mathfrak{h}} = (\dot{\mathcal{H}}_N, \dot{\mathcal{H}}_2 \partial_1 \Phi_1, \dot{\mathcal{H}}_3 \partial_1 \Phi_1),$$

$$\dot{\mathcal{H}}_N = \dot{\mathcal{H}}_1 - \dot{\mathcal{H}}_2 \partial_2 \Psi - \dot{\mathcal{H}}_3 \partial_3 \Psi, \quad \dot{\mathcal{H}}_{\tau_i} = \dot{\mathcal{H}}_1 \partial_i \Psi + \dot{\mathcal{H}}_i, \quad i = 2, 3,$$

$$[\partial_1 q'] = (\partial_1 q')|_{x_1=0} - (\mathcal{H}' \cdot \partial_1 \mathcal{H}')|_{x_1=0}, \quad \hat{v}_N = v_1' - v_2' \partial_2 \Psi - v_3' \partial_3 \Psi,$$

$$g_2 := \mathfrak{R} = \frac{1}{2} |\widetilde{\mathcal{H}}|^2 = \frac{1}{2} |\dot{\mathcal{H}} + \tilde{\varphi} \partial_1 \mathcal{H}'|^2,$$

$$\begin{pmatrix} \chi \\ \Xi \end{pmatrix} := -\frac{\widetilde{\Psi}}{\partial_1 \Phi_1'} \partial_1 \{ \mathbb{V}(\mathcal{H}', \Psi) \}, \quad g_3 := -\tilde{\varphi} \, (\text{div } \hat{\mathfrak{h}})|_{\omega_T},$$

$$\hat{\mathfrak{h}} = (\widehat{\mathcal{H}}_N, \mathcal{H}_2' \partial_1 \Phi_1, \mathcal{H}_3' \partial_1 \Phi_1), \quad \widehat{\mathcal{H}}_N = \mathcal{H}_1' - \mathcal{H}_2' \partial_2 \Psi - \mathcal{H}_3' \partial_3 \Psi,$$

and the vector $\hat{a}$ appearing in the definition of the source term $f$ in (247a) depends on the space-time gradients of $U'$, $\Psi'$ as well as on the vector of "mean values" $U^*$, but its concrete form is of no interest. Moreover, (244) and (245) are rewritten as (cf. (43), (44))

$$\text{div } \dot{h} = r \quad \text{in } Q_T^+, \tag{248}$$

$$\dot{H}_N - H_2' \partial_2 \tilde{\varphi} - H_3' \partial_3 \tilde{\varphi} + \tilde{\varphi} \, \partial_1 \widehat{H}_N = 0 \quad \text{on } \omega_T, \qquad \dot{H}_1 = 0 \quad \text{on } \omega_T^+, \tag{249}$$

with

$$r = -\frac{\widetilde{\Psi}}{\partial_1 \Phi_1'} \text{div } \hat{h}, \quad \hat{h} = (\widehat{H}_N, H_2' \partial_1 \Phi_1, H_3' \partial_1 \Phi_1), \quad \widehat{H}_N = H_1' - H_2' \partial_2 \Psi - H_3' \partial_3 \Psi.$$

As in (85), system (247a) can be rewritten in terms of the vector $\dot{\mathcal{U}} = (\dot{q}, \dot{u}, \dot{h}, \dot{S})$ with a corresponding new source term $\tilde{f}$ which, in view of Lemmata 2 and 3, can be estimated as

$$\|\tilde{f}_\gamma\|^2_{H^1_{tan,\gamma}(Q_T^+)} \leq C \|f_\gamma\|^2_{H^2_{tan,\gamma}(Q_T^+)} \leq C \|\widetilde{\Psi}_\gamma\|^2_{H^2_{tan,\gamma}(Q_T^+)} \leq C \|\tilde{\varphi}_\gamma\|^2_{H^{3/2}(\omega_T)}. \tag{250}$$

Here and below $C$ stays for different constants depending on Sobolev's norms of the solutions $(U, \mathcal{H}, \varphi)$ and $(U', \mathcal{H}', \varphi')$. The role of the coefficients for the reduced system for the vector $\dot{\mathcal{U}} = (\dot{q}, \dot{u}, \dot{h}, \dot{S})$ is played by the solution $(U, \mathcal{H}, \varphi)$. But, it is only important that the boundary matrix for this system calculated at the boundary is the matrix $\mathcal{E}_{12}$ (see (58)). Note also that the "coefficients" $\partial_1 \dot{v}_N$ and $\partial_1 \widehat{H}_N$ in (247c) and (249) are unimportant in the process of getting a priori estimates for $(\dot{U}, \dot{\mathcal{H}}, \tilde{\varphi})$ whereas the role of the rest coefficients in the boundary conditions (247c)–(247e), (249) is played by $(U', \mathcal{H}')$.

That is, problem (247) considered as a problem for $(\dot{U}, \dot{\mathcal{H}}, \tilde{\varphi})$ has the same form and same properties as the linear problem (85). It should be only noted that in the process of reduction of problem (247) to that with $g_2 = g_3 = 0$, $\chi = 0$ and $\Xi = 0$ we choose the zero "shifting" function for $\dot{H}$. But, then, differently from (53), we will have the non-zero $r$ in (248) for the reduced problem which is a counterpart of problem (52)–(54). However, in [35] equation (53) was used only for estimating the normal derivative $\partial_1 h_1$ of the noncharacteristic unknown $h_1$ through the tangential derivatives $\partial_2 h_2$ and $\partial_3 h_3$. In our case with the non-zero $r$ in (248), we have

$$\|\partial_1 \dot{h}_{1\gamma}\|_{L^2(Q^+)} \leq \|\partial_2 \dot{h}_{2\gamma}\|_{L^2(Q^+)} + \|\partial_3 \dot{h}_{3\gamma}\|_{L^2(Q^+)} + C\|\tilde{\varphi}_\gamma\|_{L^2(\omega)}^2,$$

but the last $L^2$ norm in the above inequality does not affect the derivation of a basic priori estimate for problem (247) (see [35]).

Thus, as in Theorem 15, we can derive for problem (247) the a priori estimate (86) with $g_1 = g_4 = g_5 = 0$:

$$\gamma \left( \|\dot{\mathcal{U}}_\gamma\|_{H^1_{tan,\gamma}(Q^+_T)}^2 + \|\dot{\mathfrak{H}}_\gamma\|_{H^1_\gamma(Q^-_T)}^2 + \|(\dot{q}_\gamma, \dot{u}_{1\gamma}, \dot{h}_{1\gamma})|_{\omega_T}\|_{H^{1/2}_\gamma(\omega_T)}^2 \right.$$
$$\left. + \|\dot{\mathfrak{H}}_\gamma|_{\omega_T}\|_{H^{1/2}_\gamma(\omega_T)}^2 + \|\tilde{\varphi}_\gamma\|_{H^{3/2}_\gamma(\omega_T)}^2 \right)$$
$$\leq \frac{C}{\gamma} \left( \|\tilde{f}_\gamma\|_{H^1_{tan,\gamma}(Q^+_T)}^2 + \|\chi_\gamma, \Xi_\gamma\|_{H^2_\gamma(Q^-_T)}^2 + \|g_{2\gamma}\|_{H^{3/2}_\gamma(\omega_T)}^2 + \|g_{3\gamma}\|_{H^2_\gamma(\omega_T)}^2 \right). \quad (251)$$

Taking into account Lemmata 2 and 3 and the exact form of the source terms $\chi$ and $\Xi$, we have

$$\|\chi_\gamma, \Xi_\gamma\|_{H^2_\gamma(Q^-_T)}^2 \leq C\|\tilde{\varphi}_\gamma\|_{H^{3/2}_\gamma(\omega_T)}^2. \quad (252)$$

The nonlinear term $g_2$ can be estimated as follows:

$$\|g_{2\gamma}\|_{H^{3/2}_\gamma(\omega_T)}^2 \leq C \left( \|\widetilde{\mathcal{H}} \cdot (\dot{\mathcal{H}}_\gamma + \tilde{\varphi}_\gamma \partial_1 \mathcal{H}')\|_{H^{1/2}_\gamma(\omega_T)}^2 + \sum_{j=0,2,3} \|Z_j \widetilde{\mathcal{H}} \cdot (\dot{\mathcal{H}}_\gamma + \tilde{\varphi}_\gamma \partial_1 \mathcal{H}')\|_{H^{1/2}_\gamma(\omega_T)}^2 \right)$$
$$\leq C \left( \|\dot{\mathfrak{H}}_\gamma|_{\omega_T}\|_{H^{1/2}_\gamma(\omega_T)}^2 + \|\tilde{\varphi}_\gamma\|_{H^{1/2}_\gamma(\omega_T)}^2 \right). \quad (253)$$

At last, using the decomposition

$$(\operatorname{div} \hat{\mathfrak{h}})|_{\omega_T} = (\operatorname{div} \mathfrak{h}')|_{\omega_T} - \sum_{j=2}^3 \partial_1 \mathcal{H}'_j|_{\omega_T} \partial_j \tilde{\varphi},$$

the fact that

$$\operatorname{div} \mathfrak{h}' = 0 \quad \text{in } Q^-_T,$$

where $\mathfrak{h}' = (\mathcal{H}'_N, \mathcal{H}'_2 \partial_1 \Phi'_1, \mathcal{H}'_3 \partial_1 \Phi'_1)$, $\mathcal{H}'_N = \mathcal{H}'_1 - \mathcal{H}'_2 \partial_2 \Psi' - \mathcal{H}'_3 \partial_3 \Psi'$, and the Leibniz rule

$$Z_k(\tilde{\varphi} \partial_j \tilde{\varphi}) = (Z_k \partial_j \tilde{\varphi})\tilde{\varphi} + (\partial_j \tilde{\varphi})Z_k \tilde{\varphi},$$
$$Z_m Z_k(\tilde{\varphi} \partial_j \tilde{\varphi}) = (Z_m Z_k \partial_j \tilde{\varphi})\tilde{\varphi} + (Z_k \partial_j \tilde{\varphi})Z_m \tilde{\varphi} + (Z_m \partial_j \tilde{\varphi})Z_k \tilde{\varphi} + (Z_m Z_k \tilde{\varphi})Z_j \tilde{\varphi} \quad (254)$$

$(m, k = 0, 2, 3)$, we estimate the source term $g_3$,

$$\|g_{3\gamma}\|_{H^2_\gamma(\omega_T)}^2 \leq C \sum_{j=2}^3 \|\tilde{\varphi}_\gamma \partial_j \tilde{\varphi}\|_{H^2_\gamma(\omega_T)}^2 \leq C\|\tilde{\varphi}_\gamma\|_{H^1_\gamma(\omega_T)}^2, \quad (255)$$

by treating the terms in parentheses in (254) as coefficients.

It follows from (250)–(253), (255)

$$\gamma \left( \|\dot{\mathcal{U}}_\gamma\|_{H^1_{tan,\gamma}(Q^+_T)}^2 + \|\dot{\mathfrak{H}}_\gamma\|_{H^1_\gamma(Q^-_T)}^2 + \|\dot{\mathfrak{H}}_\gamma|_{\omega_T}\|_{H^{1/2}_\gamma(\omega_T)}^2 + \|\tilde{\varphi}_\gamma\|_{H^{3/2}_\gamma(\omega_T)}^2 \right)$$
$$\leq \frac{C}{\gamma} \left( \|\dot{\mathfrak{H}}_\gamma|_{\omega_T}\|_{H^{1/2}_\gamma(\omega_T)}^2 + \|\tilde{\varphi}_\gamma\|_{H^{3/2}_\gamma(\omega_T)}^2 \right). \quad (256)$$

Absorbing the norms in the right-hand side of inequality (256) for $\gamma$ large enough, we get the estimate

$$\|\dot{\mathcal{U}}_\gamma\|_{H^1_{tan,\gamma}(Q^+_T)}^2 + \|\dot{\mathfrak{H}}_\gamma\|_{H^1_\gamma(Q^-_T)}^2 + \|\tilde{\varphi}_\gamma\|_{H^{3/2}_\gamma(\omega_T)}^2 \leq 0$$

which implies $\widetilde{U} = 0$, $\widetilde{\mathcal{H}} = 0$ and $\tilde{\varphi} = 0$, i.e., the uniqueness of the smooth solution $(U, \mathcal{H}, \varphi)$.

## Appendix A. Properties of anisotropic Sobolev spaces

The next theorems deal with the product of two functions in anisotropic Sobolev spaces.

**Theorem 39.** *Let $n \geq 2$ and $\gamma \geq 1$. Moreover let $m \geq 1$ be an integer and $s = \max\left\{m, \left[\frac{n+1}{2}\right] + 2\right\}$. For any $u \in H^m_{*,\gamma}(\mathbb{R}^n_+)$ and $v \in H^s_{*,\gamma}(\mathbb{R}^n_+)$ one has $uv \in H^m_{*,\gamma}(\mathbb{R}^n_+)$. Moreover, there exists a constant $C$ independent of $\gamma$ such that*

$$\gamma^{s-(n+1)/2}||uv||_{H^m_{*,\gamma}(\mathbb{R}^n_+)} \leq C||u||_{H^m_{*,\gamma}(\mathbb{R}^n_+)}||v||_{H^s_{*,\gamma}(\mathbb{R}^n_+)}, \qquad \forall \gamma \geq 1. \tag{257}$$

*Proof.* See [22], Theorem 34. □

Let us define the space

$$W^{1,\infty}_*(\mathbb{R}^n_+) = \{u \in L^\infty(\mathbb{R}^n_+) \,:\, Z_i u \in L^\infty(\mathbb{R}^n_+),\, i = 1\ldots, n\}, \tag{258}$$

equipped with its natural norm. We have the following Moser-type inequalities.

**Theorem 40.** *Let $n \geq 2$ and $\gamma \geq 1$. If $m$ is 1 or even, for all functions $u$ and $v$ in $H^m_{*,\gamma}(\mathbb{R}^n_+) \cap L^\infty(\mathbb{R}^n_+)$ one has*

$$\|uv\|_{H^m_{*,\gamma}(\mathbb{R}^n_+)} \leq C(\|u\|_{H^m_{*,\gamma}(\mathbb{R}^n_+)}\|v\|_{L^\infty(\mathbb{R}^n_+)} + \|u\|_{L^\infty(\mathbb{R}^n_+)}\|v\|_{H^m_{*,\gamma}(\mathbb{R}^n_+)}), \qquad \forall \gamma \geq 1. \tag{259}$$

*If $m \geq 3$ is odd, for all functions $u$ and $v$ in $H^m_{*,\gamma}(\mathbb{R}^n_+) \cap W^{1,\infty}_*(\mathbb{R}^n_+)$ one has*

$$\|uv\|_{H^m_{*,\gamma}(\mathbb{R}^n_+)} \leq C(\|u\|_{H^m_{*,\gamma}(\mathbb{R}^n_+)}\|v\|_{W^{1,\infty}_*(\mathbb{R}^n_+)} + \|u\|_{W^{1,\infty}_*(\mathbb{R}^n_+)}\|v\|_{H^m_{*,\gamma}(\mathbb{R}^n_+)}), \qquad \forall \gamma \geq 1. \tag{260}$$

*Proof.* See [24], Theorem B.3. □

**Theorem 41.** *Let $n \geq 2$ and $\gamma \geq 1$. For every integer $m \geq \left[\frac{n+1}{2}\right] + 1$ the continuous imbedding $H^m_{*,\gamma}(\mathbb{R}^n_+) \hookrightarrow C^0_B(\mathbb{R}^n_+)$ holds. Moreover, there exists a constant $C$ such that for every $u \in H^m_{*,\gamma}(\mathbb{R}^n_+)$*

$$\gamma^{m-(n+1)/2}||u||_{L^\infty(\mathbb{R}^n_+)} \leq C||u||_{H^m_{*,\gamma}(\mathbb{R}^n_+)} \qquad \forall \gamma \geq 1. \tag{261}$$

From Theorems 40 and 41 we get

**Corollary 42.** *For every even integer $m \geq \left[\frac{n+1}{2}\right] + 1$, for all functions $u$ and $v$ in $H^m_{*,\gamma}(\mathbb{R}^n_+)$ one has*

$$\gamma^{1/2}\|uv\|_{H^m_{*,\gamma}(\mathbb{R}^n_+)} \leq C(\|u\|_{H^m_{*,\gamma}(\mathbb{R}^n_+)}\|v\|_{H^{[\frac{n+1}{2}]+1}_{*,\gamma}(\mathbb{R}^n_+)} + \|u\|_{H^{[\frac{n+1}{2}]+1}_{*,\gamma}(\mathbb{R}^n_+)}\|v\|_{H^m_{*,\gamma}(\mathbb{R}^n_+)}), \qquad \forall \gamma \geq 1. \tag{262}$$

*For every odd integer $m \geq \left[\frac{n+1}{2}\right] + 2$ and for all functions $u$ and $v$ in $H^m_{*,\gamma}(\mathbb{R}^n_+)$ one has*

$$\gamma^{1/2}\|uv\|_{H^m_{*,\gamma}(\mathbb{R}^n_+)} \leq C(\|u\|_{H^m_{*,\gamma}(\mathbb{R}^n_+)}\|v\|_{H^{[\frac{n+1}{2}]+2}_{*,\gamma}(\mathbb{R}^n_+)} + \|u\|_{H^{[\frac{n+1}{2}]+2}_{*,\gamma}(\mathbb{R}^n_+)}\|v\|_{H^m_{*,\gamma}(\mathbb{R}^n_+)}), \qquad \forall \gamma \geq 1. \tag{263}$$

A version of (259), (260) only involving conormal derivatives is given in the following theorem.

**Theorem 43.** *Let $m \geq 1$ be an integer. If $u$ and $v$ are in $H^m_{tan,\gamma}(\mathbb{R}^n_+) \cap L^\infty(\mathbb{R}^n_+)$ then $uv \in H^m_{tan,\gamma}(\mathbb{R}^n_+)$ and there exists a constant $C$ such that*

$$||uv||_{H^m_{tan,\gamma}(\mathbb{R}^n_+)} \leq C(\|u\|_{H^m_{tan,\gamma}(\mathbb{R}^n_+)}\|v\|_{L^\infty(\mathbb{R}^n_+)} + \|u\|_{L^\infty(\mathbb{R}^n_+)}\|v\|_{H^m_{tan,\gamma}(\mathbb{R}^n_+)}), \qquad \forall \gamma \geq 1. \tag{264}$$

*Proof.* The result is proved by induction. For $m = 1$ it is obvious; assuming it is true for $m - 1$ let us take $\alpha$ with $|\alpha| = m$. By Leibniz's rule we have

$$\|Z^\alpha(uv)\|_{L^2(\mathbb{R}^n_+)} \leq C \sum_{\beta \leq \alpha} \|Z^{\alpha-\beta}u\, Z^\beta v\|_{L^2(\mathbb{R}^n_+)} = I_1 + I_2,$$

where we have denoted

$$I_1 = C\left(\|v\,Z^\alpha u\|_{L^2(\mathbb{R}^n_+)} + \|u\,Z^\alpha v\|_{L^2(\mathbb{R}^n_+)}\right), \qquad I_2 = C\sum_{\beta \in K_1(\alpha)} \|Z^{\alpha-\beta}u\, Z^\beta v\|_{L^2(\mathbb{R}^n_+)},$$

$$K_1(\alpha) = \{\beta \leq \alpha,\; 1 \leq |\beta| \leq m-1\}.$$

It is clear that $I_1$ may be readily estimated by the right-hand side of (264). As for $I_2$, from the Hölder's inequality we get

$$\|Z^{\alpha-\beta}u\,Z^\beta v\|_{L^2(\mathbb{R}^n_+)} \leq C \|Z^{\alpha-\beta}u\|_{L^{2m/|\alpha-\beta|}(\mathbb{R}^n_+)} \|Z^\beta v\|_{L^{2m/|\beta|}(\mathbb{R}^n_+)}$$

because $|\alpha-\beta|/2m + |\beta|/2m = 1/2$. Here we apply the interpolation formula (see (B.4) in [24])

$$\|Z^\delta u\|_{L^{2m/|\delta|}(\mathbb{R}^n_+)} \leq C \|u\|_{L^\infty(\mathbb{R}^n_+)}^{1-|\delta|/m} \textstyle\sum_{1\leq|\sigma|\leq m} \|Z^\sigma u\|_{L^2(\mathbb{R}^n_+)}^{|\delta|/m}, \qquad 1 \leq |\delta| \leq m-1,$$

and obtain

$$\begin{aligned} \|Z^{\alpha-\beta}u\,Z^\beta v\|_{L^2(\mathbb{R}^n_+)} &\leq C \|u\|_{L^\infty(\mathbb{R}^n_+)}^{|\beta|/m} \|u\|_{H^m_{tan,\gamma}(\mathbb{R}^n_+)}^{1-|\beta|/m} \|v\|_{L^\infty(\mathbb{R}^n_+)}^{1-|\beta|/m} \|v\|_{H^m_{tan,\gamma}(\mathbb{R}^n_+)}^{|\beta|/m} \\ &\leq C (\|u\|_{L^\infty(\mathbb{R}^n_+)} \|v\|_{H^m_{tan,\gamma}(\mathbb{R}^n_+)} + \|u\|_{H^m_{tan,\gamma}(\mathbb{R}^n_+)} \|v\|_{L^\infty(\mathbb{R}^n_+)}). \end{aligned} \tag{265}$$

Adding over $\alpha,\beta$ completes the proof. $\qquad\square$

From Theorems 41 and 43 we get

**Corollary 44.** *For every integer $m \geq 1$, for all functions $u$ and $v$ in $H^m_{tan,\gamma}(\mathbb{R}^n_+) \cap H^{[\frac{n+1}{2}]+1}_{*,\gamma}(\mathbb{R}^n_+)$ one has*

$$\gamma^{1/2}\|uv\|_{H^m_{tan,\gamma}(\mathbb{R}^n_+)} \leq C(\|u\|_{H^m_{tan,\gamma}(\mathbb{R}^n_+)}\|v\|_{H^{[\frac{n+1}{2}]+1}_{*,\gamma}(\mathbb{R}^n_+)} + \|u\|_{H^{[\frac{n+1}{2}]+1}_{*,\gamma}(\mathbb{R}^n_+)}\|v\|_{H^m_{tan,\gamma}(\mathbb{R}^n_+)}), \quad \forall \gamma \geq 1. \tag{266}$$

Instead of Theorem 43 may be convenient the following one.

**Theorem 45.** *Let $m \geq 1$ be an integer and $s = \max\left\{m, \left[\frac{n+1}{2}\right]+4\right\}$. If $u \in K^m_{*,\gamma}(\mathbb{R}^n_+)$ and $v \in H^s_{*,\gamma}(\mathbb{R}^n_+)$ then $uv \in H^m_{tan,\gamma}(\mathbb{R}^n_+)$ and there exists a constant $C$ such that*

$$\gamma^{s-(n+1)/2}\|uv\|_{H^m_{tan,\gamma}(\mathbb{R}^n_+)} \leq C\|u\|_{K^m_{*,\gamma}(\mathbb{R}^n_+)}\|v\|_{H^s_{*,\gamma}(\mathbb{R}^n_+)}, \qquad \forall \gamma \geq 1. \tag{267}$$

*If $m = 2$ the same result holds with $s = \left[\frac{n+1}{2}\right]+3$ and $\|u\|_{H^2_{tan,\gamma}(\mathbb{R}^n_+)}$ instead of $\|u\|_{K^2_{*,\gamma}(\mathbb{R}^n_+)}$.*

Finally, we give some other lemmata used in the proof of Section 8.

**Lemma 46.** *Let $\sigma \geq [(n+1)/2]+3$ and let $A$ be a matrix-valued function such that $A \in H^\sigma_{*,\gamma}(\mathbb{R}^n_+)$ and $A = 0$ if $x_1 = 0$. Then, for each regular enough vector-valued function $u$*

$$\|A\partial_1 u\|_{L^2(\mathbb{R}^n_+)} \leq c\|A\|_{H^\sigma_{*,\gamma}(\mathbb{R}^n_+)}\|Z_1 u\|_{L^2(\mathbb{R}^n_+)}. \tag{268}$$

*Proof.* See [24], Lemma B.9. $\qquad\square$

**Lemma 47.** *Let $m \geq 2$. Let $A \in H^m_*(\mathbb{R}^n_+)$ be a matrix-valued function such that $A = 0$ if $x_1 = 0$ and let $M$ be defined by*

$$M(x_1, x') = A(x_1, x')/\sigma(x_1),$$

*so that $A\partial_1 u = MZ_1 u$. Then*

$$\|M\|_{H^{m-2}_{*,\gamma}(\mathbb{R}^n_+)} \leq c\|A\|_{H^m_{*,\gamma}(\mathbb{R}^n_+)}.$$

*Proof.* See [34]. $\qquad\square$

## Appendix B. Some estimates

### B.1. **Commutator estimates.**

**Lemma 48.** *If $s > 1$ and $\alpha$ is a $n$-uple of length $|\alpha| \leq s$, there exists $C > 0$ such that for all $u$ and $a$ in $H^s$ with $\nabla u$ and $\nabla a$ in $L^\infty$*

$$\|\,[\partial^\alpha, a\,\nabla]\,u\,\|_{L^2} \leq C\,(\|\nabla a\|_{L^\infty}\|u\|_{H^s} + \|\nabla u\|_{L^\infty}\|a\|_{H^s}).$$

*If $s > 1$ and $\alpha$ is a $n$-uple of length $|\alpha| \leq s$ there exists $C > 0$ such that for all $u$ in $H^{s-1} \cap L^\infty$ and $a$ in $H^s$ with $\nabla a$ in $L^\infty$*

$$\|\,[\partial^\alpha, a]\,u\,\|_{L^2} \leq C\,(\|\nabla a\|_{L^\infty}\|u\|_{H^{s-1}} + \|u\|_{L^\infty}\|a\|_{H^s}).$$

*If $s > n/2 + 1$ and $\alpha$ is a $n$-uple of length $|\alpha| \le s$, there exists $C > 0$ such that for all $u \in H^{|\alpha|-1}$ and $a \in H^s$*

$$\| [\partial^\alpha, a] u \|_{L^2} \le C \|u\|_{H^{|\alpha|-1}} \|a\|_{H^s}.$$

*Proof.* See e.g. [4, 17]. $\qquad\square$

## B.2. Moser-type calculus inequalities.

**Lemma 49.** *For all $s > 0$ there exists $C > 0$ such that for all $u$ and $v$ in $H^s \cap L^\infty$*

$$\|u\,v\|_{H^s} \le C(\|u\|_{H^s}\|v\|_{L^\infty} + \|u\|_{L^\infty}\|v\|_{H^s}). \tag{269}$$

*Proof.* See e.g. [4]. $\qquad\square$

## Appendix C. Adaptation of the result of [35] to the case with outer boundaries

Unlike the reduced linearized problem from [35] formulated in the whole space domain, problem (60) contains the additional boundary conditions (60f),

$$v_1 = 0 \quad \text{on } \omega_T^+, \qquad \mathcal{H}_2 = \mathcal{H}_3 = 0 \quad \text{on } \omega_T^-, \tag{270}$$

on the outer boundaries $\omega_T^+$ and $\omega_T^-$. In [35], for the reduced linearized problem the basic energy a priori estimate was derived and the existence of solutions was proved by using a "hyperbolic regularization" of the elliptic system (60b). Namely, this regularization was inspired by a corresponding problem in relativistic MHD [39] containing the vacuum electric field $E$ as the additional unknown obeying the vacuum Maxwell equations. Introducing the small parameter of regularization $\varepsilon$ and the new auxiliary unknown $E^\varepsilon$, here we just complete the regularized problem from [35] for the unknown $(\mathcal{U}^\varepsilon, V^\varepsilon, \varphi^\varepsilon)$, with $V^\varepsilon = (\mathcal{H}^\varepsilon, E^\varepsilon)$, by adding the boundary conditions (270) written for $v_1^\varepsilon$, $\mathcal{H}_2^\varepsilon$ and $\mathcal{H}_3^\varepsilon$:

$$\widehat{\mathcal{A}}_0 \partial_t \mathcal{U}^\varepsilon + \sum_{j=1}^3 (\widehat{\mathcal{A}}_j + \mathcal{E}_{1j+1}) \partial_j \mathcal{U}^\varepsilon + \widehat{\mathcal{C}} \mathcal{U}^\varepsilon = \mathcal{F} \quad \text{in } Q_T^+, \tag{271a}$$

$$\varepsilon \partial_t \mathfrak{h}^\varepsilon + \nabla \times \mathfrak{E}^\varepsilon = 0, \qquad \varepsilon \partial_t \mathfrak{e}^\varepsilon - \nabla \times \mathfrak{H}^\varepsilon = 0 \quad \text{in } Q_T^-, \tag{271b}$$

$$\partial_t \varphi^\varepsilon = u_1^\varepsilon - \hat{v}_2 \partial_2 \varphi^\varepsilon - \hat{v}_3 \partial_3 \varphi^\varepsilon + \varphi^\varepsilon \partial_1 \hat{v}_N, \tag{271c}$$

$$q^\varepsilon = \widehat{\mathcal{H}} \cdot \mathcal{H}^\varepsilon - [\partial_1 \hat{q}] \varphi^\varepsilon - \varepsilon \widehat{E} \cdot E^\varepsilon, \tag{271d}$$

$$E_{\tau_2}^\varepsilon = \varepsilon \, \partial_t (\widehat{\mathcal{H}}_3 \varphi^\varepsilon) - \varepsilon \, \partial_2 (\widehat{E}_1 \varphi^\varepsilon), \tag{271e}$$

$$E_{\tau_3}^\varepsilon = -\varepsilon \, \partial_t (\widehat{\mathcal{H}}_2 \varphi^\varepsilon) - \varepsilon \, \partial_3 (\widehat{E}_1 \varphi^\varepsilon) \quad \text{on } \omega_T, \tag{271f}$$

$$v_1^\varepsilon = 0 \quad \text{on } \omega_T^+, \qquad \mathcal{H}_2^\varepsilon = \mathcal{H}_3^\varepsilon = 0 \quad \text{on } \omega_T^-, \tag{271g}$$

$$(\mathcal{U}^\varepsilon, V^\varepsilon, \varphi^\varepsilon) = 0 \quad \text{for } t < 0, \tag{271h}$$

where

$$E^\varepsilon = (E_1^\varepsilon, E_2^\varepsilon, E_3^\varepsilon), \quad \widehat{E} = (\widehat{E}_1, \widehat{E}_2, \widehat{E}_3), \quad \mathfrak{E}^\varepsilon = (E_1^\varepsilon \partial_1 \widehat{\Phi}_1, E_{\tau_2}^\varepsilon, E_{\tau_3}^\varepsilon),$$

$$\mathfrak{e}^\varepsilon = (E_N^\varepsilon, E_2^\varepsilon \partial_1 \widehat{\Phi}_1, E_3^\varepsilon \partial_1 \widehat{\Phi}_1), \quad E_N^\varepsilon = E_1^\varepsilon - E_2^\varepsilon \partial_2 \widehat{\Psi} - E_3^\varepsilon \partial_3 \widehat{\Psi}, \quad E_{\tau_k}^\varepsilon = E_1^\varepsilon \partial_k \widehat{\Psi} + E_k^\varepsilon, \ k = 2, 3,$$

the coefficients $\widehat{E}_j$ are given functions which are chosen in [35] so that the boundary conditions (271c)–(271f) on $\omega_T$ are *maximally nonnegative* if we neglect in them the zero-order terms for $\varphi^\varepsilon$. All the other notations for $\mathcal{H}^\varepsilon$ (e.g., $\mathfrak{H}^\varepsilon$, $\mathfrak{h}^\varepsilon$) are analogous to those for $\mathcal{H}$.

The crucial role in the process of deriving the energy estimate in [35] was played by the secondary symmetrization

$$M_0^\varepsilon \partial_t W^\varepsilon + \sum_{j=1}^3 M_j^\varepsilon \partial_j W^\varepsilon + M_4^\varepsilon W^\varepsilon = 0 \tag{272}$$

of the Maxwell equations (271b), where $W^\varepsilon = (\mathfrak{H}^\varepsilon, \mathfrak{E}^\varepsilon)$,

$$M_0^\varepsilon = \frac{1}{\partial_1 \widehat{\Phi}_1} K \mathfrak{B}_0^\varepsilon K^\mathsf{T} > 0, \quad K = I_2 \otimes \hat\eta, \quad M_j^\varepsilon = \frac{1}{\partial_1 \widehat{\Phi}_1} K \mathfrak{B}_j^\varepsilon K^\mathsf{T} \quad (j = 2, 3),$$

$$M_1^\varepsilon = \frac{1}{\partial_1 \widehat{\Phi}_1} \, K \widetilde{\mathfrak{B}}_1^\varepsilon K^\mathsf{T}, \quad \widetilde{\mathfrak{B}}_1^\varepsilon = \frac{1}{\partial_1 \widehat{\Phi}_1} \Big( \mathfrak{B}_1^\varepsilon - \sum_{k=2}^{3} \mathfrak{B}_k^\varepsilon \partial_k \widehat{\Psi} \Big),$$

$$M_4^\varepsilon = K \left( \mathfrak{B}_0^\varepsilon \partial_t + \widetilde{\mathfrak{B}}_1^\varepsilon \partial_1 + \mathfrak{B}_2^\varepsilon \partial_2 + \mathfrak{B}_3^\varepsilon \partial_3 + \mathfrak{B}_0^\varepsilon B_4 \right) \left( \frac{1}{\partial_1 \widehat{\Phi}_1} \, K^T \right),$$

$$\mathfrak{B}_0^\varepsilon = \begin{pmatrix} 1 & 0 & 0 & 0 & \varepsilon\nu_3 & -\varepsilon\nu_2 \\ 0 & 1 & 0 & -\varepsilon\nu_3 & 0 & \varepsilon\nu_1 \\ 0 & 0 & 1 & \varepsilon\nu_2 & -\varepsilon\nu_1 & 0 \\ 0 & -\varepsilon\nu_3 & \varepsilon\nu_2 & 1 & 0 & 0 \\ \varepsilon\nu_3 & 0 & -\varepsilon\nu_1 & 0 & 1 & 0 \\ -\varepsilon\nu_2 & \varepsilon\nu_1 & 0 & 0 & 0 & 1 \end{pmatrix}, \ \mathfrak{B}_1^\varepsilon = \begin{pmatrix} \nu_1 & \nu_2 & \nu_3 & 0 & 0 & 0 \\ \nu_2 & -\nu_1 & 0 & 0 & 0 & -\varepsilon^{-1} \\ \nu_3 & 0 & -\nu_1 & 0 & \varepsilon^{-1} & 0 \\ 0 & 0 & 0 & \nu_1 & \nu_2 & \nu_3 \\ 0 & 0 & \varepsilon^{-1} & \nu_2 & -\nu_1 & 0 \\ 0 & -\varepsilon^{-1} & 0 & \nu_3 & 0 & -\nu_1 \end{pmatrix},$$

$$\mathfrak{B}_2^\varepsilon = \begin{pmatrix} -\nu_2 & \nu_1 & 0 & 0 & 0 & \varepsilon^{-1} \\ \nu_1 & \nu_2 & \nu_3 & 0 & 0 & 0 \\ 0 & \nu_3 & -\nu_2 & -\varepsilon^{-1} & 0 & 0 \\ 0 & 0 & -\varepsilon^{-1} & -\nu_2 & \nu_1 & 0 \\ 0 & 0 & 0 & \nu_1 & \nu_2 & \nu_3 \\ \varepsilon^{-1} & 0 & 0 & 0 & \nu_3 & -\nu_2 \end{pmatrix}, \quad \mathfrak{B}_3^\varepsilon = \begin{pmatrix} -\nu_3 & 0 & \nu_1 & 0 & -\varepsilon^{-1} & 0 \\ 0 & -\nu_3 & \nu_2 & \varepsilon^{-1} & 0 & 0 \\ \nu_1 & \nu_2 & \nu_3 & 0 & 0 & 0 \\ 0 & \varepsilon^{-1} & 0 & -\nu_3 & 0 & \nu_1 \\ -\varepsilon^{-1} & 0 & 0 & 0 & -\nu_3 & \nu_2 \\ 0 & 0 & 0 & \nu_1 & \nu_2 & \nu_3 \end{pmatrix},$$

$I_2$ is the unit matrix of order 2, the matrix $\hat{\eta}$ is defined in (55), and $\nu_i(t, x)$ are functions chosen in appropriate way (see below). It was proved in [35] that systems (271b) and (272) are equivalent, provided that the hyperbolicity condition

$$\varepsilon |\vec{\nu}| < 1 \tag{273}$$

for system (272) is satisfied, where the vector-function $\vec{\nu} = (\nu_1, \nu_2, \nu_3)$. Clearly, inequality (273) holds for any given $\vec{\nu}$ and small $\varepsilon$.

The choice of $\nu_i$ in [35] was the following:

$$\nu_1 = \hat{v}_2 \partial_2 \hat{\varphi} + \hat{v}_3 \partial_3 \hat{\varphi}, \quad \nu_k = \hat{v}_k, \quad k = 2, 3.$$

However, it was important that $\nu_i$ have this form only on the boundary $\omega_T$, i.e. at $x_1 = 0$. Therefore, here we may modify the above choice as follows:

$$\nu_1 = \chi \left( \hat{v}_2 \partial_2 \hat{\varphi} + \hat{v}_3 \partial_3 \hat{\varphi} \right), \quad \nu_k = \chi \, \hat{v}_k, \quad k = 2, 3, \tag{274}$$

where the cut-off function $\chi(x_1) \in C^\infty(-1, 0)$ is such that $\chi(0) = 1$ and $\chi(-1) = 0$. Then, for our present case with the outer boundaries $\omega^\pm = \mathbb{R} \times \Gamma_\pm$ the boundary integral for the inner boundary $\omega = \mathbb{R} \times \Gamma$ appearing in the energy identity (in $L^2$) for the symmetric systems (271a) and (272) stays the same as in [35]:

$$\int_\omega \mathcal{A}^\varepsilon \, dx' dt,$$

where

$$\mathcal{A}^\varepsilon = -\frac{1}{2} (\mathcal{E}_{12} \mathcal{U}_\gamma^\varepsilon, \mathcal{U}_\gamma^\varepsilon)|_\omega + \frac{1}{2} (M_1^\varepsilon W_\gamma^\varepsilon, W_\gamma^\varepsilon)|_\omega, \qquad \mathcal{U}_\gamma^\varepsilon = e^{-\gamma t} \mathcal{U}^\varepsilon, \ \text{etc.},$$

and thanks to the choice in (274) (see [35])

$$\mathcal{A}^\varepsilon = -q^\varepsilon u_1^\varepsilon + \varepsilon^{-1} (\mathfrak{H}_3^\varepsilon \mathfrak{E}_2^\varepsilon - \mathfrak{H}_2^\varepsilon \mathfrak{E}_3^\varepsilon) + (\hat{v}_2 \mathfrak{H}_2^\varepsilon + \hat{v}_3 \mathfrak{H}_3^\varepsilon) \mathcal{H}_N^\varepsilon + (\hat{v}_2 \mathfrak{E}_2^\varepsilon + \hat{v}_3 \mathfrak{E}_3^\varepsilon) E_N^\varepsilon, \quad \text{on } \omega.$$

But in our case the energy identity contains also the boundary integrals

$$\mathfrak{I}^+ = \frac{1}{2} \int_{\omega^+} (\mathcal{E}_{12} \mathcal{U}_\gamma^\varepsilon, \mathcal{U}_\gamma^\varepsilon)|_{\omega^+} \, dx' dt = \int_{\omega^+} q^\varepsilon v_1^\varepsilon|_{\omega^+} \, dx' dt$$

and

$$\mathfrak{I}^- = -\frac{1}{2} \int_{\omega^-} (M_1^\varepsilon W_\gamma^\varepsilon, W_\gamma^\varepsilon)|_{\omega^-} \, dx' dt$$

for the top and bottom boundaries, where it is worth noting that thanks to the choice in (274) we have $\nu_i|_{x_1=-1} = 0$ and so one can check that

$$M_1^\varepsilon|_{\omega^-} = \mathcal{B}_1^\varepsilon|_{\vec{\nu}=0} = B_1^\varepsilon = \varepsilon^{-1} \begin{pmatrix} 0 & 0 & 0 & 0 & 0 & 0 \\ 0 & 0 & 0 & 0 & 0 & -1 \\ 0 & 0 & 0 & 0 & 1 & 0 \\ 0 & 0 & 0 & 0 & 0 & 0 \\ 0 & 0 & 1 & 0 & 0 & 0 \\ 0 & -1 & 0 & 0 & 0 & 0 \end{pmatrix}.$$

Then,

$$\mathfrak{I}^- = \varepsilon^{-1} \int_{\omega^-} (E_3^\varepsilon \mathcal{H}_2^\varepsilon - E_2^\varepsilon \mathcal{H}_3^\varepsilon)|_{\omega^-} \, dx' dt.$$

Thanks to the boundary conditions (271g) the both boundary integrals above vanish: $\mathfrak{I}^\pm = 0$.

The remaining arguments are the same as in [35] and we just refer the reader to [35] for more details. For our case (with the added outer boundaries), for the existence of solutions of problem (271) it is only important to note that the number of the boundary conditions in (271g) is in the agreement with the number of incoming characteristics for the boundaries $\Gamma_\pm$ (this is easily checked by calculating the eigenvalues of the matrices $\mathcal{E}_{12}$ and $B_1^\varepsilon$). Then, we again refer to [35] for the energy a priori estimate for the regularized problem (which is the same as that for problem (271)), the proof of the existence of solutions for it and the passage to the limit as $\varepsilon \to 0$.

*E-mail address*: paolo.secchi@ing.unibs.it

*E-mail address*: trakhin@math.nsc.ru